\documentclass[11pt,reqno]{amsart}

\usepackage[T1]{fontenc}
\usepackage{algorithm} 
\usepackage{algpseudocode} 
\usepackage{float}
\usepackage{enumerate}

\usepackage{subcaption}
\usepackage{caption} 
\usepackage{graphicx} 
\newtheorem{theorem}{Theorem}[section]
\newtheorem{lemma}[theorem]{Lemma}

\newtheorem{proposition}[theorem]{Proposition}

\newtheorem{assumption}{\bf Assumption}

\theoremstyle{remark}
\newtheorem{remark}{Remark}[section]

\newtheorem{example}{\bf Example}

\usepackage[svgnames]{xcolor}
\usepackage{color}
\usepackage{amssymb, amsmath, amsfonts, latexsym}

\usepackage{amsthm}

\usepackage{enumerate}
\usepackage{bm}
\usepackage{bbm} 
\usepackage{enumitem} 
\usepackage{hyperref}
\hypersetup{colorlinks=true, citecolor = teal, linkcolor=blue!80!black}

\usepackage[initials,alphabetic]{amsrefs}

\newcommand{\dd}{\mathrm{d}}
\newcommand{\ii}{\mathrm{i}}
\newcommand{\supp}{\operatorname{supp}}
\newcommand{\argmin}{\operatorname{argmin}}
\newcommand{\blambda}{\boldsymbol{\lambda}}

\newcommand{\by}{\boldsymbol{y}}

\newcommand{\bz}{\boldsymbol{z}}
\newcommand{\pv}{\operatorname{p.v.}}

\newcommand{\PP}{\mathbb P}

\newcommand{\R}{\mathbb R}
\newcommand{\C}{\mathbb C}

\newcommand{\dist}{\operatorname{dist}}

\newcommand{\sgn}{{\rm sgn}}

\allowdisplaybreaks[4]

\newcommand{\8}{\infty}

\newcommand{\cx}{\mathcal{C}}

\newcommand{\mx}{\mathcal{M}}

\newcommand{\rz}{\mathbb{R}}
\newcommand{\cz}{\mathbb{C}}

\newcommand{\pz}{\mathbb{P}}

\newcommand{\de}{\delta}
\newcommand{\si}{\sigma}
\newcommand{\ga}{\gamma}
\newcommand{\la}{\lambda}
\newcommand{\bt}{\beta}

\newcommand{\eps}{\varepsilon}

\numberwithin{equation}{section}
\numberwithin{figure}{section}

\begin{document}

\title[Lower tail large deviations for $\beta$-ensembles]{Lower tail large deviations of the extreme eigenvalue of $\beta$-ensembles}

\author{Siyu Wang}  
\address{Academy of Mathematics and Systems Science, Chinese Academy of Sciences, 100190 Beijing, China.}
\email{wang\_siyu@mail.bnu.edu.cn}
\author{Qiang Zeng}
\address{Academy of Mathematics and Systems Science, Chinese Academy of Sciences, 100190 Beijing, China.}
\email{qzeng@amss.ac.cn}


\date{November 29, 2026}


\begin{abstract}
We consider large deviations of $\bt$-ensembles and Coulomb gases. We start with a general large deviation principle for the empirical spectral measure when the external potential is defined on a closed subset of the whole space. This result covers in particular some `corner cases' overlooked in the mathematical literature and allows us to derive precise rate functions for the lower tail large deviation probabilities of the largest eigenvalue of $\beta$-ensembles and that of the spectral radius of Coulomb gases. Then we provide a unified mathematical treatment for the lower tail large deviations of the largest eigenvalue of the Hermite, Laguerre, and Jacobi ensembles that have appeared in the physics literature.  Although our focus is on the lower tail probabilities, we also prove complementary results for the upper tail. Our approach relies primarily on logarithmic potential theory. 
\end{abstract}

\maketitle

\tableofcontents

\section{Introduction}
Large deviation theory of random matrices has been extensively investigated in the past few decades in both mathematics and physics. Previous results are mainly concerned with the large deviation principles (LDPs) of the empirical spectral measure and the extreme eigenvalue for various ensembles; for an incomplete list, see \cites{BAG97, GZ00,Har12,BC14,BGH22,Au25} for the former and \cites{BDG01, AGZ10,BG13,Au16, GH20,AGH21,HM24,AB26,CRG26} for the latter. The theory has found tremendous applications in different contexts, including free probability and spin glasses, among others; see the survey \cite{Gui23} for an overview and related literature on Hermitian matrices.  On the other hand, LDP for the extreme eigenvalue usually gives exact speed and rate function only for the upper tail and provides limited information on the lower tail. Indeed, although the speed of lower tail large deviations coincides with that of the LDP of empirical measure, the rate function cannot be read off directly from those LDP results. Meanwhile, lower tail large deviation probabilities of extreme eigenvalue of random matrices have been considered in physics literature; see, e.g., \cites{DM06,VMB07, CMV16, Lac18, Lac19}.

In this paper, we introduce a unified framework to study the lower tail large deviations of the extreme eigenvalue of  $\beta$-ensembles and Coulomb gases. Let $D\subset\mathbb C$ be a closed set. For each $n\ge1$, we consider the probability measure
\begin{equation}\label{pdf}
\mathbb P(\dd\blambda)
=
\frac{1}{Z_n}
\prod_{1\le i<j\le n}|\lambda_i-\lambda_j|^\beta
\exp\Bigl\{-\beta n\sum_{i=1}^n V_n(\lambda_i)\Bigr\} 
\dd\blambda,
\end{equation}
where $\beta=\beta_n>0$ may depend on $n$, $V_n$ is an external potential on $D$ and $Z_n$ is the normalizing constant. 
Here $\dd\boldsymbol{\lambda}=\dd\lambda_1\cdots \dd\lambda_n$ denotes the product measure on $D^n$. More precisely, if $D\subset\mathbb R$, then each $\dd\lambda_i$ is the one-dimensional Lebesgue measure; if $D\not\subset\mathbb R$, then $\dd\lambda_i=\dd(\Re \lambda_i) \dd(\Im \lambda_i).$ We refer the reader to \cites{AGZ10,Ser26} for more background and related topics. Our primary goal is to derive the precise rate function (or explicit exponential rate of leading order) of the lower tail large deviation probability of extreme eigenvalue.

This model includes several classical ensembles. When $D\subset\mathbb R$, \eqref{pdf} gives the usual one-dimensional $\beta$-ensemble \cite{For10}. In particular, if $D=\mathbb R$ and $V_n(x)=x^2/2$, the density in \eqref{pdf} is that of the normalized $\beta$-Hermite ensemble (a.k.a.~Gaussian $\bt$-ensemble, or G$\bt$E). 
For $\beta\in\{1,2,4\}$, it is the joint eigenvalue density  of $H_{n}  =  \frac{1}{2\sqrt{n}}(G+G^*),$ where $G=G_{n}$  is an $n\times n$ Ginibre random matrix whose entries are independent standard Gaussian variables with values in  $\mathbb R$, $\mathbb C$, or $\mathbb H$, and $G_n^*$ denotes the transpose, complex conjugate transpose, or quaternion conjugate transpose, respectively. Here $\mathbb H$ is the set of quaternions. The matrix $H_n$ is said to belong to the Gaussian orthogonal ensemble (GOE), Gaussian unitary ensemble (GUE), or Gaussian symplectic ensemble (GSE), respectively.  The empirical spectral measure $L_n=\frac1n\sum_{i=1}^n\delta_{\lambda_i}$ of the normalized $\beta$-Hermite ensemble converges weakly a.s.~to the semicircle law \cites{Wig55, AGZ10} with density
\begin{align}\label{eq:semic}
    \rho^H(x)=\pi^{-1}\sqrt{2-x^2}\mathbf 1_{[-\sqrt2,\sqrt2]}(x).
\end{align}
When $D=\mathbb R_+=[0,\8)$ and
\begin{align}\label{def:vnl}
    V_n^L(x)
    =
    \frac{p}{2n}x
    -
    \left(
        \frac{p-n+1}{2n}
        -
        \frac{1}{\beta n}
    \right)\log x,
    \qquad x>0,\quad p\ge n,
\end{align}
 the density in \eqref{pdf} defines the $\beta$-Laguerre ensemble, for which we assume
\begin{equation}\label{eq:Lagpar}
\frac np \to \gamma\in[0,1].
\end{equation}
For $\beta\in\{1,2,4\}$, this is the eigenvalue density of the scaled Wishart (or sample covariance) matrix $  W_n=p^{-1}XX^*,$ where $X=X_{n\times p}$ is an $n\times p$ standard Gaussian matrix over $\mathbb R$, $\mathbb C$, and $\mathbb H$, respectively. 
 For $\gamma>0$,  the empirical measure of the scaled $\beta$-Laguerre ensemble converges weakly a.s.~to the Marchenko--Pastur law \cites{MP67,DE06,TT21} with density 
\begin{align}\label{eq:hgam}
\rho_{\gamma}^L(x)
=
\frac{1}{2\pi\gamma x}
\sqrt{(x-\gamma_{-})(\gamma_{+}-x)}\mathbf 1_{[\gamma_-,\gamma_+]}(x),
\quad
\gamma_{\pm}=(1\pm\sqrt{\gamma})^2.
\end{align}
When $D=[0,1]$ and
\begin{align} \label{def:vnj}
V_n^J(x)
    =
    -
    \left(
        \frac{p_1-n+1}{2n}
        -
        \frac{1}{\beta n}
    \right)\log x
    -
    \left(
        \frac{p_2-n+1}{2n}
        -
        \frac{1}{\beta n}
    \right)\log(1-x),
    \qquad 0<x<1,
\end{align}
the density in \eqref{pdf} defines the $\beta$-Jacobi ensemble with $p_1,p_2\geq n$, for which we assume
\begin{equation}\label{eq:Jacpar}
\frac{n}{p_1} \to \gamma\in[0,1],
\qquad
\frac{p_1}{p_2} \to \sigma\in[0,\infty),
\qquad \gamma\sigma\leq1.
\end{equation}
For $\beta\in\{1,2,4\}$, it is the eigenvalue density of the MANOVA matrix \cites{Con63,Mui82,For06}
$J_n = XX^*(XX^*+YY^*)^{-1},$
where $X=X_{n\times p_1}$ and $Y=Y_{n\times p_2}$ are independent standard Gaussian matrices over $\mathbb R$, $\mathbb C$, or $\mathbb H$, respectively. For  $\gamma\sigma>0$, the empirical measure of the $\beta$-Jacobi ensemble converges weakly a.s.~to the Wachter law  \cites{Wa78,DP12} with density 
\begin{align}\label{eq:hgamsig}
\rho^J_{\gamma,\sigma}(x)
=
\frac{1+\sigma}{2\pi\sigma\gamma}
\frac{\sqrt{(x-u_{-})(u_{+}-x)}}{x(1-x)}\mathbf 1_{[u_{-},u_{+}]}(x),
\quad 
u_{\pm}
=
\frac{\sigma}{1+\sigma}
\left(
\sqrt{1-\frac{\sigma\gamma}{1+\sigma}}
\pm
\sqrt{\frac{\gamma}{1+\sigma}}
\right)^2.
\end{align}
Unless otherwise stated, we assume \eqref{eq:Lagpar} for the $\beta$-Laguerre ensemble and \eqref{eq:Jacpar} for the $\beta$-Jacobi ensemble throughout.
For related results on $\beta$-ensembles in more general parameter regimes, including dimension dependent $\beta=\beta_n$ and nonstandard scalings of the ensemble parameters, see \cites{Jia13,JL15, Nag14,TT21,Ma23J, Ma23L, GN25, MMW25} for an incomplete list.
When $D=\mathbb C$, \eqref{pdf} becomes a Coulomb gas. In the case $\beta=2$, it includes the complex Ginibre and complex elliptic Ginibre ensembles; see \cites{PH98,Ser26}.

Our main ideas in this paper can be summarized as follows:
\begin{enumerate}[label=\textup{(\Roman*)}]
    \item\label{mainid1} The high dimensional limit of the lower tail large deviation probability of the extreme eigenvalue can be justified using the LDP of empirical measure with constrained external potential on proper closed subsets of $\rz/\cz$.
    \item\label{mainid2} The rate function of lower tail large deviations can be computed (explicitly in scenarios of our interest) via the logarithmic potential theory.
\end{enumerate}
A natural approach to derive the lower tail large deviation rate function is to minimize the rate function of the empirical measure LDP over the constrained region or interval. Our perspective here is slightly different in that we do not need a full LDP for the empirical measure. Instead, we require only the local LDP lower bound for two measures, which enables us to consider the lower tail probability for some irregular regions where a full LDP does not hold; see Remark \ref{re:assub}. To this end, for \ref{mainid1} we derive a general LDP for the empirical measure that distinguishes the conditions required in the upper and lower bounds. Moreover, it covers some classical examples with boundary singularities that are not treated before. For example, the squared real Wishart matrix ($\bt=1$ and $p=n$ in \eqref{def:vnl}) has been long overlooked before and is now justified; see the discussion in Section \ref{se:ldpln}. For \ref{mainid2}, we provide a unified approach to the results for classical $\bt$-ensembles obtained in the physics literature. 
Let us elaborate these aspects in the following.


Let $\mathcal M_1(D)$ denote the space of Borel probability measures on $D$ endowed with the topology of weak convergence. For $\mu\in\mathcal M_1(D)$, define the noncommutative (or free) entropy
\begin{align}\label{def:sigmu}
\Sigma(\mu)= \begin{cases}\iint \log |x-y|\mu( \dd x) \mu(\dd y), & \text { if } \int \log (|x|+1) \mu( \dd x)<\infty, \\ -\infty, & \text { otherwise. }
\end{cases}
\end{align}
The quantity $-\Sigma(\mu)$ is also known as the logarithmic energy of $\mu$. 
For a lower semicontinuous function $V:D\to\mathbb R\cup\{+\infty\}$, define
\begin{align}\label{defiv}
I_V(\mu)= \begin{cases}\int V \dd \mu-\frac{1}{2} \Sigma(\mu)-C_V, & \text { if } \int V \dd \mu<\infty, \\ \infty, & \text { otherwise. }\end{cases}
\end{align}
where $C_V=\inf_{\mu \in \mathcal M_1(D)}\left\{\int V \dd \mu-\frac{1}{2} \Sigma(\mu)\right\} \in \rz$ under some assumptions to be explained below. 

When the external potential does not depend on $n$, Ben Arous and Guionnet proved that the empirical measure $L_n=n^{-1}\sum_{i=1}^n\delta_{\lambda_i}$ satisfies an LDP with speed $n^2$ and good rate function $I_V$ for GOE in the seminal work \cite{BAG97}. Similar results were obtained by Ben Arous--Zeitouni \cite{BAZ98} and Petz--Hiai \cite{PH98} for the Ginibre ensembles. LDPs with more general external potentials have been considered in many other settings; see, e.g., \cites{HP98, ES08,Har12, HM13, CGZ14,AGZ10, Ma23J, Ma23L}. These results usually require $\bt n\gg \log n$ and focus on special ensembles. In the `very high' temperature regime $\bt n\to c\in [0,\8)$, there are also works when $D$ is the whole space; see \cites{GZ19,LW20,GM22}. The regime $\bt n\to\8$ was covered by \cite{DLR20} for lower semicontinuous $V$ independent of $n$ that may take value $+\8$. The work \cite{GZ19} also contains results for potentials depending on $n$ defined on a Polish space in the regime $\bt n\to\8$. However, it requires a uniform lower bound for the potentials (termed `stable sequence' there) w.r.t.~a chosen reference measure, which may be unavailable in practice.  

In this paper, our primary interest on LDP of the empirical measure is to extend the above well known results to the setting where the external potentials $V_n$ are defined on a proper subset of $\rz$ or $\cz$ so that it provides a rigorous way to compute the limit of lower tail probability of the extreme eigenvalue. To this end, we give assumptions that are straightforward (but may be technical) to verify for the examples of interest.
Under suitable convergence and growth assumptions on $V_n$, we prove that the empirical measure satisfies an LDP with speed $\beta n^2$ and rate function $I_V$, where the potential $V$ is obtained as the limit of $V_n$. In particular, this allows different approximations of the potential $V$ with the rate function unchanged.
Although allowing $\bt n\to\8$ is a byproduct of our method, it is crucial for our development that $V_n$ depends on $n$ and may be unbounded below, and that $D$ is a proper subset of $\cz$. As we shall show below, the general setting could cause some subtle issues for the LDP. 

\begin{theorem}\label{ldpln}
Let $\boldsymbol\lambda=(\lambda_1,\dots,\lambda_n)\in D^n$ have density \eqref{pdf}.
Suppose that Assumption \ref{ass:ldp} holds and that every $\mu\in \mathcal M_1(D)$  with $I_V(\mu)<\8$  satisfies Assumption \ref{ass:d}. Then the empirical measure $L_n=\frac1n\sum_{i=1}^n\delta_{\lambda_i}$ satisfies an LDP on $\mathcal M_1(D)$ with speed $\beta n^2$ and good rate function $I_V$ defined in \eqref{defiv}.
\end{theorem}

The assumptions \ref{ass:ldp} and \ref{ass:d} are given in Section \ref{sec:ex} where they are also shown to be general enough to include and extend many previous results. However, we stress that Theorem \ref{ldpln} is tailored to suit our following applications, and generality is not our focus here compared to some works in this direction.

Let us explain how we apply this result. For the lower tail probability of extreme eigenvalue, the large deviation rate is given by the minimum of the rate function $I_V$ over the constraint set determined by the lower tail event. A natural approach is to add hard walls to the potential. Nevertheless, it seems that this idea has seldom been rigorously addressed in the literature, as existing works on large deviations of the empirical measure mainly focused on the case that $V$ is real valued on the whole space; see however \cites{BG13,GZ19,DLR20} which are close to our setting but still insufficient for our direct applications. Our current framework allows $V_n$ to be unbounded below and to take values in $\rz\cup \{+\8\}$ through a decomposition scheme, which is the main novelty of our approach; see Remark \ref{r:a3unit} for a comparison with \cite{GZ19}. 
As in the classic setting,  Theorem \ref{ldpln} is also a building block for the extreme eigenvalue LDP in Proposition \ref{pr:maxldp}. Moreover, these LDPs are cornerstones for some related questions; see, e.g., \cite{WZschatten}.

After Theorem \ref{ldpln} is established, seeking the lower tail rate function amounts to solving the minimization problem with specific potentials that jump to infinity at finite points. 
Besides, we emphasize that the rate function $I_V$ of the empirical measure LDP is the only input to evaluate the lower tail rate function; the precise eigenvalue joint density is hidden in the picture. This yields unification of different densities with the same $I_V$.

We start to investigate the lower tail large deviation probabilities. In most scenarios, we actually consider the slightly more general two sided constrained problem, i.e., 
$$\pz(\zeta_1\le \la_i\le \zeta_2, i=1,\dots,n)$$ 
which trivially reduces to the lower tail probability by setting $\zeta_1=-\8$. For $\beta$-ensembles, Proposition \ref{pr:lowgap} provides a general variational characterization of this problem. For the three classical $\beta$-ensembles, explicit lower tail formulas 
have appeared in the physics literature \cites{DM06,DM08,VMB07,RKP12}. Here we consider the rescaled $\beta$-Jacobi ensemble that rigorously unifies and extends the previous physics results. 
Following \cites{Nag14,Ma23J}, we define
\begin{align}\label{defvarphi}
y=\varphi_n(x)=\frac{(p_1+p_2)x-p_1}{\sqrt{p_1 n}},
\quad
x=\varphi_n^{-1}(y)=\frac{p_1+\sqrt{p_1 n}y}{p_1+p_2}.	
\end{align}
Let $D_{\gamma,\sigma}=[-\gamma^{-1/2},(\sigma\sqrt\gamma)^{-1}]$
with $1/0=+\infty$. Up to additive constants,
$V_n^J\circ\varphi_n^{-1}$ converges locally uniformly on the interior of
$D_{\gamma,\sigma}$ to
\begin{equation}\label{eq:jpot}
V_{\gamma,\sigma}^J(x)=
\begin{cases}
-\frac{1-\gamma}{2\gamma}\log(1+\sqrt\gamma x)
-\frac{1-\sigma\gamma}{2\sigma\gamma}
 \log(1-\sigma\sqrt\gamma x),
&\gamma\sigma>0,\\
-\frac{1-\gamma}{2\gamma}\log(1+\sqrt\gamma x)
+\frac{x}{2\sqrt\gamma},
&\sigma=0<\gamma,\\
\frac{1+\sigma}{4}x^2,
&\gamma=0.
\end{cases}
\end{equation}
Define the values at finite endpoints by interior limits and set $V_{\gamma,\sigma}^J=+\infty$ outside $D_{\gamma,\sigma}$.
We are ready to state the main result.

\begin{theorem}
\label{thm:jaclow}
Let $(\lambda_1,\dots,\lambda_n)$ be the eigenvalues of $\beta$-Jacobi ensemble with potential \eqref{def:vnj} and assume that $\beta n\to\infty$. 
Let $-\gamma^{-1/2}\leq\zeta_1<\zeta_2 \leq(\sigma\sqrt\gamma)^{-1}$ with the above limiting interpretations when $\gamma=0$ or $\sigma=0$.  Then
\begin{equation}\label{eq:jwall}
 \lim_{n\to\infty}\frac1{\beta n^2}
 \log\PP\left(
   \zeta_1\leq
   \frac{(p_1+p_2)\lambda_i-p_1}{\sqrt{p_1n}}
   \leq\zeta_2,
    1\leq i\leq n
 \right)
 =-I_{V_{\gamma,\sigma}^J}(\mu_{\zeta_1,\zeta_2}^J),
\end{equation}
where
\begin{align}\label{eq:con-J}
  \mu_{\zeta_1,\zeta_2}^J
 =\argmin_{\mu\in\mathcal M_1([\zeta_1,\zeta_2])}
 I_{V_{\gamma,\sigma}^J}(\mu).
\end{align}
The support $[s,t]$ and density of $\mu_{\zeta_1,\zeta_2}^J$ are given in
Proposition \ref{p:jacmin}.  Moreover,
\begin{equation}\label{eq:Jac_t}
 I_{V_{\gamma,\sigma}^J}(\mu_{\zeta_1,\zeta_2}^J)
 =\begin{cases}
 I_{V^H}\left(
   \mu_{\sqrt{(1+\sigma)/2} \zeta_1, 
        \sqrt{(1+\sigma)/2} \zeta_2}^H\right),
 &\gamma=0,\quad\sigma\geq0,\\
 I_{V_\gamma^L}(\mu_{\zeta_1,\zeta_2}^L),
 &\sigma=0,\quad0<\gamma\leq1,\\
 \Phi_{\gamma,\sigma}^J(s,t),
 &\sigma\gamma>0.
 \end{cases}
\end{equation}
Here the rate functions in the first two cases reduce to those of Hermite and Laguerre ensembles from Examples \ref{ex:her2} and \ref{ex:Lag2}, respectively. The explicit expression for
$\Phi_{\gamma,\sigma}^J$ is given in
\eqref{eq:jacgen} below.

In particular, when $\gamma=1$ and $0<\sigma\leq1$, 
\begin{align*}
 I_{V_{1,\sigma}^J}(\mu_{\zeta_1,\zeta_2}^J)
 ={}&-\frac12\log\frac{(1+\sigma)(t-s)}{4}
 -\frac{(1-\sigma)(1+3\sigma)}{4\sigma^2}
 \log\frac{\sqrt{1+\sigma}
 \bigl(\sqrt{1-\sigma s}+\sqrt{1-\sigma t}\bigr)}{2}\\
 &+\frac{(1-\sigma)^2}{4\sigma^2}
 \log\frac{2\sqrt{(1-\sigma s)(1-\sigma t)}}
 {(1-\sigma)\sqrt{1+\sigma}
 \bigl(\sqrt{1-\sigma s}+\sqrt{1-\sigma t}\bigr)},
\end{align*}
where the singular terms at $\sigma=1$ are defined by continuity.
\end{theorem}

\begin{figure}[!ht]
 \centering
 \begin{subfigure}[t]{0.32\textwidth}
 \centering
 \includegraphics[width=\linewidth]{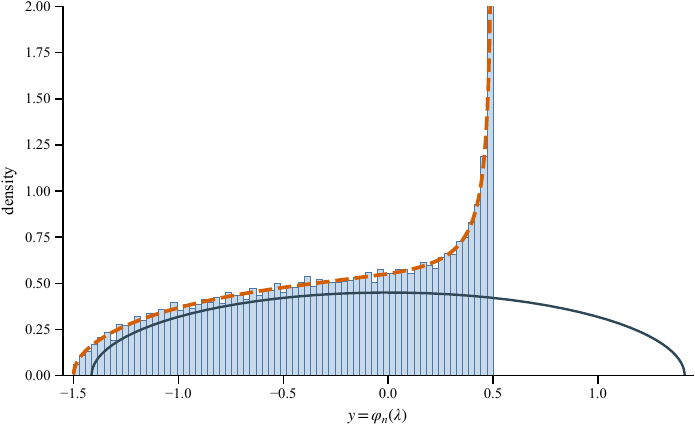}
 \caption{$p_1=20000$, $p_2=20000$.}
 \end{subfigure}\hfill
 \begin{subfigure}[t]{0.32\textwidth}
 \centering
 \includegraphics[width=\linewidth]{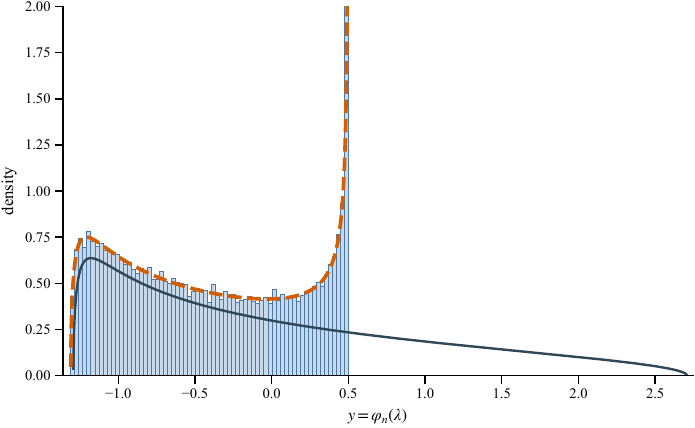}
 \caption{$p_1=200$, $p_2=20000$.}
 \end{subfigure}\hfill
 \begin{subfigure}[t]{0.32\textwidth}
 \centering
 \includegraphics[width=\linewidth]{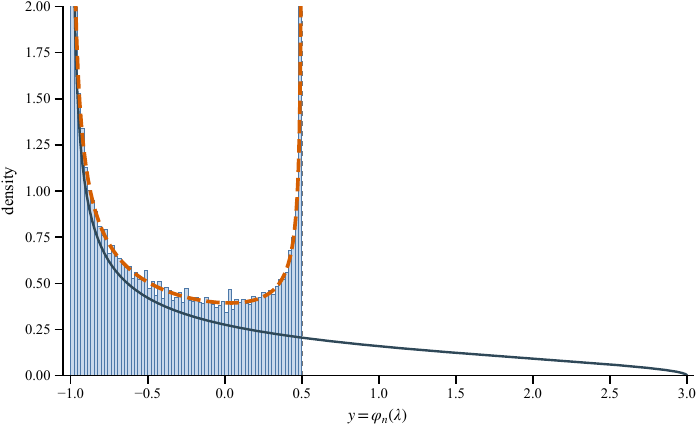}
 \caption{$p_1=100$, $p_2=20000$.}
 \end{subfigure}
 \begin{subfigure}[t]{0.32\textwidth}
 \centering
 \includegraphics[width=\linewidth]{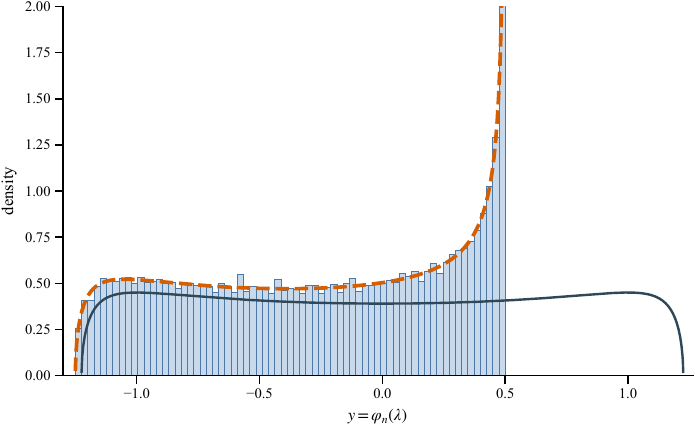}
 \caption{$p_1=200$, $p_2=200$.}
 \end{subfigure}\hfill
 \begin{subfigure}[t]{0.32\textwidth}
 \centering
 \includegraphics[width=\linewidth]{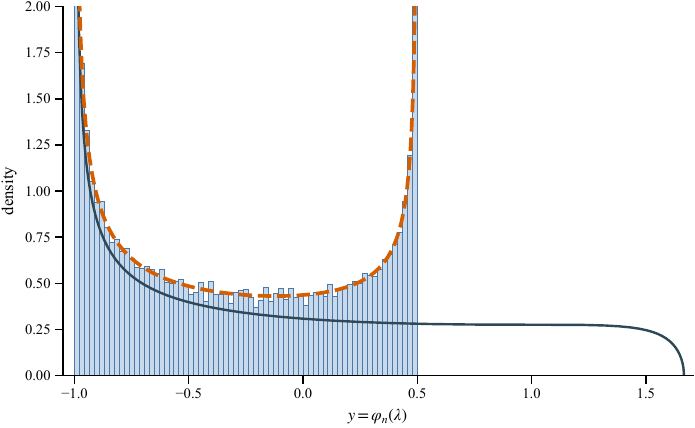}
 \caption{$p_1=100$, $p_2=200$.}
 \end{subfigure}\hfill
 \begin{subfigure}[t]{0.32\textwidth}
 \centering
 \includegraphics[width=\linewidth]{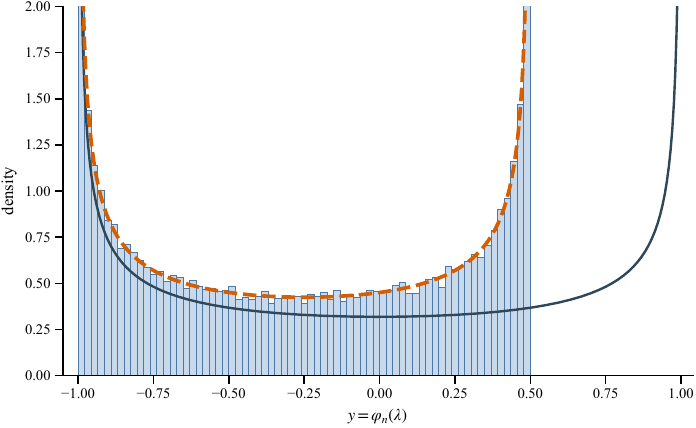}
 \caption{$p_1=100$, $p_2=100$.}
 \end{subfigure}
\caption{Constrained and unconstrained densities for the rescaled $\beta$-Jacobi ensemble with wall $\zeta=0.5$ when $\beta=2$ and $n=100$. The blue histogram is obtained by Metropolis--Hastings sampling \cite{MH70} under the constraint $\max_{1\le i\le n} \frac{(p_1+p_2)\lambda_i-p_1}{\sqrt{p_1 n}}\le 0.5$. The orange dashed curve is density of $\mu_\zeta^J$, and the dark solid curve is the unconstrained density. The illustrated parameters are as follows: (A) $\ga=0,\si=1$; (B) $\ga=0.5,\si=0$; (C) $\ga=1,\si=0$; (D) $\ga=0.5,\si=1$; (E) $\ga=1, \si=0.5$; (F) $\ga=\si=1$.
}

\label{fig:jacobi}
\end{figure}

The various parameter regimes are illustrated in Figure \ref{fig:jacobi}. 

\begin{remark}
Note that \eqref{eq:jwall} can be written equivalently as
$$
 \lim_{n\to\infty}\frac1{\beta n^2}
 \log\PP\left(
 \frac{p_1+\sqrt{p_1n} \zeta_1}{p_1+p_2}
 \leq\lambda_i\leq
 \frac{p_1+\sqrt{p_1n} \zeta_2}{p_1+p_2}, 1\leq i\leq n
 \right)
 =-I_{V_{\gamma,\sigma}^J}(\mu_{\zeta_1,\zeta_2}^J).
$$
For $\gamma\sigma>0$ and fixed $0\leq\xi_1<\xi_2\leq1$, continuity of $(a,b)\mapsto I_{V_{\gamma,\sigma}^J}(\mu_{a,b}^J)$ yields
\begin{equation}\label{e:Jacwal}
 \lim_{n\to\infty}\frac1{\beta n^2}
 \log\PP(\xi_1\leq\lambda_i\leq\xi_2, 1\leq i\leq n)
 =-I_{V_{\gamma,\sigma}^J}(\mu_{\varphi(\xi_1), \varphi(\xi_2)}^J)
 =-\Phi_{\gamma,\sigma}^J(s,t),
\end{equation}
where $\varphi$ is given in \eqref{eq:joremp} and $[s,t]=\supp(\mu_{\varphi(\xi_1), \varphi(\xi_2)}^J)$. 
In the notation of \cite{RKP12}, $N=n$ and $M_i=p_i$.
Their limiting parameters are therefore
$$
\alpha_1=\lim_{n\to\infty}\left(\frac{p_1}{n}-1\right)
 =\frac{1-\gamma}{\gamma},\qquad
\alpha_2=\lim_{n\to\infty}\left(\frac{p_2}{n}-1\right)
 =\frac{1-\gamma\sigma}{\gamma\sigma}.
$$
The action $S$ in \cite{RKP12}*{Eq.~(17)} satisfies
\begin{equation}\label{e:JacRKP}
 \Phi_{\gamma,\sigma}^J(s,t)
 =\frac12\bigl[S(\varphi^{-1}(s), \varphi^{-1}(t))-S(u_-,u_+)\bigr].
\end{equation}
The unconstrained endpoints $u_\pm$ are given in \eqref{eq:hgamsig}.
Thus \eqref{e:JacRKP} identifies \eqref{e:Jacwal} with the rate
formula in \cite{RKP12}*{Section~4.5}. The critical cases
$\gamma=1$ or $\gamma\sigma=1$ are understood by continuity.

When $\gamma\sigma=0$, the original empirical measure $n^{-1}\sum_{i=1}^n\delta_{\lambda_i}$ converges weakly a.s.~to $\delta_{a_\sigma}$, where $a_\sigma=\sigma/(1+\sigma)$
(see \cite{DP12} for fixed $\beta>0$).
The empirical measure LDP in Example \ref{ex:JacempLDP} and \eqref{eq:jwall} give, for every fixed $x\neq a_\sigma$,
$$
 \lim_{n\to\infty}\frac1{\beta n^2}
 \log\PP(\lambda_{\max}\leq x)
 =\begin{cases}
 -\infty,&x<a_\sigma,\\
 0,&x>a_\sigma.
 \end{cases}
$$
The scaling \eqref{defvarphi} resolves this degeneracy.
By \eqref{eq:Jac_t}, the rescaled rates reduce to the Hermite rate when $\gamma=0$ and to the Laguerre rate when $\sigma=0<\gamma$. Since the rescaled external potentials $V_n\circ \varphi_n^{-1}$ approximate the Hermite, Laguerre and Jacobi potentials in \eqref{eq:jpot} (up to scaling) respectively, Theorem \ref{thm:jaclow} unifies the lower tail large deviations for all the three classical ensembles.
\end{remark}

The paper is organized as follows. We present the general theory and first examples in Section \ref{sec:ex}. Theorem \ref{ldpln} is proved in Section \ref{se:ldpln}. Proposition \ref{pr:maxldp} provides an LDP for the extreme eigenvalue in Section \ref{se:maxldp}, whose proof is deferred to Appendix \ref{app:maxldp}.  The assumptions of these LDPs are  discussed thoroughly together with examples in the respective sections. In Section \ref{se:maxlow}, Proposition \ref{pr:lowgap} gives a general principle on confinement probability, and in Section \ref{se:srlow}, we prove the planar analog Proposition \ref{thm:lowspe} for Coulomb gases. Section \ref{sec:ex} contains empirical measure LDP, upper and lower tail large deviations of extreme eigenvalue of the Hermite, Laguerre and Jacobi ensembles. 
In Section \ref{sec:gen}, we give proofs for the examples and Theorem \ref{thm:jaclow}. 

\noindent\textbf{Acknowledgments.} We thank Jiaoyang Huang for helpful discussions and suggestions.

\noindent\textbf{Statement on AI use.} 
The main ideas  \ref{mainid1} and \ref{mainid2} of this work are human generated. When implementing these ideas, we employed large language models (ChatGPT and DeepSeek with versions from Summer 2025 to Summer 2026 during this project) to assist grammar check, figure creation, literature search, and some of the technical proofs, especially evaluation of some involved integrals. The main results on lower tail large deviations had been obtained by May 2026 along with other results including \cite{WZschatten}, and were presented in several conferences and seminars. Then we realized that the LDPs for the empirical measure and the largest eigenvalue in standard references did not cover all classical examples. Subsequently we explored reasonably weak assumptions to include all cases of interest by many-round interactions with ChatGPT 5.5 \& 5.6. Therefore, the most important contributions of LLMs are in Theorem \ref{ldpln} and Proposition \ref{pr:maxldp}. More specifically, Assumption \ref{A3}, the use of Fatou's lemma from \cite{WFatou} for proving LDP upper bound and the use of $U$-statistics for proving LDP lower bound are proposed by AI, which simplify our original arguments and are adopted. The pinning argument for Proposition \ref{pr:maxldp} was suggested by AI with stronger assumptions, and then we found that the result holds in the current more general setting. All the AI assisted proofs and computations are verified and rewritten by us, and we take full responsibility for all the mathematical results and other text in this paper.

\section{General theory and examples}\label{sec:ex}

In this section, we establish LDPs of the empirical measure and extreme eigenvalue for external potentials defined on subsets, and develop a general theory to study lower tail large deviations of extreme eigenvalue. Let us start with some notation. Throughout, for two sequences $(x_n)_{n\ge1}$ and $(a_n)_{n\ge1}$, we write $x_n=O(a_n)$ if $|x_n|\le C|a_n|$ for some $C>0$ and all sufficiently large $n$, and $x_n=o(a_n)$ if $x_n/a_n\to0$.
We write $x_n=\Theta(a_n)$ if there exist constants $c,C>0$ such that $c|a_n|\le |x_n|\le C|a_n|$ for all sufficiently large $n$. Unless specified otherwise, we denote by $C_{a}$ a constant depending only on (the list of parameters) $a$.
Let $D\subset \cz$ be a closed set.  For $\mu\in\mathcal M_1(D)$ and $\varepsilon>0$, let
\begin{equation}\label{eq:BLball}
B(\mu,\varepsilon)=\{\nu\in\mathcal M_1(D): d(\mu,\nu)<\varepsilon\},
\end{equation}
where $d$ denotes the bounded Lipschitz metric. 
For a measurable map $T:\cz\to\cz$ and a probability measure $\mu$, we write $T_\#\mu =\mu\circ T^{-1}$ for the push-forward measure of $\mu$ under $T$.

\subsection{Preliminaries on logarithmic potential theory}

As mentioned above, a key step in our approach is to identify the minimizing measure. To this end, we employ the logarithmic potential theory \cite{ST97} and  record two  results. The first provides a characterization of the density and support of the one-dimensional constrained equilibrium measure. It is an immediate consequence of Theorems IV.1.10, IV.1.11, and IV.2.4, together with Corollary IV.1.12 in \cite{ST97}. Recall that $e^{-V}$ is an admissible weight on $D$ if $V:D\to(-\8,\8]$ is lower semicontinuous, $V(z)<\8$ on a set of positive logarithmic capacity and $\lim_{|z|\to\8, z\in D}V(z)-\log |z| = \8$ whenever $D$ is unbounded. In this case, by \cite{ST97}*{Theorem I.1.3} $I_V(\mu)$ admits a unique minimizer $\mu_V$ with compact support and $C_V\in \rz$, where $I_V$ and $C_V$ were defined in \eqref{defiv}.

\begin{lemma}\label{remr}
Let $D=[u_1,u_2]$, where $u_1=-\infty$ and $u_2=\infty$ are allowed, and let $\mu_D=\argmin_{\nu\in\mathcal M_1(D)} I_V(\nu).$ Assume that $w=e^{-V}$ is an admissible weight on $D$, and write $\mathcal S_D=\supp(\mu_D).$ 
Suppose that $\mathcal S_D=[v_1,v_2]\subset D$, that $w>0$ is absolutely continuous on $[v_1,v_2]$, and that $ w'\in L^p([v_1,v_2])$ for some $p>1.$ Then $\mu_D$ is locally absolutely continuous on $(v_1,v_2)$, and its density is given by
\begin{align}\label{eq:densd}
\sigma_D(x)
=
\frac{\pi-\ell(x)}{\pi^2\sqrt{(v_2-x)(x-v_1)}},
\quad
\ell(x):=\pv\int_{v_1}^{v_2}
\sqrt{(v_2-y)(y-v_1)} \frac{V'(y)}{x-y} \dd y.
\end{align}

Assume in addition that either $V$ is convex on $(u_1, u_2)$, or $x V^{\prime}(x)$ increases on $(u_1, u_2)$ when $\left[u_1, u_2\right) \subset[0, \infty)$. Then $\mathcal S_D$ is a compact interval, say $\mathcal S_D=[a,b]\subset [u_1,u_2].$
The support endpoints satisfy
\begin{align}
\frac1\pi\int_a^b V'(x)\sqrt{\frac{x-a}{b-x}}\dd x
&\begin{cases}=1,&b<u_2,\\ \leq1,&b=u_2,\end{cases}\label{eqb}\\
\frac1\pi\int_a^b V'(x)\sqrt{\frac{b-x}{x-a}}\dd x
&\begin{cases}=-1,&a>u_1,\\ \geq-1,&a=u_1.\end{cases}\label{eqa}
\end{align}
If $V$ is convex, these conditions characterize the support $[a,b]$. The candidate measure can also be verified using \eqref{eq:densd} and the equilibrium criterion in
Lemma \ref{lemr}.

Finally, if $D=\mathbb R$, $V$ is even, $V'$ exists on $(0,\infty)$, and $xV'(x)$ is positive and increasing on $(0,\infty)$, then $\mathcal S_D=[-a,a],$ where $a>0$ is characterized by 
\begin{align}\label{eqb2}
 \frac{2}{\pi}\int_0^1\frac{at V'(at)}{\sqrt{1-t^2}} \dd t =1. 
\end{align}
\end{lemma}

The second result gives a necessary and sufficient criterion for the extreme measure. In practice, once a candidate $\mu_\zeta$ is proposed, it suffices to check equilibrium condition \eqref{equcon}.
\begin{lemma}[Theorems I.1.3 and I.3.3 in \cite{ST97}] \label{lemr}
Let $e^{-V}$ be an admissible weight on the closed set $D$ and $\mu_V=\argmin_{\nu\in\mx_1(D)} I_V(\nu)$ the equilibrium measure associated with $V$ on $D$. Suppose that $\sigma\in \mathcal M_1(D)$ has compact support $S:=\supp(\sigma)$ and finite noncommutative entropy.
Then $\sigma=\mu_V$ if and only if there exists a constant $F_V\in\mathbb R$ such that
\begin{equation}\label{equcon}
V(z)- \int_D \log|z-t| \sigma(\dd t)
\begin{cases}
=F_V, & \text{q.e. on } S,\\
\ge F_V, & \text{q.e. on } D.
\end{cases}
\end{equation}
In particular, whenever \eqref{equcon} holds, one necessarily has $S=\supp(\mu_V)$.
\end{lemma}

The function $x\mapsto -\int \log |x-y|\mu(\dd y)$ is called the logarithmic potential of the compactly supported Borel probability measure $\mu$.
By \cite{ST97}*{Remark I.1.7}, measures of finite logarithmic energy do not charge sets of zero logarithmic capacity. Hence the equality in \eqref{equcon} holds $\mu_V$ almost everywhere and
\begin{equation}\label{eq:calsig}
\Sigma(\mu_V)
=
\int_D V(z) \mu_V(\dd z)-F_V.
\end{equation}
Moreover, if $\supp\mu_V=[\alpha,\beta]$ with $\alpha<\beta$,
averaging the equality in \eqref{equcon} against the arcsine law gives
\begin{equation}\label{eq:arcavg}
F_V
=\frac1\pi\int_\alpha^\beta
\frac{V(x)}{\sqrt{(\beta-x)(x-\alpha)}}\dd x
-\log\frac{\beta-\alpha}{4}.
\end{equation}
Here we used
$$
\int_\alpha^\beta
\frac{\log|x-y|}{\pi\sqrt{(\beta-x)(x-\alpha)}}\dd x
=\log\frac{\beta-\alpha}{4},
\qquad y\in[\alpha,\beta].
$$
\begin{remark}\label{rem:qepw}
In the setting of Lemma \ref{lemr}, assume that $D\subset\mathbb R$ is a closed interval with nonempty interior. 
Suppose that $V$ is finite and continuous in the interior and that its boundary values are its interior limits, possibly $+\infty$. Then the quasi-everywhere relations in \eqref{equcon} actually hold everywhere:
\begin{align}\label{eq:vlogmu}
 V(x)-\int\log|x-y| \mu_V(\dd y)
 \begin{cases}
  =F_V,&x\in \supp(\mu_V),\\
 \geq F_V,&x\in D.
 \end{cases}
\end{align}
Indeed, we first show that $x\mapsto\int\log|x-y| \mu_V(\dd y)$ is continuous on $D$.  \cite{ST97}*{Theorem I.4.4} gives continuity outside $\supp(\mu_V)$. 
Since $e^{-V}$ is a continuous admissible weight, \cite{ST97}*{Theorem I.5.1} applies.
If $D=\mathbb R$, each support point belongs to the  boundary of both components of $\mathbb C\setminus D$.
Thus Theorem I.5.1(ii) applies.
If $D$ is a bounded closed interval or a closed half line, $\mathbb C\setminus D$ is connected and every support point is a regular boundary point for the Dirichlet problem. Theorem I.5.1(iv) therefore applies.  Thus the logarithmic potential is continuous on $D$.
Together with the continuity and boundary limits of $V$, this
extends \eqref{equcon} to \eqref{eq:vlogmu}.
\end{remark}

\subsection{LDP of empirical measure with general potentials} \label{se:ldpln}

Our aim is to prove Theorem \ref{ldpln} in this section. We first recall the definition of large deviation principle (LDP). Let $(X_n)$ be random variables taking values in a Polish space $\mathbb X$, and let $a_n\to\infty$. 
We say that $(X_n)$ satisfies an LDP on $\mathbb X$ with speed $a_n$ and rate function $I:\mathbb X\to[0,\infty]$ if $I$ is lower semicontinuous and for every open set $O\subset\mathbb X$ and every closed set $F\subset\mathbb X$,
$$
\liminf_{n\to\infty}\frac{1}{a_n}\log\PP(X_n\in O)\ge -\inf_{x\in O}I(x),
\quad
\limsup_{n\to\infty}\frac{1}{a_n}\log\PP(X_n\in F)\le -\inf_{x\in F}I(x).
$$
A rate function $I$ is said to be good if the sublevel set $\{x\in X: I(x)\le M\}$ is compact for all $M\in \rz_+$.

To motivate our assumptions, let us examine the Laguerre ensemble considered in, e.g., \cites{HP98,Joh00}. In our notation, the model in \cite{HP98}*{Eq.~(4)} has density
\eqref{pdf} on $\mathbb R_+^n$ with fixed $\beta>0$ and
$$
 V_n(x)=\frac{p}{\beta n}Q(x)-\frac{\kappa_n}{\beta n}\log x,
 \qquad x>0.
$$
Here $Q$ is continuous on $\mathbb R_+$ and $Q(x)/\log x\to\infty$ as $x\to\infty$. 
The assumptions of \cite{HP98}*{Theorem~2.1} (see also \cite{HP00}*{Theorem 5.5.1}) become
\begin{equation}\label{eq:hp98}
 \frac np\to\gamma\in(0,1],\qquad
 \kappa_n\geq0,\qquad \frac{\kappa_n}{p}\to\kappa>0.
\end{equation}
For the Laguerre potential \eqref{def:vnl}, we have $Q(x)=\beta x/2$ and $\kappa_n=\beta(p-n+1)/2-1$, so $\kappa_n/p\to\beta(1-\gamma)/2$.
Thus the strict positivity in \eqref{eq:hp98} holds only when $\gamma<1$. 
In the real case $\beta=1$ (in fact for all $\bt<2$), we even have $\kappa_n<0$ when $p=n$. 
Then $V_n(x)\to -\8$ as $x\downarrow0$, but existing literature usually requires that the external potential $V_n$ take values in $(-\8,+\8]$ on $D$; see, e.g., \cite{ST97}.
In the complex case $\beta=2$ (and for all $\bt\ge 2$), $\kappa_n\ge p-n\geq0$, but $\kappa_n/p\to0$ when $\gamma=1$.
Neither critical case is therefore covered by \cite{HP98}*{Theorem~2.1}.
The strict positivity assumption $\kappa>0$ is used in \cite{HP98}*{Lemma~3.1} to obtain uniform convergence of the truncated kernels near zero.
As a result, \cite{HP98}*{Theorem~4.1} and the subsequent discussion of complex Wishart matrices exclude the critical case $\gamma=1$. See Example \ref{x:lagemp} below and \cite{WZschatten}*{Lemma~3.11}.  On the other hand, the general assumptions on $V_n$ as given in \cites{Joh00,BG13} require uniform convergence of $V_n$ to $V$ on (compact subsets of) $D=[0,\8)$. This is possible only if $\kappa_n=0$ or $p-n = \frac2\bt-1$, yielding a restrictive condition on $\bt$.

These boundary difficulties motivate the following decomposition. The finite dimensional Laguerre and Jacobi potentials in \eqref{def:vnl} and \eqref{def:vnj} contain logarithmic terms whose coefficients may be of order $(\beta n)^{-1}$.  
In certain parameter regimes, these terms prevent $V_n$ from admitting a lower semicontinuous extension to the boundary.  
We therefore isolate them as a one particle weight and write
\begin{equation}\label{eq:fulldeco}
 V_n(z)=W_n(z)-\frac1{\beta n}\log\omega_n(z).
\end{equation}
We remark that the singularity does not have to blow up with speed $\log z$ at $z=0$ even though we write $\log \omega_n$. For instance, as long as 0 is a boundary point of $D$, we may take $\omega_n(z)= e^{1/|z|}$ for $z\neq 0$ which is not even integrable near 0.
Under this representation, the density \eqref{pdf} has the equivalent weighted form
\begin{equation}\label{eq:wpdf}
\mathbb P_n(\dd\boldsymbol\lambda)
=
\frac1{Z_n}
\prod_{1\leq i<j\leq n}|\lambda_i-\lambda_j|^{\beta}
\exp\left\{-\beta n\sum_{i=1}^nW_n(\lambda_i)\right\}
\prod_{i=1}^n\omega_n(\lambda_i) \dd\boldsymbol\lambda.
\end{equation}
This formulation allows the lower semicontinuity, convergence, and growth of $W_n$ and the contribution of $\omega_n$ to be treated separately. In observance of the eigenvalue density formulas of Laguerre and Jacobi ensembles, the decomposition \eqref{eq:fulldeco} is a natural choice.

\begin{assumption}\label{ass:ldp}\rm
Let $D\subset\R$ or $D\subset\C$ be closed, and consider the   density \eqref{eq:wpdf} on $D^n$. 
\begin{enumerate}[label=(A\arabic*),start=0]
\item\label{A0}
$\beta n\to+\infty$.

\item\label{A1}
The functions $W_n:D\to\mathbb R\cup\{+\infty\}$ and $V:D\to\mathbb R\cup\{+\infty\}$ are lower semicontinuous.  
Let $\mathcal O$ denote the interior of $D\cap\{V<+\infty\}$ and assume that $\mathcal O$ is nonempty.  
Moreover,
\begin{enumerate}[label=\textup{(\roman*)}]
\item for every sequence $z_n\in D$ with $z_n\to z\in D\setminus\mathcal O$, $\liminf_{n\to\infty}W_n(z_n)\geq V(z);$
\item $W_n \to V$ uniformly on compact subsets of $\mathcal O$ and $V$ is continuous on $\mathcal O$.
\end{enumerate}

\item\label{A2}
If $D$ is unbounded, assume that
\begin{equation*}
 \liminf_{|z|\to\infty}
 \inf_{n\geq1}\frac{W_n(z)}{\log(1+|z|)}>1,
 \qquad
 \liminf_{|z|\to\infty}
 \frac{V(z)}{\log(1+|z|)}>1.
\end{equation*}

\item\label{A3}
The functions $\omega_n:D\to(0,+\infty)$ are measurable and finite almost everywhere.  Moreover,
\begin{enumerate}[label=\textup{(\roman*)}]
\item for every compact subset $K$ of $\mathcal O$,
$\log\omega_n(z) =o(\beta n)$ uniformly for a.e.~$z\in K;$
\item there exists a sequence $\theta_n\to0$ such that $\theta_n\beta n\to+\infty$, and
\begin{equation*}
 \limsup_{n\to\infty}
 \frac1{\beta n^2}\log\mathcal R_n(\theta_n)\leq0,
\end{equation*}
where
\begin{align}\label{eq:def-frac}
 \mathcal R_n(\theta)
 =\int_{D^n}
 &\prod_{1\leq i<j\leq n}|z_i-z_j|^{\theta\beta}
 \exp\left\{-\theta\beta n
       \sum_{i=1}^nW_n(z_i)\right\}
 \prod_{i=1}^n\omega_n(z_i) \dd\boldsymbol z.
\end{align}
\end{enumerate}
\end{enumerate}
\end{assumption}

We allow that $\bt$ depends on $n$, but keep this dependence implicit to make the notation simple. Let us make some remarks on the conditions.

\begin{remark}
In Condition \ref{A1}, semicontinuity of $V$ and $W_n$ includes exotic potentials like $V_n(x)=\frac12x^2+1+\sin \frac1x$ for $x>0$ and $V_n(0)=0$; see also the examples in Remarks \ref{rem:mxpow} and \ref{re:lamindiv} below. The one particle weight $\omega_n$ is assumed to be only a.e.~finite so that we can always set $\omega_n(z)=1$ if $\lim_{x\to z} \omega_n(x)=+\8$. This is especially convenient when $z\in \partial \mathcal O$.
We can also enlarge $D$ to be $\tilde D$ by setting $V(x)=+\8$ on $\tilde D\setminus D$.  In practice, we usually take $V$ to be continuous on $D$ and finite on the interior of $D$, which coincides with Assumption \ref{C1} below.

For $z_n\to z$ and $w_n\to w$, condition \ref{A1} gives
$$
 \liminf_{n\to\infty}
 \left\{
 \frac12W_n(z_n)+\frac12W_n(w_n)
 -\frac12\log|z_n-w_n|
 \right\}
 \geq
 \frac12V(z)+\frac12V(w)-\frac12\log|z-w|.
$$
After truncating the pair kernel from above, this is precisely the lower limit needed to apply Fatou's lemma for varying measures \cite{WFatou} in the local LDP upper bound.

The seemingly more natural condition $\liminf_{n\to\8} W_n(z)\ge V(z)$ for all $z\in D$ is not enough. 
Indeed, consider $\beta=1,$ $D=[0,1]$, $W_n(z)=-\log n \mathbf{1}_{[\frac1{\log n},\frac2{\log n}]} (z)$, $\omega_n(z)\equiv 1$ and $V(z)\equiv 0$. 
Then Condition \ref{A1} holds with \ref{A1}(i) replaced by  $\liminf_{n\to\8} W_n(z)\ge V(z)$ for all $z\in D$, but $-1=\int W_n(x)\dd x<\int V(x)\dd x$ and Fatou's lemma does not apply here. 
Moreover, by change of variable, we see that $$Z_n\geq n^{n^2} (\log n)^{-n(n+1)/2} \int_{[0,1]^n}\prod_{i<j}|u_i-u_j| \prod_i \dd u_i,$$ which yields $n^{-2}\log Z_n \to +\8$ contradicting Proposition \ref{pr:logzn} below.
\end{remark}

\begin{remark}\label{r:a3unit}
If $\omega_n\equiv1$, then \ref{A3} follows from \ref{A0}--\ref{A2}.  
Indeed, \ref{A3}(i) is immediate, while \ref{A1}  and \ref{A2} give, for some $\varepsilon,C>0$ and all sufficiently large $n$,
$$
 W_n(z)\geq(1+\varepsilon)\log(1+|z|)-C,
 \qquad z\in D.
$$
Let $\theta_n=(\beta n)^{-1/2}$.  
Using $|z-w|\leq(1+|z|)(1+|w|)$, we obtain
\begin{align*}
 \mathcal R_n(\theta_n)
 &\leq e^{C\theta_n\beta n^2}
 \left[\int_D(1+|z|)^{-\theta_n\beta(1+\varepsilon n)}
 \dd z\right]^n
 \leq e^{C\theta_n\beta n^2+Cn}.
\end{align*}
Therefore $\limsup_{n\to\infty}(\beta n^2)^{-1}\log\mathcal R_n(\theta_n)\leq0$ and \ref{A3}(ii) follows. In fact, a primary reason for introducing \ref{A3} is to control the non-trivial one particle weight $\omega_n$ using a negligible portion of the potential $W_n$ and logarithmic interaction. At the same time, it allows us to assume only $\bt n\to\8$ without any growth rate, which is a byproduct of our approach.

\end{remark}

\begin{remark} 
We can draw connection to \cite{GZ19} through the weight $\omega_n$. Indeed, when $\omega_n=\omega$ does not depend on $n$ and is integrable, then $\pi(\dd x) \propto \omega(x)\dd x$ gives a reference probability measure on $D$ and we may apply the framework of \cite{GZ19} (the lower limit assumption there still needs to be checked, though). However, for general $\omega_n$ it could be impossible to construct a probability measure $\pi$ and a sequence of potentials uniformly bounded below (termed stable sequence there) based on our potentials $W_n$. We give an example for illustration.

Consider \eqref{eq:wpdf} with $D=[0,1]$, $\beta_n=n^{-1}(\log(n+1))^{1/4}$, $W_n=0$ and  $\omega_n(x)=1+|I_n|^{-1}\mathbf1_{I_n}(x),$ where $I_n=((n+1)^{-1},n^{-1}]$.  It is easy to verify \ref{A0}--\ref{A3}(i) and Assumption \ref{ass:d}.  Using $|x-y|\leq1$ on $D$, we have for every $\theta>0$,
\begin{equation}\label{eq:gzfrac}
 \mathcal R_n(\theta)
 =\int_{D^n}\prod_{1\leq i<j\leq n}|x_i-x_j|^{\theta\beta_n}
       \prod_{i=1}^n\omega_n(x_i)\dd\boldsymbol x
 \leq\left(\int_D\omega_n(x)\dd x\right)^n=2^n.
\end{equation}
Thus \ref{A3}(ii) holds.
By Theorem \ref{ldpln}, we have an LDP of the empirical measure $n^{-1}\sum_{i=1}^n \delta_{\lambda_i}$ at speed
$\beta_n n^2$.

Next, we show that \cite{GZ19}*{Corollary 1.3} does not apply directly to this model at speed $\beta_n n^2$. For contradiction, suppose that the corollary applied at speed $\beta_n n^2$ with a fixed probability measure $\pi$ on $D$. Then the density of $(\lambda_1,\ldots,\lambda_n)$ can be written as 
$$\dd\mathbb P_n=\widetilde Z_n^{-1} e^{-\beta_n n^2 H_n} \dd\pi^{\otimes n}$$
for some function $H_n$ of $n$ variables.
Stability condition (S) requires $H_n\geq-C$ and \cite{GZ19}*{Theorem~1.2} ensures a finite limit for $\frac1{\beta_n n^2}\log \widetilde Z_n$. 
Consequently, there exist $C'>0$ and $N$ such that for $n>N$
\begin{equation}\label{eq:gzdom}
 \mathbb P_n(E_n)
 =\widetilde Z_n^{-1}\int_{E_n}e^{-\beta_n n^2 H_n}\dd\pi^{\otimes n}
 \leq e^{C'\beta_n n^2}\pi^{\otimes n}(E_n).
\end{equation}
However, this bound fails. Indeed, let $m=\lfloor n/\sqrt{\log(n+1)}\rfloor$ and consider $E_n=\bigl\{\#\{i:\lambda_i\in I_n\}\geq m\bigr\}.$
Let $X_1,\ldots,X_n$ be independent, with $X_1,\ldots,X_m$ uniform on $I_n$ and $X_{m+1},\ldots,X_n$ uniform on $J=[1/2,1]$.
Then by Jensen's inequality
\begin{align*}
 \mathbb P_n(E_n)
 \geq \mathbb P_n(I_n^m\times J^{n-m})
 \geq \frac{|J|^{n-m}}{Z_n}
 \mathbb E e^{\beta_n\sum_{i<j}\log|X_i-X_j|}
 \geq \frac{|J|^{n-m}}{Z_n} e^{\beta_n\sum_{i<j}\mathbb E\log|X_i-X_j|}.
\end{align*}
It follows that
$$ \log\mathbb P_n(E_n)
 \geq -\log Z_n+(n-m)\log|J|
      +\beta_n\binom m2\log|I_n|-C\beta_n n^2.$$
Taking $\theta=1$ in \eqref{eq:gzfrac} gives
$Z_n=\mathcal R_n(1)\leq 2^n$.
Together with $m^2\log(1/|I_n|)=O(n^2)$ and $n=o(\beta_n n^2)$, we get
\begin{equation}\label{eq:gzeventlb}
 \log\mathbb P_n(E_n)\geq -C_1 \beta_n n^2.
\end{equation}
On the other hand, the intervals $I_n$ are disjoint, hence $\sum_n\pi(I_n)\leq1$.
Since $\sum_n 1/n=\infty$, we can choose $n_k\uparrow\infty$ such that $\pi(I_{n_k})\leq1/n_k$. We deduce as $k\to\infty$
\begin{equation*}*\label{eq:gzeventub}
 \frac{1}{\beta_{n_k} n_k^2}\log\mathbb P_{n_k}(E_{n_k})
 \leq C'+\frac{1}{\beta_{n_k} n_k^2 }\log\left[\binom{n_k}{m_k}\pi(I_{n_k})^{m_k}\right]
 \leq C'+\frac{n_k\log 2-m_k\log n_k}{\beta_{n_k} n_k^2}
  \to -\infty,
\end{equation*}
contradicting the lower bound
\eqref{eq:gzeventlb} at $n=n_k$. 
\end{remark}

\begin{remark}
The classic setting of \cite{AGZ10}*{Theorem 2.6.1} is a special case of Assumption \ref{ass:ldp}. 
There, $\beta>0$ is fixed and $V$ is a fixed continuous potential on $\mathbb R$. 
Under normalization \eqref{pdf}, the growth condition reads
$$
\liminf_{|x|\to\infty}
\frac{V(x)}{\log|x|}
>
1\vee\beta^{-1}.
$$
Taking $V_n=W_n=V$ and $\omega_n\equiv1$, condition \ref{A0} is immediate, \ref{A1} follows from the continuity of $V$, and the above growth condition implies \ref{A2}. 
Condition \ref{A3} then follows from the preceding remark.

Extending the argument of \cite{AGZ10} to varying $\beta$ and $V_n$, particularly when $\beta\to0$, requires the one particle estimate
\begin{align}\label{bd:logvn}
\limsup_{n\to\infty}
\frac1{\beta n}
\log\int_D e^{-\beta V_n(x)}\dd x
\leq0.
\end{align}
This estimate does not follow from $\beta n\to\infty$. 
Indeed, for
$D=\mathbb R$ and $V_n(x)=x^2/2$,
$$
\log\int_{\mathbb R}e^{-\beta x^2/2}\dd x
=
\frac12\log\frac{2\pi}{\beta}.
$$
The stronger condition $\beta n/\log n\to\infty$ ensures \eqref{bd:logvn}, whereas $\beta n\to\infty$ alone does not. 
To work under the latter condition, we use \eqref{eq:factor} to place a fraction $\theta_n$ of the interaction and confinement in $\mathcal R_n(\theta_n)$. 
By \ref{A3}(ii), this factor is negligible at speed $\beta n^2$. 
The proof in \cite{Har12} also uses a one particle estimate corresponding to \eqref{bd:logvn}. 
For varying $\beta$, this estimate does not follow from $\beta n\to\infty$. 
Even for a quadratic potential, the stronger condition $\beta n/\log n\to\infty$ is needed when $\beta\to0$.
\end{remark}

\begin{remark}\label{le:uniq} 
When $D$ is unbounded, the growth condition of $V$ in Assumption \ref{A2} is equivalent to  $$\liminf_{|z|\to\infty} \left\{ V(z) - a \log|z| \right\} > -\infty$$ for some  $a>1$.
Thus \cite{ST97}*{Theorem I.1.3} gives a unique minimizer $\mu=\argmin_{\nu\in\mathcal M_1(D)}I_V(\nu)$ with compact support. 
One may relax the assumption \ref{A2} to allow unbounded support for $\mu$ (see, e.g., \cite{Har12}); however, we do not pursue this direction, since our technical tools from \cite{ST97} require compact support of $\mu$.
\end{remark}

\begin{assumption} \label{ass:d}\rm
Let $\mu\in\mathcal M_1(D)$ be given with $I_V(\mu)<\infty.$ Assume that there exists a sequence of probability measures $\mu_k(\dd z)=g_k(z)\dd z, g_k\in L^\infty(D),$ such that $\supp \mu_k$ is a compact subset of the interior of $D \cap\{V<\infty\},$
$$
\mu_k\Rightarrow\mu,
\qquad 
\liminf_{k\to\infty}\Sigma(\mu_k)
\ge
\Sigma(\mu),
\qquad \text{and} \qquad
\limsup _{k \rightarrow \infty} \int V \dd \mu_k \leq \int V \dd \mu.
$$
\end{assumption}

This assumption is motivated by \cite{PH98}; see also \cite{CGZ14}. It is more like a condition on the closed set $D$ and is needed for lower bounds. By Lemma \ref{le:symreg}, it is satisfied for reasonably good probability measures on convex sets and annuli, which include all the scenarios we consider. In short, we require regularized approximating measures supported in an open set where $V$ is finite, which in turn indicates that $D^\circ$ has strictly positive Lebesgue measure. 
Here by default, the interior of $D$ and Lebesgue measure are compatible with the ambient space of $D$. 
Since we usually restrict integrals to smaller domains when proving lower bounds, this condition is natural although other choices are possible; see, e.g., \cite{AGZ10}*{Eq.~(2.6.27)}. 

\begin{remark}\label{re:assub}
To see the necessity of Assumption \ref{ass:d},  consider the Coulomb gas \eqref{eq:wpdf} on $D=\mathbb{D}\cup\{z\in\mathbb C:|z-3|=1\}$ with fixed $\beta>0$, $W_n=V=0$, and $\omega_n\equiv1$, where $\mathbb{D}=\{z\in\cz: |z|\le 1\}$ is the closed unit disk. Then Assumption \ref{ass:ldp} holds, where \ref{A3} follows from Remark \ref{r:a3unit}.
Since the circle $\{|z-3|=1\}$ has zero planar Lebesgue measure, all particles lie in $\{|z|\leq1\}$ a.s. In contrast, the uniform probability measure $\mu$ on that circle satisfies $\Sigma(\mu)=0$ and $I_0(\mu)<\infty$.
Let 
\[
f(z)=\frac12\min\{1,\operatorname{dist}(z,\mathbb{D}) \}.
\]
Then $f$ is bounded Lipschitz with $\int f\dd\mu=1/2$ and $\int f\dd L_n=0$ a.s. It follows that $d(L_n,\mu)\geq1/2$. Therefore, we have $\mathbb P_n(L_n\in B(\mu,1/4))=0$ for every $n$. 
Consequently,
$$
\liminf_{n\to\infty}\frac1{\beta n^2}
\log\mathbb P_n(L_n\in B(\mu,1/4))
=-\infty<-\inf_{\nu\in B(\mu,1/4)}I_0(\nu),
$$
violating the LDP lower bound. 
Here Assumption \ref{ass:d} fails because the uniform measure $\mu$ cannot be approximated by measures with densities on compact subsets of $D^\circ$.
 
Similar situations occur in one dimension. Indeed, let $D=[0,1]\cup(2+\mathcal C)$, where $\mathcal C=\left\{\sum_{k=1}^{\infty}\frac{2a_k}{3^k}: a_k\in\{0,1\}\right\}$ is the Cantor set.
Consider \eqref{eq:wpdf} with fixed $\beta>0$, $W_n=V=0$, and $\omega_n\equiv1$. Then Assumption \ref{ass:ldp} again holds.
Since $|\mathcal C|=0$, all particles lie in $[0,1]$ a.s. Let $\mu$ be the law of $2+\sum_{k\geq1}2B_k/3^k$, where the $B_k$ are independent Bernoulli variables with parameter $1/2$. Namely, $\mu$ is the Cantor measure on $2+\cx$.
For two independent samples $X,Y$ from $\mu$, let $T$ be the first index at which their digit sequences differ.
Then $\mathbb P(T=k)=2^{-k}$ and $|X-Y|\geq3^{-T}$. It follows that
$$
-\Sigma(\mu)=\mathbb E[-\log|X-Y|]
\leq(\log3)\sum_{k=1}^{\infty}k2^{-k}
=2\log3<\infty.
$$
Thus $I_0(\mu)<\infty$. Choosing $f(x)=\frac12\min\{1,\operatorname{dist}(x,[0,1])\}$ yields $d(L_n,\mu)\geq1/2$ a.s. 
Consequently,
$$
\liminf_{n\to\infty}\frac1{\beta n^2}
\log\mathbb P_n(L_n\in B(\mu,1/4))
=-\infty<-\inf_{\nu\in B(\mu,1/4)}I_0(\nu).
$$
Assumption \ref{ass:d} fails again as every measure with a density on $D$ is supported on $[0,1]$ and cannot approximate $\mu$ weakly.

On the other hand, our framework still allows us to evaluate the lower tail large deviation rate function on the `regular part' of the region in these examples. For instance, the eigenvalue density on $[0,1]$ in the second example is just the complex Jacobi ensemble with $p_1=p_2=n$, which even has exact formula of confined probability for finite $n$.
\end{remark}

\begin{lemma} \label{le:symreg}
Let $D$ be a closed subset of $\R$ or $\C$, and let $E\subset D$ satisfy one of the following conditions:
\begin{enumerate}[label=\textup{(\alph*)}]
\item $E$ is closed and convex with nonempty interior in the ambient space of $D$;
\item $E=\{z\in\C:r_1\leq|z|\leq r_2\}$, where
$0<r_1<r_2\leq\infty$.
\end{enumerate}
Let $V:D\to\R\cup\{+\infty\}$ be lower semicontinuous and bounded below.
Assume that $V$ is finite on the interior of $E$ and continuous on $E\cap\{V<+\infty\}$. Let $\mu\in\mx_1(D)$ satisfy
$$
\mu(E)=1,\qquad
\int_D V\dd\mu<\infty,
\qquad
\Sigma(\mu)>-\infty.
$$
Then there exists a sequence
\begin{equation}\label{eq:regseq}
\mu_k(\dd x)
=
g_k(x)\dd x,
\qquad
g_k\in C_c^\infty(E^\circ),
\end{equation}
of probability measures on $D$ such that
\begin{equation}\label{eq:reglim}
\mu_k\Rightarrow\mu,
\qquad
\liminf_{k\to\infty}\Sigma(\mu_k)\geq\Sigma(\mu),
\qquad
\limsup_{k\to\infty}\int_DV\dd\mu_k
\leq
\int_DV\dd\mu.
\end{equation}
Here $C_c^\8(E^\circ)$ denotes the set of smooth functions supported on a compact subset of $E^\circ$.
Moreover, if $E\subset\C$ satisfies $E=\bar E$ and $\mu=\bar\mu$, then the sequence may be chosen so that $\mu_k=\bar\mu_k$ for every $k\geq1$, where $\bar\mu(A)=\mu(\bar A)$ for every Borel set $A\subset\C$ and $\bar A$ is the complex conjugate of the set $A$.
\end{lemma}

\begin{proof}
By the definition of $\Sigma$ in \eqref{def:sigmu}, the condition $\Sigma(\mu)>-\infty$ implies $\int_D\log(1+|x|)\mu(\dd x)<\infty$. 
Using \eqref{e:pairbd}, $\iint_{D^2} \log_+|x-y| \mu(\dd x)\mu(\dd y)<\8$; the negative part of the logarithmic kernel is also integrable as $\Sigma(\mu)>-\infty$. Thus
\begin{equation}\label{eq:logabs}
\iint_{D^2}|\log|x-y|| \mu(\dd x)\mu(\dd y)<\infty.
\end{equation}
Since $\mu(E)=1$ and $\int_DV\dd\mu<\infty$, $\mu\bigl(E\cap\{V<+\infty\}\bigr)=1$.
Choose increasing compact sets $K_m\subset E\cap\{V<+\infty\}$ with $\mu(K_m)>0$ and $\mu(K_m)\to1$. If $E\subset\C$ is invariant under complex conjugation and $\mu=\bar\mu$, then also $\mu\bigl(\{z\in E:V(\bar z)<+\infty\}\bigr)=1.$
Thus the conjugation invariant set $\{z\in E:V(z)<+\infty,\ V(\bar z)<+\infty\}$ has full $\mu$-measure, and we may choose $K_m=\bar K_m$.
Let $\mu^{(m)}=\mu|_{K_m}/\mu(K_m)$.  Since $V$ is bounded below, dominated convergence and \eqref{eq:logabs} give, as $m\to\infty$,
\begin{equation}\label{eq:trclim}
\mu^{(m)}\Rightarrow\mu,
\qquad
\int_DV\dd\mu^{(m)} \to \int_DV\dd\mu,
\qquad
\Sigma(\mu^{(m)}) \to \Sigma(\mu).
\end{equation}
We next move the compact supports into the interior of $E$.
\begin{enumerate}[label=\textup{(\roman*)}]
\item[(a)] If $E$ is convex, choose $x_0\in E^\circ$ and let
$T_\varepsilon(x)=(1-\varepsilon)x+\varepsilon x_0$ for $0<\varepsilon<1$. In the symmetric two-dimensional setting, take $x_0\in E^\circ\cap\R$. Then $|T_\varepsilon(x)-T_\varepsilon(y)|=(1-\varepsilon)|x-y|$, while for every compact $K\subset E\cap\{V<+\infty\}$, $\sup_{x\in K}|T_\varepsilon(x)-x|=\varepsilon\sup_{x\in K}|x_0-x|\to0$ as $\varepsilon\to0$.

\item[(b)] If $E=\{z\in\C:r_1\leq|z|\leq r_2\}$ with $r_1>0$, choose $r_0\in(r_1,r_2)$ and let $T_\varepsilon(z)=(1-\varepsilon+\varepsilon r_0/|z|)z$. For $z=re^{\mathrm i\theta}$ and $w=se^{\mathrm i\phi}$,
\begin{align*}
|T_\varepsilon(z)-T_\varepsilon(w)|^2
&=(1-\varepsilon)^2(r-s)^2+2\bigl((1-\varepsilon)r+\varepsilon r_0\bigr)
\bigl((1-\varepsilon)s+\varepsilon r_0\bigr)
\bigl(1-\cos(\theta-\phi)\bigr)
\\
&\geq(1-\varepsilon)^2|z-w|^2.
\end{align*}
Thus $|T_\varepsilon(z)-T_\varepsilon(w)|\geq(1-\varepsilon)|z-w|$, and
$\sup_{z\in K}|T_\varepsilon(z)-z|=\varepsilon\sup_{z\in K}|r_0-|z||\to0$ for every compact
$K\subset E\cap\{V<+\infty\}$.
\end{enumerate}

In both cases, $T_\varepsilon(K_m)$ is a compact subset of the
interior of $E$. Continuity of $V$ on $E\cap\{V<+\infty\}$ gives
$$
\sup_{x\in K_m}|V(T_\varepsilon(x))-V(x)| \to 0
\qquad\text{as }\varepsilon\downarrow0.
$$
For fixed $m$, let $\widetilde\mu_{m,\varepsilon}=(T_\varepsilon)_\#\mu^{(m)}$.
The preceding estimates imply
\begin{align}
\Sigma(\widetilde\mu_{m,\varepsilon})
\geq
\Sigma(\mu^{(m)})+\log(1-\varepsilon), \quad
\int_DV\dd\widetilde\mu_{m,\varepsilon}
\to
\int_DV\dd\mu^{(m)}, \quad \text{ as } \varepsilon\downarrow0.
\label{eq:pullV}
\end{align}
For
$
0<\delta<\frac12
\dist\bigl(\supp\widetilde\mu_{m,\varepsilon},(E^\circ)^c\bigr),
$
define
\begin{equation*}
\varphi_\delta(u)
=
c_\delta
\exp\left\{
-\frac{1}{\delta^2-|u|^2}
\right\}
\mathbf 1_{\{|u|<\delta\}}.
\end{equation*}
Choose $c_\delta>0$ so that $\int\varphi_\delta(u)\dd u=1$. Here $u\in\R$ in dimension one and $u\in\C$ in dimension two.
Set
$\mu_{m,\varepsilon,\delta}=\varphi_\delta*\widetilde\mu_{m,\varepsilon}$.
Then $\mu_{m,\varepsilon,\delta}$ has a bounded smooth density with compact support in the interior of $E$.
Since $V$ is continuous on a neighborhood of this compact support,
\begin{equation}\label{eq:molim}
\mu_{m,\varepsilon,\delta}
\Rightarrow
\widetilde\mu_{m,\varepsilon},
\qquad
\int_DV\dd\mu_{m,\varepsilon,\delta}
\to
\int_DV\dd\widetilde\mu_{m,\varepsilon}
\end{equation}
as $\delta\downarrow0$. Moreover, by the concavity of $\Sigma$ and Jensen's inequality,
\begin{equation}\label{eq:jenmu}
\Sigma(\mu_{m,\varepsilon,\delta})\geq\Sigma(\widetilde\mu_{m,\varepsilon}).
\end{equation}
Combining \eqref{eq:trclim}, \eqref{eq:pullV},
\eqref{eq:molim}, and \eqref{eq:jenmu}, using a diagonal choice, we find that
$
\mu_k:=\mu_{m_k,\varepsilon_k,\delta_k}
$
satisfies \eqref{eq:regseq} and \eqref{eq:reglim}.
Since $E^\circ\subset\operatorname{int}(D\cap\{V<+\infty\})$,
the supports of these $\mu_k$ verify Assumption \ref{ass:d} on $D$.

Finally, suppose $E=\bar E$ and $\mu=\bar\mu$. Choose $K_m=\bar K_m$ as above; take $x_0\in E^\circ\cap\R$ in the convex case. Then $\mu^{(m)}=\overline{\mu^{(m)}}$ and
$T_\varepsilon(\bar z)=\overline{T_\varepsilon(z)}$. It follows that $\widetilde\mu_{m,\varepsilon}=\overline{\widetilde\mu_{m,\varepsilon}}$. Since $\varphi_\delta$ is radial,
$\mu_{m,\varepsilon,\delta}=\overline{\mu_{m,\varepsilon,\delta}}$. We conclude that $\mu_k=\bar\mu_k$ for every $k\geq1$.
\end{proof}

The decomposition \eqref{eq:fulldeco} allows singular points at which the original potential fails
\begin{align}\label{bd:vnv}
 \liminf_{n\to\infty}V_n(z_n)\geq V(z)
 \qquad\text{whenever }z_n\to z.
\end{align}
Let us return to the motivational $\beta$-Laguerre ensemble with $n/p\to1$.
The potentials $V_n^L$ in \eqref{def:vnl} converge to $x/2$ locally uniformly on $(0,\infty)$.
If $\beta(p-n+1)<2$, then $V_n^L(x)\to-\infty$ as $x\downarrow0$. Clearly, one can choose a sequence $x_n\downarrow0$ such that $V_n^L(x_n)\to-\infty$.
Thus \eqref{bd:vnv} fails at $z=0$. 
This also explains why \cite{HP98} proved the empirical measure LDP of the Laguerre ensemble only for $n/p \to\gamma\in(0,1)$.  In our framework, we write
$$
 W_n^L(x)
 =\frac{p}{2n}x-\frac{p-n+1}{2n}\log x,
 \qquad
 V_n^L(x)=W_n^L(x)+\frac1{\beta n}\log x.
$$
With $W_n^L(0)=+\infty$, the bound $W_n^L(x)\geq px/(2n)$ for $0<x<1$ gives \ref{A1}(i) at $0$.

In the normalization \eqref{def:vnl}, when $\gamma=0$, the empirical measure in the original coordinates converges to the Dirac mass at $1$; see the moment estimates in \cite{Dum03}*{Theorem~6.3.1} for fixed $\beta>0$. The work \cite{JL15} proved an LDP with speed $n^2$ for $n^{-1}\sum_{i=1}^n\delta_{\sqrt{\beta p/(2n)}(\lambda_i-1)}$ when $\beta>0$ is fixed and $n^2/p\to0$. When $\beta$ varies with $n$, $\beta n/\log n\to\infty$, and $n/p\to\gamma\in[0,1)$,  an LDP with speed $\beta n^2$ for $L_n=n^{-1}\sum_{i=1}^n\delta_{\sqrt{p/n}(\lambda_i-1)}$ was established in \cite{Ma23L}.
The following example requires only $\beta n\to\infty$ and also covers the critical case $\gamma=1$.
For $0<\gamma\leq1$, let $C_{V_\gamma^L}=((1-\gamma)^2\log(1-\gamma)+\gamma)/(4\gamma^2)$, $D_\gamma=[-\gamma^{-1/2},\infty)$ and
$$
V_\gamma^L(x)=\frac{x}{2\sqrt\gamma}
-\frac{1-\gamma}{2\gamma}\log(1+\sqrt\gamma x),
\qquad x>-\gamma^{-1/2},
$$
with the endpoint value defined by the interior limit. For $\gamma=0$, let $C_{V_0^L}=3/8$, $D_0=\mathbb R$ and $V_0^L(x)=x^2/4$. 
For $0<\gamma\leq1$, also define $$\widehat V_\gamma^L(x)=x/(2\gamma)-(1-\gamma)(2\gamma)^{-1}\log x, \qquad x>0$$ and define $\widehat V_\gamma^L(0)$ by its interior limit.
\begin{example}\label{x:lagemp}
For the $\beta$-Laguerre ensemble with potential \eqref{def:vnl}, assume that $\beta n\to\infty$.
\begin{itemize}
    \item[(i)] For $\gamma\in[0,1]$, the rescaled empirical measure $L_n=n^{-1}\sum_{i=1}^n\delta_{\sqrt{p/n}(\lambda_i-1)}$ satisfies an LDP on $\mathcal M_1(\mathbb R)$ with speed $\beta n^2$ and good rate function
\begin{equation}\label{eq:lagfun}
 I_{V_\gamma^L}(\mu)
 =\int V_\gamma^L\dd\mu-\frac12\Sigma(\mu)-C_{V_\gamma^L},
\end{equation}
where $I_{V_\gamma^L}(\mu)=+\infty$ if $\mu(D_\gamma)<1$.
\item[(ii)] For $\gamma \in (0,1]$, the unscaled empirical measure $\widehat L_n=n^{-1}\sum_{i=1}^n\delta_{\lambda_i}$ satisfies an LDP on $\mathcal M_1(\mathbb R_+)$ at the same speed with good rate function $I_{\widehat V_\gamma^L}$. The two rate functions satisfy
\begin{equation}\label{eq:loremp}
 I_{\widehat V_\gamma^L}(\mu)
 =I_{V_\gamma^L}(\varphi_\#\mu),\qquad
 \varphi(x)=\frac{x-1}{\sqrt\gamma}.
\end{equation}
\end{itemize}
\end{example}

For the Jacobi ensemble, \cite{HP06} proves the empirical measure LDP for $\beta=2$ under $n/p_1\to\gamma\in(0,1]$ and $p_1/p_2\to\sigma\in(0,\gamma^{-1}]$.
For fixed $\beta>0$ and $\gamma \sigma>0$, \cite{ES08} also gives  the empirical measure LDP under conditions $\beta(p_1-n+1)=2$ if $\gamma=1$ and $\beta(p_2-n+1)=2$ if $\gamma\sigma=1$. No additional conditions are needed when $\gamma<1$ and $\gamma\sigma<1$.
For fixed $\beta>0$ and $\gamma\sigma=0$, \cite{Nag14} proves almost sure semicircle and Marchenko--Pastur limits for rescaled empirical measures under additional parameter conditions, as well as LDPs for weighted spectral measures. 
For varying $\beta$, \cite{Ma23J} establishes an LDP with speed $\beta n^2$ for  $L_n=n^{-1}\sum_{i=1}^n\delta_{\varphi_n(\lambda_i)}$ under $\beta n/\log n\to\infty$,  $n/p_1\to\gamma\in[0,1)$, and $p_1/p_2\to\sigma\geq0$ with $\gamma\sigma<1$.
Here $\varphi_n$ is defined in \eqref{defvarphi}.
Using the decomposition \eqref{eq:fulldeco}, we obtain the following example, which requires only $\beta n\to\infty$ and also covers $\gamma=1$ or $\gamma\sigma=1$.
Let
$$
 C_{V_{\gamma,\sigma}^J}
 =\frac14\left[g(\gamma)+g(\gamma\sigma)
 -g\left(\frac{\sigma\gamma}{1+\sigma}\right)+\log(1+\sigma)\right],
$$
where $g(x)=((1-x)^2\log(1-x)+x)/x^2$ for $0<x\leq1$, and $g(0)=3/2$.
For $\gamma\sigma>0$, define
$$
 \widehat V_{\gamma,\sigma}^J(x)
 =-\frac{1-\gamma}{2\gamma}\log x
  -\frac{1-\sigma\gamma}{2\sigma\gamma}\log(1-x),\qquad 0<x<1.
$$
\begin{example}\label{ex:JacempLDP}
For the $\beta$-Jacobi ensemble with potential \eqref{def:vnj}, assume that $\beta n\to\infty$. 
\begin{itemize}
\item[(i)] For $\gamma\in[0,1]$ and $\gamma\sigma \in [0,1],$ the rescaled empirical measure $L_n=n^{-1}\sum_{i=1}^n\delta_{\varphi_n(\lambda_i)}$ satisfies an LDP on $\mathcal M_1(\mathbb R)$ with speed $\beta n^2$ and good rate function
\begin{equation}\label{eq:jcrdef}
 I_{V_{\gamma,\sigma}^J}(\mu)
 =\int V_{\gamma,\sigma}^J\dd\mu-\frac12\Sigma(\mu)-C_{V_{\gamma,\sigma}^J},
\end{equation}
where $\varphi_n$ and $V_{\gamma,\sigma}^J$ are defined in \eqref{defvarphi} and \eqref{eq:jpot}, respectively.
\item[(ii)] For $\gamma\in(0,1]$ and $\gamma\sigma \in (0,1]$, the original empirical measure $\widehat L_n=n^{-1}\sum_{i=1}^n\delta_{\lambda_i}$ satisfies an LDP on $\mathcal M_1([0,1])$ with speed $\beta n^2$ and good rate function $I_{\widehat V_{\gamma,\sigma}^J}$ defined by \eqref{defiv}.
Moreover,
\begin{equation}\label{eq:joremp}
 I_{\widehat V_{\gamma,\sigma}^J}(\mu)
 =I_{V_{\gamma,\sigma}^J}(\varphi_\#\mu),\qquad
 \varphi(x)=\frac{(1+\sigma)x-\sigma}{\sigma\sqrt\gamma}.
\end{equation}
\end{itemize}
\end{example}

The proofs of Examples \ref{x:lagemp} and \ref{ex:JacempLDP} are given in Section \ref{sec:gen}. 
The proof of Theorem \ref{ldpln} follows the standard strategy for Coulomb gases and $\beta$-ensembles; see, e.g., \cites{AGZ10,BAZ98,PH98}. 
The main additional difficulties are that both $\beta$ and $V_n$ may depend on $n$, and that $V_n$ need not be lower semicontinuous. 
When $\beta$ varies, standard discretization proofs of the lower bound require $\beta n/\log n\to\infty$.

We first prove the local large deviation bounds for the unnormalized measures $\overline{\mathbb P}_n := Z_n \PP,$ where $Z_n$ and $\PP$ are given in \eqref{eq:wpdf}. For $\theta\in(0,1)$,  define the measure
$$
R_{n,\theta}(\dd\boldsymbol z)
 =e^{-\theta\beta n^2H_n(\boldsymbol z)}
   \prod_{i=1}^n\omega_n(z_i) \dd\boldsymbol z,
$$
where 
$
 H_n(\boldsymbol z)
 =\frac1n\sum_{i=1}^nW_n(z_i)
   -\frac1{n^2}\sum_{1\leq i<j\leq n}\log|z_i-z_j|.
$
Its total mass is $\mathcal R_n(\theta)$ from \eqref{eq:def-frac}.
By \ref{A3}(ii), it is finite for $\theta=\theta_n$ and all sufficiently large $n$. Moreover,
\begin{equation}\label{eq:factor}
 \overline{\mathbb P}_n(\dd \bz) =e^{-(1-\theta)\beta n^2H_n(\boldsymbol z)} R_{n,\theta}(\dd \bz).
\end{equation}
The following local bounds use the balls defined in \eqref{eq:BLball}.
\begin{proposition}
\label{pr:locup} Under Assumption \ref{ass:ldp}, let $\boldsymbol\lambda=(\lambda_1,\dots,\lambda_n)\in D^n$ have density \eqref{pdf}. 
For every $\mu\in\mathcal M_1(D)$, 
$$\lim_{\varepsilon \rightarrow 0} \limsup_{n \rightarrow \infty} \frac{1}{\beta n^2} \log \overline{\mathbb P}_n\left(L_n\in B(\mu,\varepsilon)\right) \leq-\mathcal E_V(\mu),$$
where $L_n = \frac{1}{n}\sum_{i=1}^n \delta_{\lambda_i}$ and 
\begin{align}\label{defme}
\mathcal E_V(\mu)
=
\int_D V \dd\mu-\frac12\Sigma(\mu).
\end{align}
\end{proposition}
\begin{proposition}
\label{pr:loclow}
Under Assumption \ref{ass:ldp}, let $\boldsymbol\lambda=(\lambda_1,\dots,\lambda_n)\in D^n$ have density \eqref{pdf}. For every measure $\mu\in\mathcal M_1(D)$ satisfying Assumption \ref{ass:d},
$$
\lim_{\varepsilon\downarrow0}
\liminf_{n\to\infty}
\frac{1}{\beta n^2}
\log
\overline{\mathbb P}_n
\left(
L_n\in B(\mu,\varepsilon)
\right)
\geq
-\mathcal E_V(\mu).
$$
\end{proposition}

The following result reduces the problem to that on a compact set.
\begin{lemma}\label{le:tight}
    Under Assumption \ref{ass:ldp}, let $\boldsymbol\lambda=(\lambda_1,\dots,\lambda_n)\in D^n$ have density \eqref{pdf}. Then the empirical measure $L_n=\frac1n\sum_{i=1}^n\delta_{\lambda_i}$  is exponentially tight on $\mathcal M_1(D)$ with speed $\beta n^2$. Namely, there exist compact sets $\mathcal K_M\subset\mathcal M_1(D)$, $M>0$, such that
    $$
    \lim_{M\to\infty}\limsup_{n\to\infty}
    \frac1{\beta n^2}\log\mathbb P_n(L_n\notin\mathcal K_M)
    =-\infty.
    $$
\end{lemma}

The following result computes the limiting free energy.
\begin{proposition}\label{pr:logzn}
	Under Assumption \ref{ass:ldp}, let $\mu_V=\argmin_{\mu\in \mathcal M_1(D)} I_V(\mu)$. Assume further that $\mu_V$ satisfies Assumption \ref{ass:d} on $D.$ Then 
	\begin{equation}\label{eq:znlim}
\lim_{n\to\infty}
\frac1{\beta n^2}\log Z_n
=
-\inf_{\mu\in\mathcal M_1(D)}\mathcal E_V(\mu)
=
-\mathcal E_V(\mu_V),
\end{equation}
where
$$Z_n=\int_{D^n} \prod_{1 \leq i<j \leq n}\left|\lambda_i-\lambda_j\right|^\beta \exp \left\{-\beta n \sum_{i=1}^n V_n\left(\lambda_i\right)\right\} \dd \blambda.$$
\end{proposition}

\begin{proof}[Proof of Theorem \ref{ldpln}]
It was shown in \cite{AGZ10}*{Lemma 2.6.2} that $I_V$ is a good rate function. 
Combining \eqref{eq:znlim} with $\overline{\mathbb P}_n=Z_n\mathbb P_n$, by Propositions \ref{pr:locup} and \ref{pr:loclow} we obtain 
\begin{align*}
&\lim_{\varepsilon \rightarrow 0} \limsup_{n \rightarrow \infty} \frac{1}{\beta n^2} \log \mathbb P_n\left(L_n\in B(\mu,\varepsilon)\right)
=\lim_{\varepsilon \rightarrow 0} \liminf_{n \rightarrow \infty} \frac{1}{\beta n^2} \log \mathbb P_n\left(L_n\in B(\mu,\varepsilon)\right)\\
={} &
-\mathcal E_V(\mu)+\inf_{\nu\in\mathcal M_1(D)}\mathcal E_V(\nu)
=
-I_V(\mu).
\end{align*}
This shows that $\{L_n\}_{n\geq1}$ satisfies a weak LDP  with speed $\bt n^2$.
Together with the exponential tightness established in Lemma \ref{le:tight}, we obtain the full LDP.  This completes the proof.
\end{proof}

\begin{proof}[Proof of Proposition \ref{pr:locup}]
For $n\geq2$ and $x\neq y$, let
$$
 f_n(x,y)
 =\frac{n}{2(n-1)}\bigl(W_n(x)+W_n(y)\bigr)
  -\frac12\log|x-y|,
$$ 
and put $f_n(x,x)=+\infty$.  Also let
$$
 f(x,y)=\frac12V(x)+\frac12V(y)-\frac12\log|x-y|,
 \qquad f(x,x)=+\infty.
$$
Then
\begin{equation}\label{eq:hnpair}
 n^2H_n(\boldsymbol z)
 =\sum_{1\leq i\neq j\leq n}f_n(z_i,z_j)
\end{equation}
and
$
\mathcal E_V(\mu)
=\iint_{D^2} f(x,y) \mu(\dd x)\mu(\dd y).
$
For $M>0$, let
$
f_{n,M}=f_n\wedge M,
f_M=f\wedge M
$
and 
$$
J_{n,M}(\nu)
=
\iint_{D^2}
f_{n,M}(x,y) \nu(\dd x)\nu(\dd y),
\qquad
J_M(\nu)
=
\iint_{D^2}
f_M(x,y) \nu(\dd x)\nu(\dd y).
$$

We first record two consequences of Assumptions \ref{A1} and \ref{A2}.  
For every compact $K\subset D$, there exist $N_K$ and $C_K<+\infty$ such that
\begin{equation}\label{eq:wnlb}
 \inf_{n\geq N_K}\inf_{z\in K}W_n(z)\geq-C_K.
\end{equation}
Indeed, otherwise one could choose $n_j\to\infty$ and $z_j\in K$ with $W_{n_j}(z_j)\to-\infty$ and then contradict Assumption \ref{A1}(i) after extracting a convergent subsequence.  
Moreover, Assumption \ref{A2} gives the required estimate outside a sufficiently large compact set, whereas \eqref{eq:wnlb} gives it on the remaining compact set.
Thus there are constants $c,C>0$ such that, for all sufficiently large $n$ and all $z\in D$,
$$
 W_n(z)-\log(1+|z|)\geq c\log(1+|z|)-C.
$$
Combining this bound with \eqref{e:pairbd},  
we obtain
\begin{equation}\label{eq:fncoer}
 f_n(x,y)
 \geq\frac c2\log(1+|x|)+\frac c2\log(1+|y|)-C_0,
 \qquad x,y\in D
\end{equation}
for a constant $C_0$ and all sufficiently large $n$.

 We next establish the lower semicontinuity property of the varying kernels by showing that 
\begin{equation}\label{eq:vkli}
\liminf_{n\to\infty}J_{n,M}(\nu_n)
\geq
J_M(\nu).
\end{equation}
Indeed, let $x_n\to x$ and $y_n\to y$. If $x\neq y$, Assumption \ref{A1} gives
$$
\liminf_{n\to\infty}f_n(x_n,y_n)
\geq f(x,y).
$$
If $x=y$, then $-\log|x_n-y_n|\to+\infty$.
The lower bound \eqref{eq:wnlb} on the potential terms gives the same inequality.  
Therefore,
\begin{align} \label{eq:fnm}
\liminf_{n\to\infty}f_{n,M}(x_n,y_n)
\geq f_M(x,y).	
\end{align}
Equation \eqref{eq:fncoer} implies the uniform lower bound $f_{n,M}\geq-C_0$. 
Since $\nu_n\otimes\nu_n\Rightarrow\nu\otimes\nu$, Fatou's lemma for varying measures \cite{WFatou}*{Theorem 1.1} and \eqref{eq:fnm} give \eqref{eq:vkli}.

For every empirical measure $L_n$, $(L_n\otimes L_n)(\{(x,x):x\in D\})=\frac1n$, a.s. Since $f_{n,M}(x,x)=M$,
$$
\sum_{1\leq i\neq j\leq n}
f_n(\lambda_i,\lambda_j)
\geq
\sum_{1\leq i\neq j\leq n}
f_{n,M}(\lambda_i,\lambda_j)
=
n^2J_{n,M}(L_n)-nM.
$$
Using \eqref{eq:factor}, we obtain
\begin{align*}
 \overline{\mathbb P}_n(L_n\in B(\mu,\varepsilon))
 \leq{}&
 \exp\{(1-\theta_n)\beta nM\}
 \exp\left\{-(1-\theta_n)\beta n^2
       \inf_{\nu\in B(\mu,\varepsilon)}J_{n,M}(\nu)\right\}
 \mathcal R_n(\theta_n).
\end{align*}
By \ref{A3}(ii) and $\theta_n\to0$,
\begin{align}
\limsup_{n\to\infty}
\frac1{\beta n^2}
\log
\overline{\mathbb P}_n
\bigl(L_n\in B(\mu,\varepsilon)\bigr)
\leq
-
\liminf_{n\to\infty}
\inf_{\nu\in B(\mu,\varepsilon)}
J_{n,M}(\nu).
\label{eq:upeps}
\end{align}
We now claim that
\begin{equation}\label{eq:locinf}
\lim_{\varepsilon\downarrow0}
\liminf_{n\to\infty}
\inf_{\nu\in B(\mu,\varepsilon)}
J_{n,M}(\nu)
\geq
J_M(\mu).
\end{equation}
Suppose otherwise.
Then there would exist $\delta>0$, a sequence $\varepsilon_j\downarrow0$, integers $n_j\to\infty$, and $\nu_j\in B(\mu,\varepsilon_j)$ such that
$J_{n_j,M}(\nu_j) \leq J_M(\mu)-\delta.$
Since $d(\nu_j,\mu)\leq\varepsilon_j$, one has $ \nu_j\Rightarrow\mu.$
This contradicts \eqref{eq:vkli}. 

Combining \eqref{eq:upeps} and \eqref{eq:locinf}, we obtain
$$
\lim_{\varepsilon\downarrow0}
\limsup_{n\to\infty}
\frac1{\beta n^2}
\log
\overline{\mathbb P}_n
\bigl(d(L_n,\mu)\leq\varepsilon\bigr)
\leq
-J_M(\mu).
$$
Since $f_M \uparrow f$ and $f$ is bounded  below, monotone convergence yields
$ J_M(\mu) \to  \iint_{D^2} f(x,y) \mu(\dd x)\mu(\dd y),$
and letting $M\to\infty$ completes the proof.
\end{proof}

\begin{proof}[Proof of Proposition \ref{pr:loclow}]
If $\mathcal E_V(\mu)=+\infty$, then the desired lower bound is trivial. We therefore assume that
$
\mathcal E_V(\mu)<+\infty.
$

Fix $\varepsilon>0$ and $\delta>0$. By Assumption \ref{ass:d} for the target measure $\mu$ on $D$, there exists a probability measure
$\nu(\dd z)=g(z) \dd z, g\in L^\infty(D),$
such that $K:=\supp\nu$ is a compact subset of the interior of $D\cap\{V<+\infty\}$ and
\begin{equation}\label{eq:nuener}
d(\nu,\mu)<\frac{\varepsilon}{2},
\qquad \mathcal E_V(\nu)
\leq
\mathcal E_V(\mu)+\delta.
\end{equation}
Let $X_1,X_2,\ldots$ be i.i.d.~with law $\nu$ and put $L_n^X=\frac1n\sum_{i=1}^n\delta_{X_i}$.  
Since $D$ is a separable metric space, Varadarajan's theorem
\cite{Dudley02}*{Theorem 11.4.1} yields
\begin{align}\label{eq:lnx}
	L_n^X=\frac1n\sum_{i=1}^n\delta_{X_i}
\Rightarrow \nu
\qquad\text{a.s.}
\end{align}
By Assumption \ref{A1}(ii), $W_n\to V$ uniformly on the compact set $K$, so
\begin{align*}
\left|
\int_D W_n \dd L_n^X
-
\int_D V \dd L_n^X
\right|
&\leq
\sup_{x\in K}|W_n(x)-V(x)|
 \to 0.
\end{align*}
Moreover, the continuity of $V$ on $K$ implies
\begin{equation}\label{eq:wslln}
 \frac1n\sum_{i=1}^nW_n(X_i)\to\int V \dd\nu
 \qquad\text{a.s.}
\end{equation}

Let $\psi(x,y)=\log|x-y|, x,y\in D.$
We have
\begin{align}\label{eq:psiint}
\mathbb E\left|\psi(X_1,X_2)\right|
=
\iint_{D^2}
\left|\log|x-y|\right|
 \nu(\dd x)\nu(\dd y)
\leq
\|g\|_\infty^2
\iint_{K^2}
\left|\log|x-y|\right|
\dd x\dd y<+\infty.
\end{align}
Consider the $U$-statistic of degree $2$ associated with the symmetric kernel $\psi$:
$$
U_n
=
\binom{n}{2}^{-1}
\sum_{1\leq i<j\leq n}
\psi(X_i,X_j)
=
\frac{2}{n(n-1)}
\sum_{1\leq i<j\leq n}
\log|X_i-X_j|.
$$
By the strong law of large numbers for $U$-statistics \cite{Lee1990}*{Page 122 Theorem 3},  the integrability condition \eqref{eq:psiint} yields
$U_n \to\mathbb E\bigl[\psi(X_1,X_2)\bigr]=\Sigma(\nu),$ a.s.
It follows that
\begin{equation}\label{eq:islln}
\frac1{n^2}
\sum_{1\leq i<j\leq n}
\log|X_i-X_j|
 \to 
\frac12\Sigma(\nu)
\qquad\text{a.s.}
\end{equation}
Combining \eqref{eq:wslln} and \eqref{eq:islln},
\begin{equation}\label{eq:hslln}
H_n(X_1,\ldots,X_n)=\frac1n\sum_{i=1}^nW_n(X_i) -\frac1{n^2}\sum_{1\leq i<j\leq n}\log|X_i-X_j| \to  \mathcal E_V(\nu) \qquad\text{a.s.}
\end{equation}
Define
$$
A_n
=
\left\{
d(L_n^X,\nu)<\frac{\varepsilon}{2},
\quad
H_n(X_1,\ldots,X_n)
\leq
\mathcal E_V(\nu)+\delta
\right\}.
$$
Here, by abuse of notation $X_i$'s are coordinate variables and thus $A_n$ is a deterministic set.
By \eqref{eq:nuener}, on $A_n$ one has
$ d(L_n^X,\mu) \leq d(L_n^X,\nu)+d(\nu,\mu) < \varepsilon.$
Thus
$
A_n
\subset
\left\{
d(L_n^X,\mu)<\varepsilon
\right\}.
$ Put
$M_\nu=1\vee\|g\|_\infty$.  By Assumption
\ref{A3}(i), for all sufficiently large $n$,
$\omega_n\geq e^{-\de\beta n}$ almost everywhere on $K$.  Thus
\begin{align*}
 \overline{\mathbb P}_n\bigl(L_n\in B(\mu,\varepsilon)\bigr)
 &\geq
 M_\nu^{-n}e^{-\de\beta n^2}
 \int_{A_n}e^{-\beta n^2H_n(\boldsymbol x)}
          \nu^{\otimes n}(\dd\boldsymbol x)\\
 &\geq
 M_\nu^{-n}e^{-\de\beta n^2}
 e^{-\beta n^2(\mathcal E_V(\nu)+\delta)}
 \nu^{\otimes n}(A_n).
\end{align*}
Equations \eqref{eq:lnx} and \eqref{eq:hslln} imply $\nu^{\otimes n}(A_n) \to 1.$ Consequently, Assumption \ref{A0} and \eqref{eq:nuener} give
$$
\liminf_{n\to\infty}
\frac{1}{\beta n^2}
\log
\overline{\mathbb P}_n
\left(
d(L_n,\mu)\leq\varepsilon
\right)
\geq
-\mathcal E_V(\nu)-2\delta
\geq
-\mathcal E_V(\mu)-3\delta.
$$
Since $\delta>0$ is arbitrary, letting $\delta\downarrow0$ and then $\varepsilon\downarrow0$ proves the local lower bound.
\end{proof}

\begin{proof}[Proof of Lemma \ref{le:tight}]
Choose a probability measure $\nu_0$ with bounded density whose compact support is contained in the interior of $D\cap \{V<+\infty\}$. Then $\mathcal E_V(\nu_0)<+\infty$, and Proposition \ref{pr:loclow} yields
\begin{equation}\label{eq:znrlb}
 \liminf_{n\to\infty}\frac1{\beta n^2}\log Z_n \ge \lim_{\eps\downarrow 0} \liminf_{n\to\infty} \frac1{\beta n^2}\log \overline{\pz}_n(L_n\in B(\nu_0,\eps))
 \geq-\mathcal E_V(\nu_0)>-\infty.
\end{equation}
For $M>0$, let
\begin{equation}\label{eq:coerc}
 \mathcal K_M
 =\left\{\nu\in\mathcal M_1(D):
          \int_D\log(1+|x|) \nu(\dd x)\leq M\right\}.
\end{equation}
Since $\lim_{x\to\8}\log(1+|x|)=\8$, $\mathcal K_M$ is compact in $\mathcal M_1(D)$.

From \eqref{eq:fncoer} and \eqref{eq:hnpair}, for $\boldsymbol \la=(\la_1,\dots,\la_n)$, we have
$$
 H_n(\boldsymbol \la)
 \geq\left(1-\frac1n\right)
 \left(c\int_D\log(1+|x|) L_n(\dd x)-C_0\right)
$$
for all sufficiently large $n$.  Hence, by
\eqref{eq:factor},
\begin{align}\label{eq:utail}
 \overline{\mathbb P}_n(L_n\notin\mathcal K_M)
 \leq
 \exp\left\{-(1-\theta_n)\beta n^2
 \left(1-\frac1n\right)(cM-C_0)\right\}
 \mathcal R_n(\theta_n).
\end{align}
Using Assumption \ref{A3}(ii) and
\eqref{eq:znrlb},
$$
 \limsup_{n\to\infty}\frac1{\beta n^2}
 \log\mathbb P_n(L_n\notin\mathcal K_M)
 \leq-cM +C_0+\mathcal E_V(\nu_0).
$$
Letting $M\to\infty$ proves exponential tightness.
\end{proof}

\begin{proof}[Proof of Proposition \ref{pr:logzn}]
For every $M>0$,
\begin{align}\label{eq:zsplit}
 Z_n
 =\overline{\mathbb P}_n(D^n)
 =\overline{\mathbb P}_n(L_n\in\mathcal K_M)
  +\overline{\mathbb P}_n(L_n\notin\mathcal K_M),
\end{align}
where $\mathcal K_M$ is defined in
\eqref{eq:coerc}.  Fix $\delta>0$ and $0<F<\infty$.  For every
$\mu\in\mathcal K_M$, Proposition \ref{pr:locup}
provides $r_\mu>0$ such that
\begin{equation}\label{eq:uprad}
 \limsup_{n\to\infty}\frac1{\beta n^2}
 \log\overline{\mathbb P}_n
 \bigl(d(L_n,\mu)\leq r_\mu\bigr)
 \leq- (\mathcal E_V(\mu)\wedge F)+\delta.
\end{equation}
The balls $B(\mu,r_\mu)$, $\mu\in\mathcal K_M$, form an open cover
of $\mathcal K_M$.  By compactness, there exist
$\mu_1,\ldots,\mu_m\in\mathcal K_M$ such that
$
 \mathcal K_M\subset
 \bigcup_{j=1}^mB(\mu_j,r_{\mu_j}).
$
Consequently,
\begin{align*}
 \frac1{\beta n^2}
 \log\overline{\mathbb P}_n(L_n\in\mathcal K_M)
 &\leq\frac1{\beta n^2}
 \log\sum_{j=1}^m
 \overline{\mathbb P}_n
 \bigl(d(L_n,\mu_j)<r_{\mu_j}\bigr)\\
 &\leq\frac{\log m}{\beta n^2}
 +\max_{1\leq j\leq m}\frac1{\beta n^2}
 \log\overline{\mathbb P}_n
 \bigl(d(L_n,\mu_j)\leq r_{\mu_j}\bigr).
\end{align*}
Applying \eqref{eq:uprad} and Assumption \ref{A0}, letting $n\to\8$ and then $\delta\downarrow0$, $F\to +\infty$, we obtain
\begin{equation}\label{zncup}
 \limsup_{n\to\infty}\frac1{\beta n^2}
 \log\overline{\mathbb P}_n(L_n\in\mathcal K_M)
 \leq-\inf_{\nu\in\mathcal K_M}\mathcal E_V(\nu).
\end{equation}

By \eqref{eq:utail} and Assumption \ref{A3}(ii),
$$
 \limsup_{n\to\infty}\frac1{\beta n^2}
 \log\overline{\mathbb P}_n(L_n\notin\mathcal K_M)
 \leq-cM+C_0.
$$
Combining this estimate with
\eqref{eq:zsplit} and
\eqref{zncup}, and using Assumption
\ref{A0}, we obtain
\begin{align}\label{eq:znprel}
 \limsup_{n\to\infty}\frac1{\beta n^2}\log Z_n
 \leq
 \max\left\{
 -\inf_{\nu\in\mathcal K_M}\mathcal E_V(\nu),
 -cM+C_0
 \right\}.
\end{align}
By \ref{A2} and the lower bound for $V$ on compact sets, the
argument of \eqref{eq:fncoer} gives
$$
\mathcal E_V(\nu)\ge c \int \log (1+|x|)\nu(\dd x) -C_0.
$$
It follows that $\mathcal E_V(\nu)<\8$ implies $\int \log (1+|x|)\nu(\dd x)\le M$ for some large $M$.
Since the sets $\mathcal K_M$ are increasing, we find that as $M\to\8$, 
$
 \inf_{\nu\in\mathcal K_M}\mathcal E_V(\nu)
 \downarrow
 \inf_{\nu\in\mathcal M_1(D)}\mathcal E_V(\nu).
$
Letting $M\to\infty$ in \eqref{eq:znprel} gives
\begin{equation}\label{eq:znls}
 \limsup_{n\to\infty}\frac1{\beta n^2}\log Z_n
 \leq-\inf_{\mu\in\mathcal M_1(D)}\mathcal E_V(\mu).
\end{equation}

By Remark \ref{le:uniq}, $I_V$ has a unique minimizer $\mu_V$.
Since $I_V$ differs from $\mathcal E_V$ by the constant $C_V$, we have
$\inf_{\mu\in\mathcal M_1(D)}\mathcal E_V(\mu)
 =\mathcal E_V(\mu_V).$
For the converse inequality, Proposition
\ref{pr:loclow} gives
\begin{equation}\label{eq:znli}
 \liminf_{n\to\infty}\frac1{\beta n^2}\log Z_n
 \geq
 \lim_{\varepsilon\downarrow0}\liminf_{n\to\infty}
 \frac1{\beta n^2}
 \log\overline{\mathbb P}_n
 \bigl(L_n\in B(\mu_V,\varepsilon)\bigr)
 \geq-\mathcal E_V(\mu_V).
\end{equation}
Combining \eqref{eq:znls} and \eqref{eq:znli} proves
\eqref{eq:znlim}.
\end{proof}

\subsection{LDP of largest eigenvalue with general potential}\label{se:maxldp}

Early results on the largest eigenvalue LDP were proved in \cite{BDG01} for the GOE and \cite{Mai07} for the GUE. For fixed $\beta>0$ and a fixed continuous potential $V$ on $\R$, \cite{AGZ10} derives the largest eigenvalue LDP from the empirical measure LDP under additional assumptions.
For fixed $\beta>0$ and continuous potentials $V_n$ on an interval $D$ with $\sup_{x\in D}|V_n(x)-V(x)|\to0$, \cite{BG13} proves the corresponding LDP under suitable confinement and edge conditions.  Our goal is to establish a general LDP for the largest eigenvalue of $\bt$-ensembles in this section.

Define
\begin{equation}\label{eq:mxrate}
 J_V(x)=
 \begin{cases}
 V(x)- \int\log|x-y| \mu_V(\dd y)-F_V,
       &x\in D, x\geq b_V,\\
 +\infty,&\text{otherwise}.
 \end{cases}
\end{equation}
Here $\mu_V=\argmin_{\mu\in\mathcal M_1(D)} I_V(\mu)$
is the unique equilibrium measure, $b_V=\max\supp\mu_V$ and $F_V$ is characterized by
$$
 V(x)-\int\log|x-y| \mu_V(\dd y) \equiv F_V, \quad x \in \supp \mu_V.
$$

\begin{assumption}\label{ass:maximum}\rm
Suppose that
$D \subset \R $ is a closed interval with nonempty interior and that the following
conditions hold. Under the law \eqref{eq:wpdf}, let $\lambda_{\max}=\max_{i\leq n}\lambda_i$.
\begin{enumerate}[label=(C\arabic*),start=0]
\item\label{C0} $\log n=o(\beta n)$.
\item\label{C1} The set $\mathcal O$ in \ref{A1} is an interval, and $V=+\infty$ on $D\setminus\overline{\mathcal O}$. At  each finite endpoint  $x_0$ of $\overline{\mathcal O}$, assume $V(x_0)=\lim_{\mathcal O\ni x\to x_0}V(x)$, possibly $+\infty$. Moreover, $b_V\in\mathcal O$, and $J_V$ in \eqref{eq:mxrate} attains its minimum only at $b_V$.
\item\label{C2} For every compact set $K\subset D$ with $K\not\subset\mathcal O$,
\begin{equation}\label{eq:mxcomp}
 \limsup_{n\to\infty}\frac1{\beta n}
 \log\int_K e^{-\beta nW_n(x)}\omega_n(x)\dd x
 \leq-\inf_{x\in K}V(x).
\end{equation}
\item\label{C3} If $D$ is unbounded, then
\begin{equation*}
\lim_{R\to\infty}\limsup_{n\to\infty}\frac1{\beta n}
 \log\mathbb P(\lambda_{\max}>R)=-\infty.
\end{equation*}
\end{enumerate}
\end{assumption}

The following result extends the classic results of the largest eigenvalue LDP \cites{BDG01,AGZ10,BG13}. Our assumptions are weaker and cover some examples that are not known before. A key difference is that we assume only a necessary tail bound for $\la_{\max}$. The proof is given in Appendix \ref{app:maxldp}.

\begin{proposition}  \label{pr:maxldp}
Under Assumptions \ref{ass:ldp}
and \ref{ass:maximum}, let $\boldsymbol\lambda=(\lambda_1,\dots,\lambda_n)\in D^n$ have density \eqref{eq:wpdf}. Then $\lambda_{\max}:=\max_{1\leq i\leq n}\lambda_i$ satisfies an LDP on $\mathbb R$ with speed $\beta n$ and good rate function $J_V.$ In particular, for every $x>b_V$ in the interior of $D\cap\{V<+\infty\}$,
\begin{equation}\label{eq:mxrt}
 \lim_{n\to\infty}\frac1{\beta n}\log\mathbb P(\lambda_{\max}\geq x)
 =-\inf_{\substack{y\in D\\y\geq x}}J_V(y).
\end{equation}
\end{proposition}

Before going to examples, a few remarks are needed to clarify the conditions. Since we will use Gaussian $\bt$-ensemble for illustrations, let us recall the classic result with $V^H(x)=x^2/2$ and fixed $\beta>0$. A series of works \cites{BDG01,Mai07,AGZ10} prove that the largest eigenvalue satisfies an LDP with speed $\beta n$ and good rate function
\begin{equation}\label{eq:Hermax}
J_{V^H}(x)=
\begin{cases}\int_{\sqrt2}^x\sqrt{t^2-2}\dd t,&x\geq\sqrt2,\\
+\infty,&x<\sqrt2.
\end{cases}
\end{equation}
To start, we observe that the condition \ref{C3} is equivalent to exponential tightness since the empirical measure LDP has speed $\bt n^2$. For compact $K\subset\mathcal O$, \eqref{eq:mxcomp} follows from \ref{A0}, \ref{A1}(ii), and \ref{A3}(i).  Together with \ref{C2}, this gives \eqref{eq:mxcomp} for every compact $K\subset D$. 

\begin{remark}\label{re:mxendpt}
 Under Assumption \ref{ass:ldp}, suppose that
$\mathcal O$ is an interval and that the integrand
$e^{-\beta nW_n}\omega_n$ vanishes almost everywhere on every
compact subset of $D\setminus\overline{\mathcal O}$ for all
sufficiently large $n$. Then \ref{C2} holds if, for every finite
endpoint $e$ of $\overline{\mathcal O}$,
\begin{equation}\label{eq:mxend}
\lim_{\delta\downarrow0}\limsup_{n\to\infty}\frac1{\beta n}
 \log\int_{D\cap[e-\delta,e+\delta]}
 e^{-\beta nW_n(x)}\omega_n(x)\dd x\leq-V(e).
\end{equation}
Indeed, fix a compact set $K\subset D$ and let
$K_\delta=K\setminus\bigcup_{e\in K\cap\partial\mathcal O}
(e-\delta,e+\delta)$ for $\delta>0$. Then
\begin{align*}
&\limsup_{n\to\infty}\frac1{\beta n}
\log\int_K e^{-\beta nW_n(x)}\omega_n(x)\dd x\\
\leq{} &\limsup_{n\to\infty}\frac1{\beta n}
\log\Biggl[\int_{K_\delta\cap\mathcal O}
e^{-\beta nW_n(x)}\omega_n(x)\dd x
+\sum_{e\in K\cap\partial\mathcal O}\int_{D\cap[e-\delta,e+\delta]}
e^{-\beta nW_n(x)}\omega_n(x)\dd x\Biggr]\\
\leq{} &\max\left\{-\inf_{K_\delta\cap\mathcal O}V,\,
\max_{e\in K\cap\partial\mathcal O}\limsup_{n\to\infty}\frac1{\beta n}
\log\int_{D\cap[e-\delta,e+\delta]}
e^{-\beta nW_n(x)}\omega_n(x)\dd x\right\},
\end{align*}
where in the last inequality we used \ref{A0}, \ref{A1}(ii), and \ref{A3}(i) to bound the first term. The endpoint term is omitted if
$K\cap\partial\mathcal O=\varnothing$.
Letting $\delta\downarrow0$ and using \eqref{eq:mxend}
prove \ref{C2}.
\end{remark}

\begin{remark}
The condition \ref{A0} used for the empirical measure LDP is insufficient for the largest eigenvalue LDP, even in the Gaussian case. For example, take $D=\mathbb R$, $W_n(x)=V^H(x)=x^2/2$, $\omega_n\equiv1$, and $\beta n=\log\log n$ for all sufficiently large $n$. Assumption \ref{ass:ldp} holds, and every $\mu\in\mathcal M_1(\mathbb R)$ with $I_{V^H}(\mu)<\infty$ satisfies Assumption \ref{ass:d} on $\mathbb R$. Thus Theorem \ref{ldpln} gives the empirical measure LDP and convergence in probability to the semicircle law on $[-\sqrt2,\sqrt2]$.
Under the normalization \eqref{pdf}, the tridiagonal model of \cite{DE03}*{Theorem~2.12} has independent diagonal entries $d_{i,n}\sim N(0,1/(\beta n))$. By \cite{LLR83}*{Theorem~1.5.3},
$$
 \frac{\sqrt{\beta n}\max_{1\leq i\leq n}d_{i,n}}
 {\sqrt{2\log n}}\xrightarrow{\mathbb P}1.
$$
Since $\beta n=\log\log n$, the Rayleigh inequality gives
$\lambda_{\max}\geq\max_{1\leq i\leq n}d_{i,n}
\xrightarrow{\mathbb P}+\infty$.
Thus \ref{A0} alone does not ensure tightness of the maximum, which motivates the stronger condition \ref{C0}.
\end{remark}

\begin{remark}\label{rem:mxpow}
Let $D\subset\mathbb R$. Under Assumptions \ref{ass:ldp} and \ref{C0}, a sufficient condition for \ref{C3} is that there exist $\varepsilon,C,R_0>0$ such that for all sufficiently large $n$,
\begin{equation}\label{eq:mxgrowth}
 V_n(x)=W_n(x)-\frac1{\beta n}\log\omega_n(x)
 \geq(1+\varepsilon)\log(1+x)-C, \qquad \text{a.e. }x\in D,\ x\geq R_0.
\end{equation}
Indeed, \eqref{eq:mxwtail} together with
\eqref{eq:fulldeco} and \eqref{eq:mxgrowth} gives
for $R\geq R_0$ and all sufficiently large $n$,
\begin{align*}
 \mathbb P(\lambda_{\max}>R)
 \leq n e^{(K+C)\beta n+o(\beta n)}
 \int_R^\infty(1+x)^{-\varepsilon\beta n}\dd x
 =\frac{n e^{(K+C)\beta n+o(\beta n)}
 (1+R)^{1-\varepsilon\beta n}}{\varepsilon\beta n-1}.
\end{align*}
By \ref{C0}, this bound implies \ref{C3}.

In view of Assumptions \ref{A2} and \ref{A3}(i), the growth condition \eqref{eq:mxgrowth} seems optimal for exponential tightness. However, since we only assume lower semicontinuity for $V_n$,  \eqref{eq:mxgrowth} is not necessary for
\ref{C3}. For example, take $D=\mathbb R$, $\beta=1$,
$W_n(x)=V(x)=x^2/2$, and
$$
 I_n=[n,n+e^{-n^4}],\qquad
 \omega_n(x)=
 \begin{cases}
 e^{nx^2/2},&x\in I_n,\\
 1,&x\notin I_n.
 \end{cases}
$$
Then $V_n=0$ on $I_n$ and $V_n(x)=x^2/2$ elsewhere, hence \eqref{eq:mxgrowth} fails. But the LDP still holds with speed $n$ and good rate function $J_{V^H}$ given in \eqref{eq:Hermax}. Indeed, Conditions \ref{A0}--\ref{A3}(i) and \ref{C0}--\ref{C2} are immediate.
For \ref{A3}(ii), take $\theta_n=n^{-1/2}$. 
Applying \eqref{e:pairbd} and \eqref{eq:Frelb} with $\alpha=2$ and $\varepsilon=1/4$, we obtain
\begin{align*}
 \mathcal R_n(\theta_n)
 \leq e^{cn^{3/2}}
 \left[\int_{\mathbb R}e^{-\sqrt n x^2/4}\omega_n(x)\dd x\right]^n
 \leq e^{cn^{3/2}}
 \left(2\sqrt\pi\,n^{-1/4}
 +e^{-n^4+n(n+1)^2/2}\right)^n,
\end{align*}
and hence $\limsup_n n^{-2}\log\mathcal R_n(\theta_n)\leq0$.
For $R$ large enough, we split the integral in \eqref{eq:mxwtail} over $I_n$ and its complement. Using  $\log(1+x)\leq x^2/4$ for $x\geq4$, we obtain
$$
 \mathbb P(\lambda_{\max}>R)
 \leq ne^{Kn+o(n)}\left[
 \frac{2}{nR}e^{-nR^2/4}+e^{-n^4}(n+2)^{n-1}\right].
$$
Thus \ref{C3} holds.

In the classic theory of largest eigenvalue LDP as in \cites{AGZ10,BG13}, the growth condition \eqref{eq:mxgrowth} is the standard assumption that is inherited from the empirical measure LDP. In the above example, we see that the largest eigenvalue LDP could hold with a growth condition weaker than \eqref{eq:mxgrowth}.
\end{remark}

\begin{remark} \label{re:lamindiv}
 The largest eigenvalue LDP in Proposition \ref{pr:maxldp}
can hold even when $\lambda_{\min}\to-\infty$ in probability,
so that $\max_i|\lambda_i|$ is not tight.
In contrast, the assumptions in \cites{AGZ10,BG13} ensure
tightness of $\max_i|\lambda_i|$.
Consider \eqref{eq:wpdf} with $D=\mathbb R$, $\beta=2$,
$W_n(x)=V(x)=V^H(x)=x^2/2$, and
\begin{align*}
	 \omega_n(x)=
 \begin{cases}
 e^{n^3+2n^{3/2}},&x\in I_n:=[-n^{1/4}-e^{-n^3},-n^{1/4}],\\
 1,&x\notin I_n.
 \end{cases}
\end{align*}

We first prove that $\lambda_{\min}\to-\infty$ in probability. Since $\omega_n=1$ outside $I_n$, we obtain
$$\mathbb P(\lambda_i\notin I_n\text{ for all }i)
=\frac1{Z_n}\int_{(\mathbb R\setminus I_n)^n}
e^{-n\sum_{i=1}^n x_i^2}
\prod_{1\leq i<j\leq n}|x_i-x_j|^2\dd\boldsymbol x
\leq\frac{Z_{n,n}^{G}}{Z_n},$$
where 
$$
Z_{k,n}^{G}
=\int_{\mathbb R^k}e^{-n\sum_{i=1}^k x_i^2}
\prod_{1\leq i<j\leq k}|x_i-x_j|^2\dd\boldsymbol x, \quad k\geq1.
$$
To bound $Z_n$ from below, first note that
\begin{align*}
\int_{I_n}e^{-nx^2}\omega_n(x)\dd x
&=e^{n^{3/2}}\int_0^1
e^{-2n^{5/4}e^{-n^3}t-ne^{-2n^3}t^2}\dd t
=e^{n^{3/2}+o(1)}.
\end{align*}
Since $\omega_n=1$ on $[-2,2]$ and
$\operatorname{dist}(I_n,[-2,2])>1$ for large $n$, we obtain
\begin{align*}
Z_n
&=\int_{\mathbb R^n}e^{-n\sum_{i=1}^n x_i^2}
\prod_{1\leq i<j\leq n}|x_i-x_j|^2
\prod_{i=1}^n\omega_n(x_i)\dd\boldsymbol x\\
&\geq n\left(\int_{I_n}e^{-nx^2}\omega_n(x)\dd x\right)
\left(\int_{[-2,2]^{n-1}}e^{-n\sum_{i=1}^{n-1}x_i^2}
\prod_{1\leq i<j\leq n-1}|x_i-x_j|^2\dd\boldsymbol x\right).
\end{align*}
Let $(\eta_1,\ldots,\eta_{n-1})$ be the $(n-1)$-particle Hermite ensemble with $\beta=2$ and potential $V^H$.
By \eqref{eq:Hermax},
$$
\limsup_{n\to\infty}\frac1{2(n-1)}
\log\mathbb P\!\left(\max_{i\leq n-1}\eta_i>2\right)
\leq-J_{V^H}(2)<0.
$$
For $0<c<J_{V^H}(2)$ and all sufficiently large $n$, it follows that
\begin{align*}
\mathbb P\!\left(\max_{i\leq n-1}|\eta_i|>
2\sqrt{\frac n{n-1}}\right)
\leq\mathbb P\!\left(\max_{i\leq n-1}|\eta_i|>2\right)
\leq2\mathbb P\!\left(\max_{i\leq n-1}\eta_i>2\right)
\leq2e^{-2c(n-1)}.
\end{align*}
Therefore, the change of variables $x_i=\sqrt{(n-1)/n}\,\eta_i$ yields
\begin{align*}
\int_{[-2,2]^{n-1}}e^{-n\sum_{i=1}^{n-1}x_i^2}
\prod_{1\leq i<j\leq n-1}|x_i-x_j|^2\dd\boldsymbol x={}&Z_{n-1,n}^{G}
\mathbb P\left(\max_{i\leq n-1}|\eta_i|
\leq2\sqrt{\frac n{n-1}}\right)\\
\geq{} & \bigl(1-2e^{-2c(n-1)}\bigr)Z_{n-1,n}^{G}
\geq\frac12Z_{n-1,n}^{G}.
\end{align*}
Using the Selberg integral formula (see, e.g., \cite{AGZ10}*{Corollary 2.5.9}) and \eqref{e:factbd}, we find
\begin{align*}
\mathbb P(\lambda_i\notin I_n\text{ for all }i)
\leq\frac2n\frac{Z_{n,n}^{G}}{Z_{n-1,n}^{G}}
e^{-n^{3/2}+o(1)}
=\frac{2\sqrt\pi\,n!}{2^{n-1}n^{n+1/2}}
e^{-n^{3/2}+o(1)}
=e^{-n^{3/2}+O(n)}\to  0,
\end{align*}
For every fixed $R>0$, we have $I_n\subset(-\infty,-R)$
for all sufficiently large $n$, and hence
$$
\mathbb P(\lambda_{\min}\geq-R)
\leq\mathbb P(\lambda_i\notin I_n\text{ for all }i)
\to 0.
$$
This proves $\lambda_{\min}\to-\infty$ in probability.

We next verify the assumptions for the LDPs of empirical measure and largest eigenvalue with good rate functions $I_{V^H}$ and $J_{V^H}$, respectively.
Conditions \ref{A0}--\ref{A3}(i) and \ref{C0}--\ref{C2} are immediate.
For \ref{A3}(ii), take $\theta_n=n^{-1/2}$.
Using $|I_n|=e^{-n^3}$, \eqref{e:pairbd}, and
\eqref{eq:Frelb}, we deduce 
\begin{align*}
\mathcal R_n(\theta_n)
&=\sum_{m=0}^n\binom nm
e^{m(n^3+2n^{3/2})}
\int_{I_n^m\times(\mathbb R\setminus I_n)^{n-m}}
e^{-\sqrt n\sum_{i=1}^n x_i^2}
\prod_{1\leq i<j\leq n}|x_i-x_j|^{2/\sqrt n}\dd\boldsymbol x\\
&\leq e^{Cn^{3/2}}\sum_{m=0}^n\binom nm
|I_n|^m e^{m(n^3+2n^{3/2})}
|I_n|^{m(m-1)/\sqrt n}
(2+n^{1/4})^{2m\sqrt n}\\
&=e^{Cn^{3/2}}\sum_{m=0}^n\binom nm
\exp\left\{
2mn^{3/2}+2m\sqrt n\log(2+n^{1/4})
-n^{5/2}m(m-1)
\right\}\\
&\leq e^{C'n^{3/2}},
\end{align*}
where the exponent inside the sum is nonpositive for $m\geq2$ and all sufficiently large $n$. Thus \ref{A3}(ii) holds. 
Lemma \ref{le:symreg} verifies Assumption \ref{ass:d} for every $\mu\in\mathcal M_1(\mathbb R)$ with $I_{V^H}(\mu)<\infty$.
Theorem \ref{ldpln} therefore gives the empirical measure LDP with speed $2n^2$ and rate $I_{V^H}$. 
Since $\omega_n(x)=1$ for $x>0$, we have $V_n(x)=x^2/2\geq2\log(1+x)-2$ for $x\geq1$.
Thus \eqref{eq:mxgrowth} holds, and Remark \ref{rem:mxpow} gives \ref{C3}. Proposition \ref{pr:maxldp} yields the LDP for $\lambda_{\max}$ with speed $2n$ and rate $J_{V^H}$ in \eqref{eq:Hermax}.
\end{remark}

\begin{remark}
The proof of Proposition \ref{pr:maxldp} uses the endpoint condition in \ref{C1} (i.e.,
$
V(x_0)=\lim_{\mathcal O\ni x\to x_0}V(x),$  $x_0\in\partial\mathcal O$)
to extend the local lower bound to the finite right endpoint.
Lower semicontinuity of $V$ alone is not enough. To see this, we fix $\beta>0$, $D=[-2,2]$, and choose
$0<c<J_{V^H}(2)$, where $J_{V^H}$ is given in \eqref{eq:Hermax}.
Define
$$
\widetilde V^H(x)=
\begin{cases}
V^H(x),&-2\leq x<2,\\
V^H(2)-c,&x=2.
\end{cases}
$$
The potential $\widetilde V^H$ is lower semicontinuous and continuous on
$\mathcal O=(-2,2)$, so \ref{A1} holds.
\ref{A0} and \ref{A2} are trivial. Remark \ref{r:a3unit} gives \ref{A3}.
Every $\mu\in\mathcal M_1(D)$ of finite
logarithmic energy satisfies $\mu(\{2\})=0$.
Thus $\int_D\widetilde V^H\dd\mu=\int_D V^H\dd\mu$, and the
two energy functionals have the same unique minimizer.
Together with \eqref{eq:calsig}, we have
$$
\mu_{\widetilde V^H}(\dd x)=\rho^H(x)\dd x,\qquad
F_{\widetilde V^H}=\frac{1+\log2}{2}.
$$
Thus the function defined by \eqref{eq:mxrate} satisfies
$$
J_{\widetilde V^H}(2)
=V^H(2)-c-\int_{\mathbb R}\log|2-y|\rho^H(y)\dd y-F_{\widetilde V^H}
=J_{V^H}(2)-c.
$$
Let $\mathbb P_n$ and $\widetilde{\mathbb P}_n$ denote the laws on $D^n$ defined by \eqref{pdf} with potentials $V^H$ and
$\widetilde V^H$, respectively.
Proposition \ref{pr:maxldp} applied to $V^H$ on $D$ gives an LDP for $\lambda_{\max}$ under $\mathbb P_n$, with rate $J_{V^H}$ on $[\sqrt2,2]$ and $+\infty$ elsewhere.
The two densities are equal almost everywhere, so $\widetilde{\mathbb P}_n=\mathbb P_n$ for every $n$. Hence
\begin{align*}
\lim_{\delta\downarrow0}\lim_{n\to\infty}\frac1{\beta n}
\log\widetilde{\mathbb P}_n(|\lambda_{\max}-2|<\delta)
&=\lim_{\delta\downarrow0}\lim_{n\to\infty}\frac1{\beta n}
\log\mathbb P_n(|\lambda_{\max}-2|<\delta)
=-J_{V^H}(2)\neq-J_{\widetilde V^H}(2).
\end{align*}
All the other conditions in Assumption \ref{ass:maximum} hold. Hence the endpoint condition in \ref{C1} cannot be omitted while keeping the rate function \eqref{eq:mxrate}.
\end{remark}

\begin{remark}\label{rem:mxunit}
If $\omega_n\equiv1$, then \ref{C0}, \ref{A1}, and \ref{A2} imply \ref{C2} and \ref{C3}. Indeed, for every compact $K\subset D$, a moment of reflection shows that $\liminf_n\inf_KW_n\geq\inf_KV$ by \ref{A1} and compactness. Using \ref{A0}, we find
\begin{align*}
 \limsup_{n\to\infty}\frac1{\beta n}
 \log\int_K e^{-\beta nW_n(x)}\dd x
 \leq\limsup_{n\to\infty}
 \left[\frac{\log(1+|K|)}{\beta n}-\inf_KW_n\right]\leq-\inf_KV.
\end{align*}
This proves \ref{C2}.
Assumption \ref{A2} also implies the growth bound \eqref{eq:mxgrowth} in Remark \ref{rem:mxpow}, which gives \ref{C3}.  Since uniform convergence of $W_n$ on compact subsets of $D$ implies \ref{A1}, together with Remark \ref{r:a3unit}, this shows that our framework covers the corresponding settings in \cites{Joh00,BG13}. In particular, the classic result \cite{AGZ10}*{Theorem 2.6.6} is recovered.
\end{remark}
   
Let us turn to examples. For $\bt$-Hermite ensemble with varying $\beta$, \cite{Pak20}*{Theorem~1.1} establishes this LDP under $\log n/n\ll\beta\ll1/\log n$. Proposition \ref{pr:maxldp} requires only $\log n/n\ll\beta$. Indeed, for $x\geq\sqrt2$,
$$
J_{V^H}(x) = \frac{x^2}{2}
-\int_{\mathbb R}\log|x-y|\rho^H(y)\dd y
-\frac{1+\log2}{2},
$$
where $\rho^H$ is given in \eqref{eq:semic}.
The Stieltjes transform of $\rho^H$ gives
$J'_{V^H}(x)=x-\int_{\mathbb R}\frac{\rho^H(y)}{x-y}\dd y=\sqrt{x^2-2}>0$ for $x>\sqrt2$.
Together with $b_{V^H}=\sqrt2$ and the continuity of $V^H$, this verifies \ref{C1}.
The proof of Example \ref{ex:her2} verifies Assumption \ref{ass:ldp} and Assumption \ref{ass:d} on $\mathbb R$ for every $\mu\in\mathcal M_1(\mathbb R)$ with $I_{V^H}(\mu)<\infty$.
Condition \ref{C0} holds by assumption.
Since $\omega_n\equiv1$, Remark \ref{rem:mxunit} gives \ref{C2} and \ref{C3}.

For the $\beta$-Laguerre ensemble with $\beta=2$, \cite{Joh00}*{Theorem~1.6} considered large deviations of the largest eigenvalue for $p=\lfloor n/\gamma\rfloor$ for $0<\gamma\leq1$ without giving explicit rate functions except for the case $p=n$.  
For $0<\gamma<1$ and fixed $\bt>0$,  the physics paper \cite{KPC10} gave the large deviation rate functions for both of the upper and lower tails of the smallest eigenvalue.
For fixed $\beta>0$, \cite{GR11}*{Lemma~5.1} and \cite{Fer08}*{Section~5.2} discussed the LDPs for both extreme eigenvalues and for the largest eigenvalue, respectively, based on the empirical measure LDP in \cite{HP98}, which does not include the critical case $\gamma=1$ as discussed after \eqref{eq:hp98}.  
Define
\begin{equation}\label{eq:Lagint}
J_{\widehat V_\gamma^L}(x)=
\begin{cases}
 \frac1{2\gamma}\int_{{\gamma_+}}^x
\frac{\sqrt{(t-\gamma_-)(t-\gamma_+)}}t\dd t,&x\geq{\gamma_+},\\
+\infty,&x<{\gamma_+},
\end{cases}
\end{equation}
and 
\begin{equation}\label{eq:lcmax}
J_{V_\gamma^L}(x)=
\begin{cases}
J_{\widehat V_\gamma^L}(1+\sqrt\gamma x),&0<\gamma\leq1,\\
J_{V^H}(x/\sqrt2),&\gamma=0.
\end{cases}
\end{equation}
where $\gamma_\pm=(1\pm\sqrt\gamma)^2$ and $J_{V^H}$ is given in \eqref{eq:Hermax}.
For the smallest eigenvalue, define
\begin{equation}\label{eq:lomin}
J_{\widehat V_\gamma^L}^{-}(x)=
\begin{cases}
\displaystyle\frac1{2\gamma}\int_x^{\gamma_-}
\frac{\sqrt{(\gamma_--t)(\gamma_+-t)}}{t}\dd t,
&0<x\leq\gamma_-,\\
+\infty,&\text{otherwise}.
\end{cases}
\end{equation}

\begin{example}\label{ex:lag-max}
For the $\beta$-Laguerre ensemble in Example \ref{x:lagemp}, assume that $\beta n/\log n\to\infty$. 
\begin{itemize}
    \item[(i)] For $\gamma \in[0,1] $,  $\sqrt{p/n}(\lambda_{\max}-1)$ satisfies an LDP on $\mathbb R$ with speed $\beta n$ and good rate function $J_{V_\gamma^L}$.
    \item[(ii)] For $\gamma\in(0,1]$, $\lambda_{\max}$ satisfies an LDP on $\mathbb R$ with speed $\beta n$ and good rate function $J_{\widehat V_\gamma^L}$.
    \item[(iii)]   For $\gamma\in(0,1)$, $\lambda_{\min}$ satisfies an LDP on $\mathbb R$ with speed $\beta n$ and good rate function $J_{\widehat V_\gamma^L}^{-}$.
\end{itemize}
\end{example}
When $\gamma=1$, the left endpoint of $\supp\rho_1^L$ in \eqref{eq:hgam} equals the hard boundary $0$. The LDP for $\lambda_{\min}$ at speed $\beta n$ is  then degenerate.
\begin{remark}
For the $\beta$-Laguerre ensemble with fixed $\beta>0$ and $n/p\to\gamma=0$, \cite{JL15}*{Theorem 2} proves an LDP for $\lambda_{\max}$ at speed $p$ in the normalization \eqref{def:vnl}. In particular, for $x>1$,
$$
\lim_{n\to\infty}\frac1{\beta p}
\log\mathbb P(\lambda_{\max}\geq x)
=-\frac{x-1-\log x}{2}.
$$
Example \ref{ex:lag-max} gives
$$
\lim_{n\to\infty}\frac1{\beta n}
\log\mathbb P \left(\lambda_{\max}\geq1+\sqrt{n/p} y\right)
=-\frac12\int_2^y\sqrt{t^2-4}\dd t, \qquad y>2.
$$
Thus the former concerns large deviations at fixed $x>1$, and the latter moderate deviations at $1+\sqrt{n/p} y$ for fixed $y>2$.
\end{remark}

For the $\beta$-Jacobi ensemble with fixed $\beta>0$, \cite{GR11}*{Lemma~6.1} states the largest eigenvalue LDP when $\gamma\in(0,1]$ and $\gamma\sigma\in(0,1)$, and the smallest eigenvalue LDP when $\gamma\in(0,1)$ and $\gamma\sigma\in(0,1]$. 
Under $\beta n\gg\log n$, \cite{MW24}*{Theorems~1 and 2} proves LDPs for the rescaled extreme eigenvalues with $\gamma\in[0,1)$ and $\gamma\sigma\in[0,1)$, relying on \cite{Ma23J}*{Theorem~1.7}.
As noted above, these works alone do not include the critical case $\gamma=1$ for general $\beta$. Define
\begin{equation}\label{eq:jomax}
J_{\widehat V_{\gamma,\sigma}^J}(x)=
\begin{cases}
\frac{1+\sigma}{2\sigma\gamma}
\int_{u_+}^x\frac{\sqrt{(t-u_-)(t-u_+)}}{t(1-t)}\dd t,
&u_+\leq x \leq 1,\\
+\infty,&\text{otherwise},
\end{cases}
\end{equation}
where  $u_\pm$ are given in \eqref{eq:hgamsig} and
\begin{equation}\label{eq:jcmax}
J_{V_{\gamma,\sigma}^J}(y)=
\begin{cases}
J_{\widehat V_{\gamma,\sigma}^J}
\left(\frac{\sigma(1+\sqrt\gamma y)}{1+\sigma}\right),
&\gamma\sigma>0,\\
J_{V_\gamma^L}(y),&\sigma=0<\gamma,\\
J_{V^H}\left(\sqrt{\frac{1+\sigma}{2}} y\right),
&\gamma=0,
\end{cases}
\end{equation}
Here $J_{V^H}$ and $J_{V_\gamma^L}$ are given in \eqref{eq:Hermax} and \eqref{eq:lcmax}, respectively. 

\begin{example}\label{ex:jac-max}
Under the assumptions of Example \ref{ex:JacempLDP}, suppose further that $\beta n/\log n\to\infty$.
\begin{itemize}
\item[(i)] For every $\gamma\in[0,1]$ and $\gamma\sigma\in[0,1)$, $((p_1+p_2)\lambda_{\max}-p_1)/\sqrt{np_1}$ satisfies an LDP on $\mathbb R$ with speed $\beta n$ and good rate function $J_{V_{\gamma,\sigma}^J}.$
\item[(ii)] For $\gamma \in (0,1] $ and $\gamma\sigma\in (0,1)$, $\lambda_{\max}$ satisfies an LDP on $\mathbb R$ with speed $\beta n$ and good rate function $J_{\widehat V_{\gamma,\sigma}^J}.$
\item[(iii)]   For $\gamma\in(0,1)$ and $\gamma\sigma\in(0,1]$, $\lambda_{\min}$ satisfies an LDP on $\mathbb R$ with speed $\beta n$ and good rate function $J^-_{\widehat V_{\ga,\si}^J } (x):= J_{\widehat V_{\gamma\sigma,\sigma^{-1}}^J}(1-x)$.
\end{itemize}
\end{example}

When $\gamma=1$, the LDP for $\lambda_{\min}$ at speed $\beta n$ is degenerate, with rate function equal to $0$ at $0$ and $+\infty$ elsewhere.
When $\gamma\sigma=1$, the LDPs for the original and rescaled largest eigenvalues at the same speed are also degenerate.
Their rate functions equal $0$ at $1$ and $\sqrt\gamma$, respectively, and $+\infty$ elsewhere.
The proofs of Examples \ref{ex:lag-max} and \ref{ex:jac-max} are given in Section \ref{sec:gen}. Clearly, the argument of Example \ref{ex:jac-max} also works for the smallest eigenvalue of Laguerre and Jacobi ensembles with trivial rate function when $\ga=1$.
The proof of Proposition \ref{pr:maxldp} is given in Appendix \ref{app:maxldp}.

\subsection{Lower tail probability of largest eigenvalue with general potential}\label{se:maxlow}

We start with the large gap probability on the real line. 

\begin{proposition} \label{pr:lowgap}
Under Assumption \ref{ass:ldp}, let $D\subset\mathbb R$ and $\boldsymbol\lambda=(\lambda_1,\dots,\lambda_n)\in D^n$ have density \eqref{pdf}. 
Let $\mu_V=\argmin_{\mu \in \mathcal{M}_1(D)} I_V(\mu)$ and, for $\zeta_1<\zeta_2,$ let $\mu_{\zeta_1,\zeta_2}=\argmin_{\mu \in \mathcal{M}_1(D\cap[\zeta_1,\zeta_2])} I_V(\mu)$. Assume further that $\mu_V$ and $\mu_{\zeta_1,\zeta_2}$ satisfy Assumption \ref{ass:d} with $D$ and $D\cap[\zeta_1,\zeta_2]$, respectively. Then 
\begin{align*}
    \lim_{n \rightarrow \infty}
    \frac{1}{\beta n^2}
    \log \mathbb{P}\bigl(\zeta_1\leq \lambda_i \leq \zeta_2, 
    1\leq i \leq n\bigr)
    =
    -
    I_V(\mu_{\zeta_1,\zeta_2}).
\end{align*}
Moreover,
\begin{align*}
    I_V(\mu_{\zeta_1,\zeta_2})
    =
    \frac12\int V \dd\mu_{\zeta_1,\zeta_2}
    +
    \frac12 F_{\zeta_1,\zeta_2}
    -
    C_V,
\end{align*}
where the constant $C_V$ is defined in \eqref{defiv}, and $F_{\zeta_1,\zeta_2}\in\mathbb R$ is characterized by
\begin{align}\label{eq:qelov}
-\int \log|x-y| \mu_{\zeta_1,\zeta_2}(\dd y)
+
V(x)
=
F_{\zeta_1,\zeta_2},
\qquad
\text{q.e. } x\in\supp(\mu_{\zeta_1,\zeta_2}).
\end{align}
Moreover, if $\supp(\mu_{\zeta_1,\zeta_2})=[s,t]$ with $s<t$ and $V$ is finite and continuous on $[s,t]$, then \eqref{eq:qelov} holds pointwise:
\begin{align}\label{eq:fzeta}
    -\int \log|x-y|\mu_{\zeta_1,\zeta_2}(\dd y)
    +V(x)
    =F_{\zeta_1,\zeta_2},
    \qquad
    x\in[s,t].
\end{align}
\end{proposition}

\begin{proof}
Let $K_{\zeta_1,\zeta_2}=D\cap[\zeta_1,\zeta_2].$
Since $\beta$ and $\{V_n\}_{n\geq1}$ satisfy Assumption \ref{ass:ldp} on $D$, so do their restrictions to $K_{\zeta_1,\zeta_2}$. By the definition in \eqref{pdf}, we have the confinement principle
$$
\mathbb P\bigl(\zeta_1\le \lambda_i\le \zeta_2, 1\le i\le n\bigr) = \frac{Z_n^K}{Z_n},
$$
where 
$$
Z_n=\int_{D^n} \prod_{1 \leq i<j \leq n}\left|\lambda_i-\lambda_j\right|^\beta \exp \left\{-\beta n \sum_{i=1}^n V_n\left(\lambda_i\right)\right\} \dd \blambda,
$$
$$Z_n^K=\int_{K_{\zeta_1,\zeta_2}^n} \prod_{1 \leq i<j \leq n}\left|\lambda_i-\lambda_j\right|^\beta \exp \left\{-\beta n \sum_{i=1}^n V_n\left(\lambda_i\right)\right\} \dd \blambda.$$
Since $\mu_V$ and $\mu_{\zeta_1,\zeta_2}$ satisfy Assumption \ref{ass:d} with $D$ and $D\cap[\zeta_1,\zeta_2]$, respectively, Proposition \ref{pr:logzn} yields
\begin{align*}
     \lim_{n\to\infty}
\frac{1}{\beta n^2}
\log\mathbb P\bigl(\zeta_1\le \lambda_i\le \zeta_2, 1\le i\le n\bigr) ={}& -\inf_{\mu \in \mathcal{M}_1(K_{\zeta_1,\zeta_2})} \mathcal E_V(\mu) +\inf _{\mu \in \mathcal{M}_1(D)} \mathcal{E}_V(\mu) \\
={}& -\inf _{\mu \in \mathcal{M}_1(K_{\zeta_1,\zeta_2})} I_V(\mu) = -I_V(\mu_{\zeta_1,\zeta_2}).
\end{align*}
To evaluate $I_V(\mu_{\zeta_1,\zeta_2})$, recall that $\mu_{\zeta_1,\zeta_2}$ is the equilibrium measure on $K_{\zeta_1,\zeta_2}$. 
By Lemma \ref{lemr}, there exists a constant $F_{\zeta_1,\zeta_2}\in\mathbb R$ such that \eqref{eq:qelov} holds. 
Applying \eqref{eq:calsig} on $K_{\zeta_1,\zeta_2}$ and using \eqref{defiv}, we obtain
\begin{align*}
I_V(\mu_{\zeta_1,\zeta_2})
&=\int V\dd\mu_{\zeta_1,\zeta_2}
-\frac12\left(\int V\dd\mu_{\zeta_1,\zeta_2}
-F_{\zeta_1,\zeta_2}\right)-C_V
=\frac12\int V\dd\mu_{\zeta_1,\zeta_2}
+\frac12F_{\zeta_1,\zeta_2}-C_V.
\end{align*}
Under the additional assumptions, Remark \ref{rem:qepw} shows that \eqref{eq:fzeta} holds for every $x\in[s,t]$.
\end{proof}

Based on this result, Theorem \ref{thm:jaclow} unifies the lower tail probabilities for all the three classical ensembles, while we still need to compute the rate functions case by case in the proof. For clarity and for comparison with the physics literature, we list the Hermite and Laguerre cases as examples below, whose proofs will be provided in the next section.
The two-sided constrained problem for the $\beta$-Hermite ensemble was studied in the physics literature in \cites{DM06,DM08}; see also \cite{MS14}.

\begin{example} \label{ex:her2}
Consider the $\beta$-Hermite ensemble on $\mathbb R$, namely \eqref{pdf} with potential $V_H(x)=x^2/2$. Assume that $\beta n \to \infty$. For $-\8\le \zeta_1<\zeta_2\le \8$,
$$
\lim_{n\to\infty}
\frac{1}{\beta n^2}
\log \mathbb P\bigl(\zeta_1\le \lambda_i \le \zeta_2,  1\leq i \leq n \bigr)
=
-I_{V_H}(\mu_{\zeta_1,\zeta_2}^H),
$$
where
\begin{align}\label{eq:hcon}
\mu_{\zeta_1,\zeta_2}^H = \argmin_{\nu\in\mathcal M_1([\zeta_1,\zeta_2])} I_{V_H}(\nu).
\end{align}
The support of $\mu_{\zeta_1,\zeta_2}^H$ is an interval $[s,t]\subset[\zeta_1,\zeta_2]$ given by
$$
[s,t]
=
\begin{cases}
[-\sqrt2,\sqrt2],
& \zeta_1\le -\sqrt2, \zeta_2\ge \sqrt2,\\
[\zeta_1,b(\zeta_1)],
& \zeta_1>-\sqrt2, \zeta_2\ge b(\zeta_1),\\
[a(\zeta_2),\zeta_2],
& \zeta_2<\sqrt2, \zeta_1\le a(\zeta_2),\\
[\zeta_1,\zeta_2],
& \text{otherwise},
\end{cases}
$$
where  $a(\zeta)=(\zeta-2\sqrt{\zeta^2+6})/3$ and  $b(\zeta)=(\zeta+2\sqrt{\zeta^2+6})/3.$
Let $m=(s+t)/2,$ and  $h=(t-s)/2.$ Then $\mu_{\zeta_1,\zeta_2}^H$ has density
$$
    \rho_{\zeta_1,\zeta_2}^H(x)
    =
    \frac{
        1-x^2+mx+\frac{h^2}{2}
    }{
        \pi\sqrt{(t-x)(x-s)}
    }
    \mathbf 1_{[s,t]}(x)
$$
and
$$
I_{V_H}(\mu_{\zeta_1,\zeta_2}^H)
=
\frac{m^2}{2}
+
\frac{h^2}{4}
-
\frac{m^2h^2}{4}
-
\frac{h^4}{32}
-
\frac12\log\frac{h}{2}
-
\frac38
-
\frac14\log2.
$$
In particular, 
\begin{align}
 I_{V_H}(\mu_{-\infty,\zeta}^H)
 &=\frac{-\zeta^4-\zeta^3\sqrt{\zeta^2+6}+36\zeta^2
       -15\zeta\sqrt{\zeta^2+6}}{108}
   -\frac12\log\frac{\zeta+\sqrt{\zeta^2+6}}{3\sqrt2}, \quad
 \zeta<\sqrt2, \notag \\ 
 I_{V_H}(\mu_{-\zeta,\zeta}^H)
 &=-\frac{\zeta^4}{32}+\frac{\zeta^2}{4}
   -\frac12\log\frac{\zeta}{\sqrt2}-\frac38, \quad 0<\zeta<\sqrt2.  \label{eq:ivhsym}
\end{align}
\end{example}

\begin{remark}\label{re:her2sd}
The optimization leading to \eqref{eq:ivhsym} is a particular case of \cite{WZschatten}*{Lemma 3.6} and has motivated our study of the Schatten ball volume problems (see also Remark 3.4 in the same paper). Indeed, letting $s=1/2$ and $p=2$ there, we have $b_2=2$ and $a_0=\sqrt2$. Hence, if $0<\zeta<\sqrt2$, the minimizing measure $\mu_{-\zeta,\zeta}^H$ is supported on $[-\zeta,\zeta]$ and has density
$$
\si_{\zeta}(x)
=
\frac{1+\zeta^2/2-x^2}{\pi\sqrt{\zeta^2-x^2}},
\qquad |x|<\zeta.
$$
Moreover, $\Sigma(\mu_{-\zeta,\zeta}^H)=\log\frac{\zeta}{2}-\frac{\zeta^4}{16}$ and $\int_{-\zeta}^\zeta |x|^2 \mu_{-\zeta,\zeta}^H(\dd x)=\frac{\zeta^2}{2}-\frac{\zeta^4}{8}$. Then \eqref{eq:ivhsym} follows from $C_{V_H}=\frac14\log2 +\frac38$. 
\end{remark}

For the $\beta$-Laguerre ensemble, \cite{VMB07} derived the largest eigenvalue left tail rate for $0<\gamma\leq1$ and \cite{KPC10} derived both left and right tail rates for the smallest eigenvalue in the physics literature.  We give a unified interval confinement rate that also covers the rescaled regime $\gamma=0$. The result holds for varying $\beta$ under $\beta n\to\infty$. 

\begin{example}\label{ex:Lag2}
Under the assumptions of Example \ref{x:lagemp}, if  $-\infty\leq\zeta_1<\zeta_2\leq\infty$ and $\zeta_2>-\gamma^{-1/2}$ for $\gamma>0$, then
\begin{equation}\label{eq:lc2ldp}
\lim_{n\to\infty}\frac1{\beta n^2}\log\mathbb P\left(\zeta_1\leq\sqrt{\frac pn}(\lambda_i-1)\leq\zeta_2, 1\leq i\leq n\right)
=-I_{V_\gamma^L}(\mu_{\zeta_1,\zeta_2}^L),
\end{equation}
where
\begin{align}\label{eq:lcon}
\mu_{\zeta_1,\zeta_2}^L
=\argmin_{\nu\in\mathcal M_1(D_\gamma\cap[\zeta_1,\zeta_2])}I_{V_\gamma^L}(\nu).
\end{align}
Its support is an interval $[s,t]\subset[\zeta_1,\zeta_2]$ given by
\begin{equation}\label{eq:lcsupp}
 [s,t]=
 \begin{cases}
 [\sqrt\gamma-2,\sqrt\gamma+2],
   &\zeta_1\leq\sqrt\gamma-2,\quad
    \zeta_2\geq\sqrt\gamma+2,\\
 [\zeta_1,b_\gamma(\zeta_1)],
   &\zeta_1>\sqrt\gamma-2,\quad
    \zeta_2\geq b_\gamma(\zeta_1),\\
 [a_\gamma(\zeta_2),\zeta_2],
   &\zeta_2<\sqrt\gamma+2,\quad
    \zeta_1\leq a_\gamma(\zeta_2),\\
 [\zeta_1,\zeta_2],&\text{otherwise}.
 \end{cases}
\end{equation}
The endpoint maps $a_\gamma$ and $b_\gamma$ are defined as follows. For $\ga=0$,
$$
 a_0(z)=\frac{z-2\sqrt{z^2+12}}3,
 \qquad
 b_0(z)=\frac{z+2\sqrt{z^2+12}}3.
$$
For $0<\gamma\leq1$, $b_\gamma(z)$ is the unique solution in $(z,\infty)$ to
$$
 1+\sqrt\gamma b_\gamma(z)
 +2(1-\gamma)\sqrt{\frac{1+\sqrt\gamma z}
 {1+\sqrt\gamma b_\gamma(z)}}
 =1+\sqrt\gamma z+2(1+\gamma).
$$
For $0<\gamma<1$, $a_\gamma(z)$ is the unique solution in $(-\gamma^{-1/2},z)$ to
$$
 1+\sqrt\gamma a_\gamma(z)
 +2(1-\gamma)\sqrt{\frac{1+\sqrt\gamma z}
 {1+\sqrt\gamma a_\gamma(z)}}
 =1+\sqrt\gamma z+2(1+\gamma).
$$
For $\gamma=1$, $a_1(z)=-1$.  The explicit trigonometric expressions are given in \eqref{eq:lrootl} and \eqref{eq:lrootr}.

Let $m=(s+t)/2$ and $h=(t-s)/2$.  
Then $\mu_{\zeta_1,\zeta_2}^L$ has density
\begin{equation}\label{eq:lcden}
 \rho_{\gamma;\zeta_1,\zeta_2}^L(x)=
 \begin{cases}
  
 \frac{2-x^2+mx+h^2/2}{2\pi\sqrt{(t-x)(x-s)}}
 \mathbf1_{[s,t]}(x),&\gamma=0,\\
  
 \frac{1+\gamma+\sqrt\gamma(m-x)
 -(1-\gamma)\frac{
 \sqrt{(1+\sqrt\gamma s)(1+\sqrt\gamma t)}}
 {1+\sqrt\gamma x}}
 {2\pi\gamma\sqrt{(t-x)(x-s)}}\mathbf1_{[s,t]}(x),
 &0<\gamma\leq1.
 \end{cases}
\end{equation}
Moreover,
\begin{equation}\label{eq:lcrate}
 I_{V_\gamma^L}(\mu_{\zeta_1,\zeta_2}^L)=
 \begin{cases}
 I_{V^H}(\mu_{\zeta_1/\sqrt2, \zeta_2/\sqrt2}^H),&\gamma=0,\\
 \Phi_\gamma^L(s,t),&0<\gamma\leq1.
 \end{cases}
\end{equation}
Here $I_{V^H}(\mu_{\zeta_1/\sqrt2, \zeta_2/\sqrt2}^H)$ is the Hermite rate function from Example \ref{ex:her2} and for $0<\gamma\leq1$, 
\begin{align}\label{eq:clphi}
 \Phi_\gamma^L(s,t)
 ={}&
 \frac{2+\sqrt\gamma(s+t)-\frac\gamma8(t-s)^2
 }{8\gamma^2}
 -\frac{1-\gamma}{4\gamma^2}
 \sqrt{(1+\sqrt\gamma s)(1+\sqrt\gamma t)}\notag\\
 &{}
 +\frac{\sqrt\gamma(s+t)-4}{8\gamma}
 -\frac12\log\frac{t-s}{2}-\frac{1-\gamma^2}{2\gamma^2}
 \log\bigl(\sqrt{1+\sqrt\gamma s}
            +\sqrt{1+\sqrt\gamma t}\bigr)\notag\\
 &+\frac{(1-\gamma)^2}{4\gamma^2}
 \log\frac{\sqrt{(1+\sqrt\gamma s)(1+\sqrt\gamma t)}}{1-\gamma}
 +\frac{\log2}{2\gamma^2},
 \end{align}
where the singular terms at $\gamma=1$ are defined by continuity. In particular,
\begin{equation}\label{eq:lcgam0}
 I_{V_0^L}(\mu_{\zeta_1,\zeta_2}^L)
 =\frac{m^2}{4}+\frac{h^2}{8}
   -\frac{m^2h^2}{16}-\frac{h^4}{128}
   -\frac12\log\frac h2-\frac38.
\end{equation}
\end{example}
\begin{remark}
Under the assumptions of Example \ref{ex:Lag2}, its
conclusion can be written in the original coordinates as
\begin{equation*}
 \lim_{n\to\infty}\frac1{\beta n^2}
 \log\PP\left(
  1+\sqrt{\frac np} \zeta_1\leq\lambda_i
  \leq1+\sqrt{\frac np} \zeta_2,
    1\leq i\leq n\right)
 =-I_{V_\gamma^L}(\mu_{\zeta_1,\zeta_2}^L).
\end{equation*}
For fixed $0\leq\xi_1<\xi_2$ and $0<\gamma\leq1$, continuity
of the constrained energy in the walls gives
\begin{equation}\label{eq:lclwal}
 \lim_{n\to\infty}\frac1{\beta n^2}
 \log\PP(\xi_1\leq\lambda_i\leq\xi_2,  1\leq i\leq n)
 =-\Phi_\gamma^L(s,t),
\end{equation}
where $\Phi_\gamma^L$ is given by \eqref{eq:clphi} and $[s,t]$
is determined by \eqref{eq:lcsupp} with
$\zeta_j=(\xi_j-1)/\sqrt\gamma$.
In particular, for $0<\xi<(1+\sqrt\gamma)^2$,
\eqref{eq:lclwal} and \eqref{eq:lrootl} give
\begin{align*}
 &\lim_{n\to\infty}\frac1{\beta n^2}
 \log\PP(\lambda_i\leq\xi,  1\leq i\leq n)\\
={}  & -\Phi_\gamma^L\Biggl(
 \frac{4(\xi+2+2\gamma)}{3\sqrt\gamma}
 \cos^2 \left[
 \frac13\arccos \left(
 -\frac{3\sqrt3(1-\gamma)\sqrt\xi}
 {(\xi+2+2\gamma)^{3/2}}\right)-\frac{2\pi}{3}\right]
 -\frac1{\sqrt\gamma},\frac{\xi-1}{\sqrt\gamma}\Biggr)\\
={}&-\frac12\Phi_-\left(
 \frac{(1+\sqrt\gamma)^2-\xi}{\gamma};
\gamma\right),
\end{align*}
where $\Phi_-$ is the rate function in \cite{VMB07}*{Eq.~(56)}. This agrees with their result for fixed $\beta$ after setting $N=n$, $M=p$, and $\tilde\lambda_i=p\lambda_i$, where $N,M$ are the dimension parameters and $\tilde\lambda_i$ the eigenvalues in \cite{VMB07}.
For $0<\gamma<1$, the equality of the two rate functions follows from rewriting
$$
 \operatorname{arcsinh}\sqrt v
 =\operatorname{arctanh}\sqrt{\frac{v}{1+v}}
 =\log\bigl(\sqrt v+\sqrt{1+v}\bigr),\qquad v>0
$$
in their expressions.
\end{remark}

\subsection{Lower tail probability of spectral radius of Coulomb gas}\label{se:srlow}
The empirical measure LDPs of real and complex Ginibre ensembles were proved in \cite{BAZ98,PH98}, and the spectral radius LDPs were established in \cite{XZ25} and later extended to elliptic ensembles in \cite{byu26}.
The lower tail large deviation probability of the spectral radius is closely related to the so called hole probability and has been studied extensively in both physics and mathematics; see, e.g., \cites{GHS88, Fo92, JLM93, APS09, CMV16, Ad18, Lac18,Lac19, Cha23,Cha24} and in particular the introduction of \cite{Cha24} for a nice review of this subject. 
We consider the probability that all eigenvalues lie in the annulus $\{z:\zeta_1\leq|z|\leq\zeta_2\}$, or equivalently that no eigenvalue lies in the hole $U=\{z:|z|<\zeta_1\}\cup\{z:|z|>\zeta_2\}$. Among the huge body of work on the hole probability problem, two papers are of particular interest to us. For fixed $\beta>0$, \cite{Ad18} studies hole probabilities for bounded holes within the equilibrium support for a class of radial potentials.
\cite{Cha23}*{Theorem~2.5(i)} extends \cite{Ad18}*{Theorem~1.2(B)} to nonradial potentials and possibly unbounded holes. These results are mainly concerned with general shapes of the holes and require stronger regularity on the external potentials, whereas our following result focuses on the annuli and needs weaker regularity conditions; see the remark afterwards for a discussion.

\begin{proposition}\label{thm:lowspe}
Under Assumption \ref{ass:ldp}, let $D\not\subset\mathbb R$ and $\boldsymbol\lambda=(\lambda_1,\dots,\lambda_n)\in D^n$ have density \eqref{pdf}.  Let $D_\zeta=\{z\in \C: \zeta_1\leq|z|\leq\zeta_2\}$. Assume further that $\mu_V=\argmin_{\mu \in \mathcal{M}_1(D)} I_V(\mu)$ and $\mu_{\zeta_1,\zeta_2}=\argmin_{\mu \in \mathcal{M}_1(D\cap D_\zeta)} I_V(\mu)$ satisfy Assumption \ref{ass:d} with $D$ and $D\cap D_\zeta$ respectively. For $0\le \zeta_1<\zeta_2\le\8$,  we have
\begin{align}\label{eq:cmz12}
\lim_{n\to\infty}
\frac{1}{\beta n^2}
\log \mathbb P(\zeta_1\le |\lambda_i|\le \zeta_2, i=1,\dots,n )
=
-I_V(\mu_{\zeta_1, \zeta_2}).
\end{align}
Assume further $D_\zeta \subset \supp(\mu_V)$. Then the minimizer can be written as  $\mu_{\zeta_1,\zeta_2} = \mu_V\big|_{D_\zeta}+\bar\mu_1+\bar\mu_2$. Here, $\bar\mu_2$ is supported on $\{|z|=\zeta_2\}$ and has density
$$
\sigma_{\zeta_2}(z)
=
\frac{1}{2\pi\zeta_2}
\int_{\{|w|>\zeta_2\}}
\frac{|w|^2-\zeta_2^2}{|w-z|^2} \mu_V(\dd w),
\qquad |z|=\zeta_2,
$$
with respect to arc length, while if $\zeta_1>0$, then $\bar\mu_1$ is supported on $\{|z|=\zeta_1\}$ and has density
$$
\sigma_{\zeta_1}(z)
=
\frac{1}{2\pi\zeta_1}
\int_{\{|w|<\zeta_1\}}
\frac{\zeta_1^2-|w|^2}{|z-w|^2}\mu_V(\dd w),
\qquad |z|=\zeta_1,
$$
with respect to arc length. When $\zeta_1=0$, we have $\bar\mu_1=0$. Moreover,
\begin{align*}
    I_V(\mu_{\zeta_1,\zeta_2})
    =
    \frac12\int V \dd\mu_{\zeta_1,\zeta_2} - \frac{1}{2}\int V  \dd \mu_V  +\frac{1}{2}\int_{\{|t|>\zeta_2\}}
\log\frac{|t|}{\zeta_2}\mu_V(\dd t).
\end{align*}
\end{proposition}
\begin{remark}
The constrained equilibrium measure and the rate function in Proposition \ref{thm:lowspe} already follow from \cite{Cha23}*{Proposition~3.3}. However, we allow the potential $V_n$ and $\beta=\beta_n$ to depend on $n$, with $\beta_n n\to\infty$. More importantly,
for annular confinement, our results require weaker conditions than those in \cites{Ad18,Cha23}.
First, the condition in \cite{Ad18}*{Theorem~1.2(A)} and Assumption \ref{ass:d} for the constrained minimizer both provide approximations for the probability lower bound. The former requires these approximations to be equilibrium measures on the complements of enlarged holes, whereas Proposition \ref{thm:lowspe} only needs arbitrary approximations. Second, \cite{Cha23}*{Theorem~2.5(i)} weakens the required potential regularity of \cite{Ad18}*{Theorem~1.2(B)} from Lipschitz continuity near the equilibrium support to H\"older continuity, while we need only continuity. 
To illustrate the weaker regularity allowed here, take $D=\mathbb C$ and let $V(0)=0$, $V(z)=1/\log(e/|z|)$ for $0<|z|\leq1$, and $V(z)=(|z|^2+1)/2$ for $|z|>1$. 
Its equilibrium measure is
$$
\mu_V(\dd z)
=\frac{\mathbf1_{\{0<|z|<1\}}}
{\pi |z|^2[\log(e/|z|)]^3} \dd \Re z \dd \Im z,
\qquad
\supp\mu_V=\{z:|z|\leq1\}.
$$
The potential $V$ is not H\"older continuous at $0\in\supp\mu_V$.
Nevertheless, its continuity and quadratic growth together with Remark \ref{r:a3unit} verify Assumption \ref{ass:ldp}.
For $0\leq\zeta_1<\zeta_2\leq1$, Lemma \ref{le:symreg} verifies Assumption \ref{ass:d} for $\mu_V$ on $\mathbb C$ and for $\mu_{\zeta_1,\zeta_2}$ on the corresponding disk or annulus.
Thus Proposition \ref{thm:lowspe} applies, whereas
\cite{Cha23}*{Theorem~2.5(i)} does not.
\end{remark}

\begin{proof}[Proof of Proposition \ref{thm:lowspe}] 
The measures $\mu_V$ and $\mu_{\zeta_1,\zeta_2}$ satisfy Assumption \ref{ass:d} with $D$ and $D\cap D_\zeta$, respectively. Proposition \ref{pr:logzn} therefore yields \eqref{eq:cmz12} by the same argument as in Proposition \ref{pr:lowgap}.

We next compute the value of the rate function $I_V(\mu_{\zeta_1,\zeta_2})$ under the condition that $D_\zeta\subset\operatorname{supp}(\mu_V).$
We first assume that $0<\zeta_1<\zeta_2<\infty$.
Decompose $\mu_V$ as
$
\mu_V
=
\mu_V^{in}
+
\mu_V^{mid}
+
\mu_V^{out},
$
where
$$
\mu_V^{in}
=
\mu_V\mathbf 1_{\{|w|<\zeta_1\}},
\qquad
\mu_V^{mid}
=
\mu_V\mathbf 1_{\{\zeta_1 \leq |w| \leq \zeta_2\}},
\qquad
\mu_V^{out}
=
\mu_V\mathbf 1_{\{|w| > \zeta_2\}}.
$$
Let $\bar\mu_1$ denote the balayage of $\mu_V^{in}$ onto the circle
$\{|z|=\zeta_1\}$, and let $\bar\mu_2$ denote the balayage of
$\mu_V^{out}$ onto the circle $\{|z|=\zeta_2\}$. Define
$
\mu_{\zeta_1,\zeta_2}
=
\bar\mu_1+\mu_V^{mid}+\bar\mu_2.
$
Consider the domains $G_1=\{|z|<\zeta_1\}$ and $G_2=\{|z|>\zeta_2\}.$
For $i=1,2$, let $g_{G_i}$ denote the Green function of $G_i$. Then
$$
g_{G_i}(z,w)
=
\log\left|
\frac{\zeta_i^2-z\overline w}
{\zeta_i(z-w)}
\right|,
\qquad
z,w\in G_i,\quad z\neq w.
$$
In particular,
$g_{G_2}(z,\infty)=\log (|z|/\zeta_2).$

By the balayage formula \cite{ST97}*{Theorem II.4.11}, the measures $\bar\mu_1$ and $\bar\mu_2$ are absolutely continuous with respect to the arc-length measures on $\{|z|=\zeta_1\}$ and $\{|z|=\zeta_2\}$, respectively. Their densities are given by
$$
\sigma_{\zeta_1}(z)
=
\frac{1}{2\pi}
\int_{\{|w|<\zeta_1\}}
\frac{\partial g_{G_1}(z,w)}
{\partial\mathbf n_1^+} 
\mu_V(\dd w) =
\frac{1}{2\pi\zeta_1}
\int_{\{|w|<\zeta_1\}}
\frac{\zeta_1^2-|w|^2}
{|z-w|^2} 
\mu_V(\dd w),
\qquad |z|=\zeta_1,
$$
\begin{align}\label{densig2}
\sigma_{\zeta_2}(z)
=
\frac{1}{2\pi}
\int_{\{|w|>\zeta_2\}}
\frac{\partial g_{G_2}(z,w)}
{\partial\mathbf n_2^+} 
\mu_V(\dd w) 
=
\frac{1}{2\pi}
\int_{\{|w|>\zeta_2\}}
\frac{|w|^2-\zeta_2^2}
{\zeta_2|w-z|^2} 
\mu_V(\dd w),
\qquad |z|=\zeta_2.
\end{align}
Here, $\mathbf n_1^+$ denotes the inward unit normal to $G_1$, while
$\mathbf n_2^+$ denotes the outward unit normal to the circle
$\{|z|=\zeta_2\}$.

Since every boundary point of a disk is regular, the balayage property \cite{ST97}*{Theorem II.4.1} yields\footnote{Literally this theorem requires $\supp \mu_V^{in}\subset G_1$; however, by an approximation argument as for \cite{ST97}*{Theorem II.4.7} it readily extends to our current setting.}
$$
\int \log|z-w|\bar\mu_1(\dd w)
=
\int \log|z-w|\mu_V^{in}(\dd w),
\qquad   |z|\geq \zeta_1.
$$
Similarly, by \cite{ST97}*{Theorem II.4.7},
$$
\int \log|z-w|\mu_V^{\mathrm{out}}(\dd w)
=
\int \log|z-w|\bar\mu_2(\dd w)+c_2(\zeta_2), 
\qquad 
|z|\le \zeta_2,
$$
where
$$
c_2(\zeta_2)
=
\int_{\{|t|>\zeta_2\}} g_{G_2}(t,\infty) \mu_V(\dd t)
=
\int_{\{|t|>\zeta_2\}}
\log\frac{|t|}{\zeta_2}\mu_V(\dd t).
$$
Consequently, for $\zeta_1\leq|z|\leq\zeta_2$,
\begin{equation}\label{eq:lps}
\int \log|z-w|\mu_{\zeta_1,\zeta_2}(\dd w)
=
\int \log|z-w|\mu_V(\dd w)-c_2(\zeta_2).
\end{equation}

Since $\mu_V$ is the equilibrium measure associated with $V$, there exists $F_V\in\mathbb R$ such that
\begin{align}\label{eq:equimuv}
V(z)-\int_{\mathbb C}\log|z-w| \mu_V(\dd w)
=
F_V,
\qquad \mbox{q.e. }  z\in \supp \mu_V.
\end{align}
By the definition of $C_V$ and the integrated equilibrium relation \eqref{eq:calsig},
\begin{align}\label{eq:cvfv}
C_V
=
\int V \dd\mu_V-\frac12\Sigma(\mu_V)
=
\frac12\int V \dd\mu_V+\frac12F_V.
\end{align}
Combining \eqref{eq:lps} and the equilibrium condition \eqref{eq:equimuv} with the assumption $D_\zeta\subset \supp \mu_V$, we obtain
\begin{align}\label{eq:vtzeta}
V(z)-\int_{\mathbb C}\log|z-w| 
\mu_{\zeta_1,\zeta_2}(\dd w)
=
F_V+c_2(\zeta_2) \qquad \mbox{q.e. }  z\in \supp \mu_{\zeta_1,\zeta_2}.
\end{align}
Lemma \ref{lemr} therefore implies that $\mu_{\zeta_1,\zeta_2}$ is the unique equilibrium measure. 

Finally, by \eqref{eq:cvfv} and \eqref{eq:vtzeta},
\begin{align*}
I_V(\mu_{\zeta_1,\zeta_2})
=
\frac12\int V \dd\mu_{\zeta_1,\zeta_2}
+\frac12F_V
+\frac12c_2(\zeta_2)
-C_V =
\frac12\int V \dd\mu_{\zeta_1,\zeta_2}
-\frac12\int V \dd\mu_V
+\frac12c_2(\zeta_2).
\end{align*}
If $\zeta_1=0$, the finiteness of the logarithmic energy implies $\mu_V(\{0\})=0$, so $\mu_V^{in}=\bar\mu_1=0$ and all terms involving the inner balayage are omitted.  
This is the desired formula, and the proof is complete.
\end{proof}

\begin{remark}
The condition $D_\zeta \subset \supp(\mu_V)$ is crucial in the above proof. Without this condition, we could still compute the balayage, but we cannot verify the equilibrium condition without knowing the support of the minimizer. 
\end{remark}

\section{Large deviations of classical \texorpdfstring{$\beta$}{beta}-ensembles}\label{sec:gen}

In this section, we prove the large deviation results for classical $\bt$-ensembles as presented in the preceding section, including LDPs of the empirical measure and the largest eigenvalue, and lower tail probability of the latter. The proofs are direct applications of Theorem \ref{ldpln}, Propositions \ref{pr:maxldp} and \ref{pr:lowgap}. Basically, we need to verify Assumptions \ref{ass:ldp}, \ref{ass:d} and \ref{ass:maximum}, and compute the rate functions for the three types of classical ensembles.

\subsection{The Hermite and Laguerre ensembles}
We start with the lower tail of the largest eigenvalue of Hermite ensemble.

\begin{proof}[Proof of Example \ref{ex:her2}]
For fixed $\beta>0$, \cite{BAG97} shows that the empirical measure of the $\beta$-Hermite ensemble satisfies an LDP with
speed $\beta n^2$ and good rate function 
\begin{align}\label{eq:herfun}
    I_{V_H}(\mu)= \int V_H (x) \mu(\dd x) -\frac{1}{2} \Sigma(\mu)-C_{V_H},
    \qquad C_{V_H}
=
\frac38+\frac14\log2.
\end{align}
The unique minimizer of $I_{V_H}$ over
$\mathcal M_1(\mathbb R)$ is the semicircle law $\mu_{sc}$ with density given in \eqref{eq:semic}.

We next verify that Theorem \ref{ldpln} applies whenever
$\beta n\to\infty$.  In the weighted representation
\eqref{eq:wpdf}, $W_{n, H}=V^H=x^2/2$ and
$\omega_{n, H}=1$.  Thus \ref{A1} and \ref{A2} hold immediately, while Remark
\ref{r:a3unit} gives \ref{A3}.
Finally, since $V^H$ is finite and continuous on $\mathbb R$, Lemma \ref{le:symreg} verifies Assumption \ref{ass:d} on $\mathbb R$ for every $\mu\in\mathcal M_1(\mathbb R)$ with $I_{V^H}(\mu)<\infty$. Theorem \ref{ldpln} therefore gives the LDP with speed $\beta n^2$ and rate function \eqref{eq:herfun} under the sole condition $\beta n\to\infty$. Hence Proposition \ref{pr:lowgap} yields, for every $\zeta_1 <\zeta_2$,
\begin{align*}
\lim_{n\to\infty} \frac{1}{\beta n^2} \log \mathbb P\bigl(\zeta_1\le \lambda_i \le \zeta_2,  1\leq i \leq n\bigr)
=- I_{V_H}( \mu_{\zeta_1,\zeta_2}^H), 
\qquad
 \mu_{\zeta_1,\zeta_2}^H =\argmin_{\nu\in\mathcal M_1([\zeta_1,\zeta_2])} I_{V_H}(\nu).
\end{align*}
The constrained energy identity in Proposition \ref{pr:lowgap} takes the form
\begin{align}\label{eq:ivh1}
I_{V_H}(\mu_{\zeta_1,\zeta_2}^H)
=
\frac12\int V_H(x) \mu_{\zeta_1,\zeta_2}^H(dx)
+
\frac12F_{\zeta_1,\zeta_2}
-
C_{V_H},	
\end{align}
where $F_{\zeta_1,\zeta_2}$ is determined by \eqref{eq:fzeta}. 

We first identify the support of the $\mu_{\zeta_1,\zeta_2}^H$, defined in \eqref{eq:hcon}. If $\zeta_1<-\sqrt{2}$ and $\zeta_2\ge\sqrt2$, then $\mu_{sc} \in\mathcal M_1([\zeta_1,\zeta_2])$. 
Therefore $\mu_{\zeta_1,\zeta_2}^H=\mu_{sc}$ and  $[s, t]= [-\sqrt{2}, \sqrt{2}].$
It remains to consider the constrained cases. 
Applying Lemma
\ref{remr} with $D=[\zeta_1,\zeta_2]$ and $V_H'(x)=x$, the support
of $\mu_{\zeta_1,\zeta_2}^H$ is a single interval
$[s,t]\subset[\zeta_1,\zeta_2]$. If only the left endpoint is constrained, then $s=\zeta_1$, while the
right endpoint $t$ is determined by
$$
    \frac1\pi\int_{\zeta_1}^t
    x\sqrt{\frac{x-\zeta_1}{t-x}}\dd x=1.
$$
Under the change of variables $x=\zeta_1+(t-\zeta_1)u$, the integral becomes
\begin{align*}
\frac1\pi\int_{\zeta_1}^t
x\sqrt{\frac{x-\zeta_1}{t-x}}\dd x   
&=
\frac{t-\zeta_1}{\pi}
\left[
\zeta_1\int_0^1 u^{1/2}(1-u)^{-1/2}\dd u
+
(t-\zeta_1)
\int_0^1 u^{3/2}(1-u)^{-1/2}\dd u
\right]\\
&=
\frac{t-\zeta_1}{\pi}
\left[
\zeta_1\cdot\frac{\pi}{2}
+
(t-\zeta_1)\cdot\frac{3\pi}{8}
\right].
\end{align*}
Therefore $t$ satisfies
$\frac{t-\zeta_1}{2}\zeta_1+\frac{3}{8}(t-\zeta_1)^2=1.$
Solving this quadratic equation and choosing the root $t>\zeta_1$, we
obtain
$$
    t
    =
    \frac{\zeta_1+2\sqrt{\zeta_1^2+6}}{3}
    =
    b(\zeta_1).
$$
Thus, if $\zeta_1>-\sqrt{2}$ and $\zeta_2 >b(\zeta_1),$ then $[s, t]=[\zeta_1, b(\zeta_1)].$ If instead $\zeta_1>-\sqrt{2}$ and $\zeta_2 <b(\zeta_1),$ then $[s, t]=[\zeta_1, \zeta_2].$
Similarly, if only the right endpoint is constrained, then $t=\zeta_2$,
while $s$ is determined by
$$
    \frac1\pi\int_s^{\zeta_2}
    x\sqrt{\frac{\zeta_2-x}{x-s}}\dd x=-1.
$$
Equivalently, $    3s^2-2\zeta_2s-\zeta_2^2-8=0.$
Choosing the root $s<\zeta_2$, we get
$$
    s
    =
    \frac{\zeta_2-2\sqrt{\zeta_2^2+6}}{3}
    =
    a(\zeta_2).
$$
Therefore, if $\zeta_1<a(\zeta_2)$, then $ [s,t]=[a(\zeta_2),\zeta_2].$ If $\zeta_1>a(\zeta_2)$, then $ [s,t]=[\zeta_1,\zeta_2].$ Combining the above cases gives the stated formula for $[s,t]$.

We next compute the density. For fixed support $[s,t]$, Lemma \ref{remr} gives
$$
\rho_{\zeta_1,\zeta_2}^H(x)
=
\frac{\pi-\ell(x)}{\pi^2\sqrt{(t-x)(x-s)}},
\qquad x\in(s,t),
$$
where
\begin{align*}
	\ell(x)
&=
\pv\int_s^t
\sqrt{(t-y)(y-s)} \frac{V_H'(y)}{x-y}\dd y
=
\pv\int_s^t
\sqrt{(t-y)(y-s)} \frac{y}{x-y}\dd y\\
& =
x \pv\int_s^t
\frac{\sqrt{(t-y)(y-s)}}{x-y}\dd y
-
\int_s^t
\sqrt{(t-y)(y-s)}\dd y
 = \pi x(x-m)
-
\frac{\pi h^2}{2},
\end{align*}
where in the last equality, we used Lemma \ref{le:Ia} with $\alpha=s$ and $\beta=t$.
It follows that
$$
\rho_{\zeta_1,\zeta_2}^H(x)=
\frac{
1-x^2+mx+\frac{h^2}{2}
}{
\pi\sqrt{(t-x)(x-s)}
},
\qquad x\in(s,t).
$$

We now compute $F_{\zeta_1,\zeta_2}$. Taking $x=t$ in
\eqref{eq:fzeta} and substituting $y=m+h\cos\theta$, we obtain
\begin{align*}
F_{\zeta_1,\zeta_2}
&=\frac{t^2}{2}
-\int_s^t\log|t-y|\rho_{\zeta_1,\zeta_2}^H(y)\dd y\\
&=\frac{t^2}{2}-\log h
-\frac1\pi\int_0^\pi\log(1-\cos\theta)
\left(1+\frac{h^2}{2}-h^2\cos^2\theta-mh\cos\theta\right)\dd\theta\\
&=\frac{t^2}{2}-\log\frac h2-\frac{h^2}{4}-mh
=\frac{m^2}{2}+\frac{h^2}{4}-\log\frac h2,
\end{align*}
where we used integration by parts and $t=m+h$ in the last step.
Since $V_H(x)=x^2/2$, a similar calculation yields
\begin{align*}
\int x^2 \mu_{\zeta_1,\zeta_2}^H(\dd x)
&=
\frac1\pi
\int_0^\pi
\bigl(m^2+2mh\cos\theta+h^2\cos^2\theta\bigr)
\left(
1+\frac{h^2}{2}
-h^2\cos^2\theta
-mh\cos\theta
\right)
\dd\theta  \\
&=
m^2+\frac{h^2}{2}
-
m^2h^2
-
\frac{h^4}{8}.
\end{align*}
Substituting the above identities into \eqref{eq:ivh1}, we obtain
\begin{align*}
I_{V_H}(\mu_{\zeta_1,\zeta_2}^H)
&=
\frac12\cdot
\frac12
\left(
m^2+\frac{h^2}{2}-m^2h^2-\frac{h^4}{8}
\right)
+
\frac12
\left(
\frac{m^2}{2}
+
\frac{h^2}{4}
-
\log\frac{h}{2}
\right)
-
\frac38-\frac14\log2\\
&=
\frac{m^2}{2}
+
\frac{h^2}{4}
-
\frac{m^2h^2}{4}
-
\frac{h^4}{32}
-
\frac12\log\frac{h}{2}
-
\frac38
-
\frac14\log2.
\end{align*}

Finally, for $\zeta<\sqrt2$ the one sided support is $[a(\zeta),\zeta]$, whereas for $0<\zeta<\sqrt2$ the symmetric support is $[-\zeta,\zeta]$.  Substituting these endpoints into the last formula and simplifying gives the two special rate functions in the statement.
\end{proof}

Next, we consider the Laguerre ensemble. The following argument in particular resolves the aforementioned issue for $\ga=1$ in empirical measure LDP.

\begin{proof}[Proof of Example \ref{x:lagemp}]
To prove (i), it suffices to verify the assumptions of Theorem \ref{ldpln}. The constant $C_{V_\gamma^L}$ is given in \cite{Ma23L}. 
Take $D=\mathbb R$ and let $V_\gamma^L=+\infty$ outside $D_\gamma$. For $y>-\sqrt{p/n}$, define
\begin{align}\label{e:Lagwgt}
 W_{n, L}(y)
 &=\frac12\sqrt{\frac pn} y
 -\frac{p-n+1}{2n}\log\left(1+\sqrt{\frac np} y\right),
 \qquad
 \omega_{n, L}(y)=\left(1+\sqrt{\frac np} y\right)^{-1}.
\end{align}
Let $W_{n,L}=+\infty$ and $\omega_{n,L}=1$ for $y\leq-\sqrt{p/n}$.
Let $(\lambda_1,\ldots,\lambda_n)$ be $\beta$-Laguerre ensemble with potential \eqref{def:vnl}. Thus $(\sqrt{p/n}(\lambda_1-1),\ldots,\sqrt{p/n}(\lambda_n-1))$
has the weighted density \eqref{eq:wpdf} with
$W_n=W_{n,L}$ and $\omega_n=\omega_{n,L}$.

Assumption \ref{A0} is immediate.
For \ref{A1}, lower semicontinuity follows from the definitions. \ref{A1}(ii) is immediate for $\gamma>0$ and follows from a second order Taylor expansion when $\gamma=0$.
Here the set $\mathcal O$ in Assumption \ref{ass:ldp} is the interior of $D_\gamma$.
Thus $D\setminus\mathcal O=(-\infty,-\gamma^{-1/2}]$ when $\gamma>0$ and is empty when $\gamma=0$.
It remains to verify \ref{A1}(i) for $\gamma>0$:
\begin{equation}\label{e:Laglim}
\liminf_{n\to\infty}W_{n,L}(y_n)\geq V_\gamma^L(y)
\qquad\text{whenever }y_n\to y\in D\setminus\mathcal O.
\end{equation}
Since $W_{n,L}(y)=+\infty$ for $y\leq-\sqrt{p/n}$, it suffices to consider $y_n>-\sqrt{p/n}$ with $y_n\to-\gamma^{-1/2}$.
For $0<\gamma<1$, $\liminf_{n\to\infty}W_{n,L}(y_n)=+\infty =V_\gamma^L(-\gamma^{-1/2})$.
For $\gamma=1$, the logarithmic term $-\frac{p-n+1}{2n}\log(1+\sqrt{n/p} y_n)$ is nonnegative for all sufficiently large $n$. Hence
$$
\liminf_{n\to\infty}W_{n,L}(y_n)
\geq\lim_{n\to\infty}\frac12\sqrt{\frac pn} y_n
=-\frac12=V_1^L(-1).
$$
Thus \eqref{e:Laglim} holds in both cases, proving \ref{A1}(i).

Assumption \ref{A2} is immediate for $V_\gamma^L$.
For $y\geq0$, \eqref{e:logpos} gives 
\begin{align}
W_{n,L}(y)
&=  \frac{p}{2n}\left(\sqrt{\frac{n}{p}} y-\log\Big(1+\sqrt{\frac{n}{p}} y\Big)\right)
+\frac{n-1}{2n}\log\Big(1+\sqrt{\frac{n}{p}} y\Big) \notag\\
&\geq \frac{y^2}{4(1+\sqrt{n/p} y)}
\geq\frac{y^2}{4(1+y)}, \label{e:Lagpos}
\end{align}
where in the last inequality we used $p\geq n$.
For $-\sqrt{p/n}<y<0$, the lower bound in \eqref{e:logneg} yields
\begin{align*}
W_{n,L}(y)
\geq\frac{n-1}{2n}\sqrt{\frac np} y+\frac{p-n+1}{4p}y^2
\geq-\frac12+\frac{1-n/p}{4}y^2
\geq\frac{y^2}{4}-\frac34,
\end{align*}
where in the last two inequalities we used $-1<\sqrt{n/p}\,y<0$ and $(n/p)y^2<1$.
Since $W_{n,L}=+\infty$ on $(-\infty,-\sqrt{p/n}]$, the two bounds above imply \ref{A2}.

By \eqref{e:Lagwgt}, $\omega_{n,L}$ is positive, finite, and measurable, with $\log\omega_{n,L}=O(1)=o(\beta n)$ locally uniformly in the interior of $D_\gamma$. Thus \ref{A3}(i) holds.
For \ref{A3}(ii), take $\theta_n=(\beta n)^{-1/2}$. By \eqref{eq:def-frac} and \eqref{e:Lagwgt},
\begin{align*}
 & \mathcal R_{n, L}(\theta_n)
 ={}\int_{(-\sqrt{p/n},\infty)^n}
 \prod_{i<j}|y_i-y_j|^{\sqrt{\beta n}/n}
 \prod_{i=1}^n e^{-\sqrt{\beta p} y_i/2}
 \left(1+\sqrt{\frac np}y_i\right)^{
 \frac{\sqrt{\beta n}(p-n+1)}{2n}-1}\dd\boldsymbol y\\
 ={}&e^{p\sqrt{\beta n}/2}
 \left(\frac pn\right)^{n/2+\sqrt{\beta n}(n-1)/4}
 \int_{(0,\infty)^n}\prod_{i<j}|x_i-x_j|^{\sqrt{\beta n}/n} \prod_{i=1}^n
 x_i^{\frac{\sqrt{\beta n}(p-n+1)}{2n}-1}
 e^{-p\sqrt{\beta n} x_i/(2n)}\dd\boldsymbol x\\
 ={}&e^{p\sqrt{\beta n}/2}
 \left(\frac pn\right)^{n/2+\sqrt{\beta n}(n-1)/4}
 \left(\frac{p\sqrt{\beta n}}{2n}\right)^{-p\sqrt{\beta n}/2}\prod_{j=1}^n
 \frac{\Gamma\bigl(1+j\sqrt{\beta n}/(2n)\bigr)
       \Gamma\bigl(\frac{\sqrt{\beta n}(p-n+j)}{2n}\bigr)}
      {\Gamma\bigl(1+\sqrt{\beta n}/(2n)\bigr)},
\end{align*}
where the second equality uses $y_i=\sqrt{p/n}(x_i-1)$ and the last equality follows from Selberg integral; see, e.g., \cite{AGZ10}*{Eq.~(2.5.10)}. 
By \eqref{bd:stir2}, we get
\begin{align*}
 &\log\left\{
 e^{p\sqrt{\beta n}/2}
 \left(\frac pn\right)^{n/2+\sqrt{\beta n}(n-1)/4}
 \left(\frac{p\sqrt{\beta n}}{2n}\right)^{-p\sqrt{\beta n}/2}
 \prod_{j=1}^n
 \Gamma\left(\frac{\sqrt{\beta n}(p-n+j)}{2n}\right)
 \right\}\\
\leq{} & \frac{p\sqrt{\beta n}}2
 +\frac{2n+\sqrt{\beta n}(n-1)}4\log\frac pn
 +Cn+C\sum_{j=1}^n
 \left(-\log\frac{\sqrt{\beta n}(p-n+j)}{2n}\right)_+
 \\
 & -\frac{p\sqrt{\beta n}}2\log\frac{p\sqrt{\beta n}}{2n} +\sum_{j=1}^n\left[
 \frac{\sqrt{\beta n}(p-n+j)-n}{2n}
 \log\frac{\sqrt{\beta n}(p-n+j)}{2n}
 -\frac{\sqrt{\beta n}(p-n+j)}{2n}\right]\\
 ={}&
 \frac{\sqrt{\beta n}(n-1)}4
 -\frac{2n+\sqrt{\beta n}(n-1)}4
 \log\frac{\sqrt{\beta n}}2
 +\frac{\sqrt{\beta n}}2\sum_{j=1}^n
 \frac{p-n+j}{n}\log\frac{p-n+j}{p}\\
 & -\frac12\sum_{j=1}^n\log\frac{p-n+j}{p}
 +Cn+C\sum_{j=1}^n
 \left(-\log\frac{\sqrt{\beta n}(p-n+j)}{2n}\right)_+.
\end{align*}
Since $p\ge n$, the weighted logarithmic sum is nonpositive, and \eqref{e:factbd} gives 
\begin{align*}
 -\sum_{j=1}^n\log\frac{p-n+j}{p}
 &\leq\sum_{j=1}^n\log\frac nj
 =n\log n-\log(n!)=O(n).
\end{align*}
As $\beta n\to\infty$, 
for all sufficiently large $n$,
$$
 \sum_{j=1}^n
 \left(-\log\frac{\sqrt{\beta n}(p-n+j)}{2n}\right)_+
 \leq\sum_{j=1}^n\log\frac{2n}{j}=O(n).
$$
By \eqref{e:gambd}, we obtain
\begin{align}
 &\log\prod_{j=1}^n
 \frac{\Gamma\bigl(1+j\sqrt{\beta n}/(2n)\bigr)}
      {\Gamma\bigl(1+\sqrt{\beta n}/(2n)\bigr)}
 =\sum_{j=1}^n
 \log\Gamma\left(1+\frac{j\sqrt{\beta n}}{2n}\right)
 -n\log\Gamma\left(1+\frac{\sqrt{\beta n}}{2n}\right) \notag \\
 \leq {}& C\sum_{j=1}^n
 \left(1+\frac{j\sqrt{\beta n}}{2n}\right)
 \log\left(2+\frac{j\sqrt{\beta n}}{2n}\right)+Cn\leq C\log(2+\sqrt{\beta n})
 \sum_{j=1}^n
 \left(1+\frac{j\sqrt{\beta n}}{2n}\right)+Cn \notag \\
 ={}& C\log(2+\sqrt{\beta n})
 \left(n+\frac{n+1}{4}\sqrt{\beta n}\right)+Cn 
 \leq Cn(1+\sqrt{\beta n})\log(2+\sqrt{\beta n}). \label{bd:gamj}
\end{align}
Combining these estimates gives, uniformly in $p\geq n$,
$$
 \log\mathcal R_{n, L}(\theta_n)
 \leq Cn(1+\sqrt{\beta n})\log(2+\sqrt{\beta n})
 =o(\beta n^2),
$$
which proves \ref{A3}(ii). 

Note that $\mu(D_\gamma)=1$ for $\mu\in\mathcal M_1(\mathbb R)$ with $I_{V_\gamma^L}(\mu)<\infty$. Using Lemma \ref{le:symreg} with $D=\mathbb R$ and $E=D_\gamma$, Assumption \ref{ass:d} holds.
Thus, part (i) follows from Theorem \ref{ldpln}.

To prove (ii), we consider $D=\mathbb R_+$, 
$$
 W_n(x)=\frac{p}{2n}x-\frac{p-n+1}{2n}\log x,
 \qquad \omega_n(x)=x^{-1},\qquad x>0,
$$
with $W_n(0)=+\infty$ and $\omega_n(0)=1$.
The relation $y=\sqrt{p/n}(x-1)$ yields $W_n(x)=W_{n,L}(y)+p/(2n)$ and
$\omega_n(x)=\omega_{n,L}(y)$. Since $n/p\to\gamma>0$, the preceding convergence and growth estimates verify \ref{A1} and \ref{A2}, and \ref{A3}(i) holds trivially.
Observing that the logarithm of the integral in \eqref{eq:def-frac} is affected up to $O(n+\theta_n\beta n^2)=o(\beta n^2)$ by the change of variables, \ref{A3}(ii) follows.
Lemma \ref{le:symreg} verifies Assumption \ref{ass:d} for every $\mu\in\mathcal M_1(\mathbb R_+)$ with $I_{\widehat V_\gamma^L}(\mu)<\infty$.
Theorem \ref{ldpln} therefore proves (ii).
By \eqref{eq:loremp},
$$
V_\gamma^L(\varphi(x))=\widehat V_\gamma^L(x)-\frac1{2\gamma},
\qquad
\Sigma(\varphi_\#\mu)=\Sigma(\mu)-\frac12\log\gamma.
$$
The resulting additive constant in the energy cancels after subtracting its minimum, which proves \eqref{eq:loremp}.
\end{proof}

The following recovers the largest eigenvalue LDP for Laguerre ensemble, where the cases $\ga=1$ and $\bt\neq 2$ seem to be overlooked before. By \cite{Ma23L}, the unique unconstrained minimizer $\widetilde\nu_\gamma^L$ of $I_{V_\gamma^L}$ over $\mathcal M_1(D_\gamma)$ has density
\begin{equation}\label{eq:lmin}
\widetilde\rho_\gamma^L(x)=
\begin{cases}
\sqrt\gamma \rho_\gamma^L(1+\sqrt\gamma x),
&0<\gamma\leq1,\\
\frac1{2\pi}\sqrt{4-x^2} \mathbf1_{[-2,2]}(x),
&\gamma=0,
\end{cases}
\end{equation}
where $\rho_\gamma^L$ is given in \eqref{eq:hgam}.

\begin{proof}[Proof of Example \ref{ex:lag-max}]
 To prove (i), use \eqref{e:Lagwgt} on $D=\mathbb R$,  with $\mathcal O$ the interior of $D_\gamma$, where $D_\gamma=[-\gamma^{-1/2},\infty)$ for $\gamma>0$ and $D_0=\mathbb R$.
The proof of Example \ref{x:lagemp} verifies Assumption \ref{ass:ldp} and Assumption \ref{ass:d} for every $\mu\in\mathcal M_1(\mathbb R)$ with $I_{V_\gamma^L}(\mu)<\infty$.
By Proposition \ref{pr:maxldp}, it remains to verify Assumption \ref{ass:maximum}. Condition \ref{C0} holds by assumption.

Suppose first that $0<\gamma\leq1$.
The density in \eqref{eq:lmin} gives $b_{V_\gamma^L}=\max\supp(\widetilde\nu_\gamma^L)=2+\sqrt\gamma$.
By \eqref{eq:mxrate},
\begin{equation}\label{eq:lomax}
J_{V_\gamma^L}(y)
=V_\gamma^L(y)
-\int_{\sqrt\gamma-2}^{\sqrt\gamma+2}
\log|y-t| \widetilde\rho_\gamma^L(t)\dd t-F_{V_\gamma^L},
\qquad y\geq2+\sqrt\gamma,
\end{equation}
For $0<\gamma\leq1$, \eqref{eq:arcavg} and
Lemma \ref{le:intlog} give
\begin{align*}
F_{V_\gamma^L}
&=\frac1\pi\int_{\sqrt\gamma-2}^{\sqrt\gamma+2}
\frac{V_\gamma^L(x)}{\sqrt{4-(x-\sqrt\gamma)^2}}\dd x\\
&=\frac1{2\pi\sqrt\gamma}\int_0^\pi
(\sqrt\gamma+2\cos\theta)\dd\theta
-\frac{1-\gamma}{2\pi\gamma}\int_0^\pi
\log(1+\gamma+2\sqrt\gamma\cos\theta)\dd\theta\\
&=\frac12-\frac{1-\gamma}{2\gamma}
\log\frac{1+\gamma+\sqrt{(1+\gamma)^2-4\gamma}}{2}=\frac12.
\end{align*}
Here we used $x=\sqrt\gamma+2\cos\theta$.
Differentiating \eqref{eq:lomax} and using the Stieltjes transform in \cite{BS10}*{Lemma~3.11}, we obtain
\begin{align*}
J'_{V_\gamma^L}(y)
=\frac{y+\sqrt\gamma}{2(1+\sqrt\gamma y)}
-\int_{\gamma_-}^{\gamma_+}
\frac{\sqrt\gamma \rho_\gamma^L(t)}{1+\sqrt\gamma y-t}\dd t
=\frac{\sqrt{(y-\sqrt\gamma)^2-4}}
{2(1+\sqrt\gamma y)}>0,
\qquad y>2+\sqrt\gamma.
\end{align*}
Here $\gamma_\pm=(1\pm\sqrt\gamma)^2$.
Since $J_{V_\gamma^L}(2+\sqrt\gamma)=0$, integration gives
\begin{equation}\label{e:Lagrate}
J_{V_\gamma^L}(y)
=\frac12\int_{2+\sqrt\gamma}^y
\frac{\sqrt{(t-\sqrt\gamma)^2-4}}{1+\sqrt\gamma t}\dd t,
\qquad y\geq2+\sqrt\gamma.
\end{equation}
 Together with $b_{V_\gamma^L}\in\mathcal O$, this verifies \ref{C1}.

For \ref{C2}, we verify \eqref{eq:mxend} at $-\gamma^{-1/2}$.
Since $-\sqrt{p/n}\to-\gamma^{-1/2}$, the integrand vanishes on every compact subset of $\mathbb R\setminus D_\gamma$ for all sufficiently large $n$. For fixed $0<\delta<\gamma^{-1/2}$, the substitution $x=1+\sqrt{n/p} y$ yields
\begin{align}\label{e:Lagsbd}
& \limsup_{n\to\infty}\frac1{\beta n}
\log\int_{-1/\sqrt{\gamma} -\delta}^{-1/\sqrt{\gamma}+\delta}
e^{-\beta nW_{n,L}(y)}\omega_{n,L}(y)\dd y
 \leq\limsup_{n\to\infty}\frac1{\beta n}
\log\int_{-\infty}^{-1/\sqrt{\gamma}+\delta}
e^{-\beta nW_{n,L}(y)}\omega_{n,L}(y)\dd y\notag\\
\leq{} & \limsup_{n\to\infty}\frac1{\beta n}
\log
\frac{2\sqrt{p/n} e^{\beta p/2}
\left[1-\sqrt{n/(\gamma p)}+\sqrt{n/p}\,\delta\right]^{\beta(p-n+1)/2}}
{\beta(p-n+1)}
\leq\frac{1-\gamma}{2\gamma}\log(\sqrt\gamma\,\delta)
+\frac1{2\gamma},
\end{align}
where we used that the integrand vanishes for $y\leq-\sqrt{p/n}$, that $p-n+1\geq1$, and that $\log^+(1/\beta)=o(\beta n)$ by \ref{C0}.
Since
$$
\lim_{\delta\downarrow0}\left[
\frac{1-\gamma}{2\gamma}\log(\sqrt\gamma\,\delta)
+\frac1{2\gamma}\right]
=-V_\gamma^L(-\gamma^{-1/2}),
$$
letting $\delta\downarrow0$ in \eqref{e:Lagsbd} proves \eqref{eq:mxend} and hence Remark \ref{re:mxendpt} gives \ref{C2}.

For \ref{C3}, 
since $\omega_{n,L}(y)\leq1$ for $y\geq1$, by \eqref{e:Lagpos} we have 
$$
W_{n,L}(y)-\frac1{\beta n}\log\omega_{n,L}(y)
\geq W_{n,L}(y)
\geq\frac{y^2}{4(1+y)}
\geq\frac y8.
$$
Then Remark \ref{rem:mxpow} yields \ref{C3}.

If $\gamma=0$, then $\mathcal O=D=\mathbb R$ and $V_0^L(x)=V^H(x/\sqrt2)$. Using the densities in \eqref{eq:lmin} and \eqref{eq:semic}, we find $b_{V_0^L}=2$ and
$$
\int\log|x-y| \widetilde\nu_0^L(\dd y)
=\frac12\log2+
\int\log|x/\sqrt2-t|\rho^H(t)\dd t.
$$
Evaluating \eqref{eq:vlogmu} at $x=2$ shows that $F_{V_0^L}=F_{V^H}-\frac12\log2$. 
Hence \eqref{eq:mxrate} yields, for $x\geq2$,
\begin{align*}
J_{V_0^L}(x)
&=V_0^L(x)-\int\log|x-y| \widetilde\nu_0^L(\dd y)-F_{V_0^L}\\
&=V^H(x/\sqrt2)-\int\log|x/\sqrt2-t|\rho^H(t)\dd t-F_{V^H}
=J_{V^H}(x/\sqrt2).
\end{align*}
For $x<2$, $J_{V_0^L}(x)=J_{V^H}(x/\sqrt2)=+\infty$.
This proves the second identity in \eqref{eq:lcmax}. Recalling \eqref{eq:Hermax}, we obtain $J_{V_0^L}(2)=0$ and $J_{V_0^L}(x)>0$ for $x>2$. Then \ref{C1} follows from continuity. 

Condition \ref{C2} is automatic. For \ref{C3}, noting that
$\omega_{n,L}(y)\leq1$ for $y\geq0$, we derive from \eqref{e:Lagpos} that
$$
W_{n,L}(y)-\frac1{\beta n}\log\omega_{n,L}(y)
\geq W_{n,L}(y)\geq\frac{y^2}{4(1+y)}.
$$
Then Remark \ref{rem:mxpow} implies \ref{C3}.
By Proposition \ref{pr:maxldp}, we conclude (i) for all $\gamma\in[0,1]$.

To prove (ii), we take $D=\mathbb R_+$, $W_n(0)=+\infty$, $\omega_n(0)=1$ and
$$
W_n(x)=\frac{p}{2n}x-\frac{p-n+1}{2n}\log x,\qquad
\omega_n(x)=x^{-1},\qquad x>0.
$$
The limiting potential is $\widehat V_\gamma^L$ and
$\mathcal O=(0,\infty)$.
The same endpoint calculation as in (i) verifies \eqref{eq:mxend} at $0$ and \ref{C2} follows. The remaining assumptions of Proposition \ref{pr:maxldp} can be checked as in (i) for $\gamma>0$.
Therefore, $\lambda_{\max}$ satisfies an LDP with speed $\bt n$ and good rate function $J_{\widehat V_\gamma^L}$ given in \eqref{eq:mxrate}. 
By the same computations leading to \eqref{e:Lagrate}, we get
\begin{align*}
J_{\widehat V_\gamma^L}(x)
=\frac1{2\gamma}\int_{\gamma_+}^x
\frac{\sqrt{(t-\gamma_-)(t-\gamma_+)}}{t}\dd t,
\qquad x\geq\gamma_+,
\end{align*}
and $J_{\widehat V_\gamma^L}(x)=+\infty$ for $x<\gamma_+$.
This is \eqref{eq:Lagint} and (ii) is proven.

For (iii), let $0<\gamma<1$ and replace $\lambda_i$ by $-\lambda_i$.
By the same argument as in (ii), $\max_i(-\lambda_i)=-\lambda_{\min}$ satisfies an LDP with speed $\beta n$ and good rate function $x\mapsto J_{\widehat V_\gamma^L}^{-}(-x)$.
The contraction principle then proves (iii), with rate function given by \eqref{eq:lomin}.
\end{proof}

The following justifies the lower tail probability of the largest eigenvalue of Laguerre ensemble.
\begin{proof}[Proof of Example \ref{ex:Lag2}]
 The proof of Example \ref{x:lagemp} verifies Assumption
\ref{ass:ldp} for the rescaled density and Assumption \ref{ass:d} for
the unconstrained and constrained minimizers. Proposition
\ref{pr:lowgap} therefore applies and gives
\eqref{eq:lc2ldp} and
\begin{align}\label{eq:defigg}
 I_{V_\gamma^L}(\mu_{\zeta_1,\zeta_2}^L)
 =\frac12\int V_\gamma^L\dd\mu_{\zeta_1,\zeta_2}^L
 +\frac12F_{\zeta_1,\zeta_2}-C_{V_\gamma^L},	
\end{align}
where $F_{\zeta_1,\zeta_2}$ is characterized by
\eqref{eq:fzeta}.
If $\zeta_1\leq\sqrt\gamma-2$ and
$\zeta_2\geq\sqrt\gamma+2$, the unconstrained minimizer is admissible; hence $[s,t]=[\sqrt\gamma-2,\sqrt\gamma+2]$ and the rate is zero. We consider the other constrained cases below.

We first assume that $0<\gamma\leq1$.  Since
$$
 (V_\gamma^L)'(x)=\frac1{2\sqrt\gamma}
 -\frac{1-\gamma}{2\sqrt\gamma(1+\sqrt\gamma x)},
$$
Lemma \ref{remr}, applied with $D=[\zeta_1,\zeta_2]$, shows that the
support is a single interval $[s,t]\subset[\zeta_1,\zeta_2]$. 
Suppose that $0<\gamma<1$.
If only the lower constraint is active, then
$s=\zeta_1$, while the right endpoint is determined by
$$
 1
 =\frac1\pi\int_{\zeta_1}^t(V_\gamma^L)'(x)
   \sqrt{\frac{x-\zeta_1}{t-x}}\dd x
 =\frac{t-\zeta_1}{4\sqrt\gamma}
 -\frac{1-\gamma}{2\gamma}
 \left(1-\sqrt{\frac{1+\sqrt\gamma\zeta_1}
 {1+\sqrt\gamma t}}\right).
$$
Thus $w=\sqrt{1+\sqrt\gamma t}$ satisfies
\begin{equation*}
 \mathcal P_{\gamma,\zeta_1}^{R}(w)
 :=w^3-(3+2\gamma+\sqrt\gamma\zeta_1)w
 +2(1-\gamma)\sqrt{1+\sqrt\gamma\zeta_1}=0.
\end{equation*}
We have
\begin{align*}
 \mathcal P_{\gamma,\zeta_1}^{R}(-\infty)&<0,
 &\mathcal P_{\gamma,\zeta_1}^{R}(0)&>0,\\
 \mathcal P_{\gamma,\zeta_1}^{R}
 \bigl(\sqrt{1+\sqrt\gamma\zeta_1}\bigr)
 &=-4\gamma\sqrt{1+\sqrt\gamma\zeta_1}<0,
 &\mathcal P_{\gamma,\zeta_1}^{R}(+\infty)&>0.
\end{align*}
Hence the largest root is the unique root for which $t>\zeta_1$, and
\begin{equation}\label{eq:lrootr}
 1+\sqrt\gamma b_\gamma(\zeta_1)
 =\frac43(3+2\gamma+\sqrt\gamma\zeta_1)\cos^2 \left[
 \frac13\arccos \left(
 -\frac{3\sqrt3(1-\gamma)\sqrt{1+\sqrt\gamma\zeta_1}}
{(3+2\gamma+\sqrt\gamma\zeta_1)^{3/2}}\right)\right].
\end{equation}
If only the upper constraint is active, then $t=\zeta_2$, while the left endpoint $s=a_\gamma(\zeta_2)\geq\zeta_1$ is determined by
$$
 -1
 =\frac1\pi\int_s^{\zeta_2}(V_\gamma^L)'(x)
   \sqrt{\frac{\zeta_2-x}{x-s}}\dd x
=\frac{\zeta_2-s}{4\sqrt\gamma}
 -\frac{1-\gamma}{2\gamma}
 \left(\sqrt{\frac{1+\sqrt\gamma\zeta_2}
 {1+\sqrt\gamma s}}-1\right).
$$
Consequently, $w=\sqrt{1+\sqrt\gamma s}$ satisfies
\begin{equation*}
 w^3-(3+2\gamma+\sqrt\gamma\zeta_2)w
 +2(1-\gamma)\sqrt{1+\sqrt\gamma\zeta_2}=0.
\end{equation*}
Its roots lie in
$(-\infty,0)$, $(0,\sqrt{1+\sqrt\gamma\zeta_2})$, and
$(\sqrt{1+\sqrt\gamma\zeta_2},\infty)$, respectively.  The smaller
positive root is therefore the unique root for which $s<\zeta_2$, and
\begin{equation}\label{eq:lrootl}
 1+\sqrt\gamma a_\gamma(\zeta_2)
 =\frac{4(3+2\gamma+\sqrt\gamma\zeta_2)}{3}\cos^2 \left[
 \frac13\arccos \left(
 -\frac{3\sqrt3(1-\gamma)\sqrt{1+\sqrt\gamma\zeta_2}}
 {(3+2\gamma+\sqrt\gamma\zeta_2)^{3/2}}\right)
 -\frac{2\pi}{3}\right].
\end{equation}
Thus, for $0<\gamma<1$, the lower wall alone is active when $\zeta_1>\sqrt\gamma-2$ and $\zeta_2\geq b_\gamma(\zeta_1)$ and in this case $[s,t]=[\zeta_1,b_\gamma(\zeta_1)]$. Similarly, the upper wall alone is active if $\zeta_2<\sqrt\gamma+2$ and $\zeta_1\leq a_\gamma(\zeta_2)$ and then $[s,t]=[a_\gamma(\zeta_2),\zeta_2]$.
In the remaining constrained cases, both walls are active and $[s,t]=[\zeta_1,\zeta_2]$.

Now let $\gamma=1$. Since $V_1^L(x)=x/2$ is increasing, the lower edge of the support is $s=\max\{\zeta_1,-1\}$. The equation for a free right endpoint reduces to $(t-s)/4=1$, yielding that
$
t=\min\{\zeta_2,s+4\}.
$
Thus $a_1(z)=-1$ and $b_1(z)=z+4$, where the latter also follows from \eqref{eq:lrootr}. Together with the unconstrained case, this proves \eqref{eq:lcsupp} for $0<\gamma\leq1$. 

We next compute the density.  Let $m=(s+t)/2$ and $h=(t-s)/2$.
Lemma \ref{remr} gives
$$
 \rho_{\gamma;\zeta_1,\zeta_2}^L(x)
 =\frac{\pi-\ell(x)}
 {\pi^2\sqrt{(t-x)(x-s)}},\qquad  \ell(x)
 =\pv\int_s^t\sqrt{(t-y)(y-s)}
   \frac{(V_\gamma^L)'(y)}{x-y}\dd y.
$$
By Lemma \ref{le:Ia},
\begin{align*}
 \ell(x)
 &=\frac{\pi(x-m)}{2\sqrt\gamma}
 -\frac{\pi(1-\gamma)}{2\gamma}
 \left(1-\frac{\sqrt{(1+\sqrt\gamma s)(1+\sqrt\gamma t)}}
 {1+\sqrt\gamma x}\right),
\end{align*}
where we used
\begin{align*}
 \pv\int_s^t\frac{\sqrt{(t-y)(y-s)}}{x-y}\dd y
 &=\pi(x-m),\\
 \pv\int_s^t
 \frac{\sqrt{(t-y)(y-s)}}{(1+\sqrt\gamma y)(x-y)}\dd y
 &=\frac{\pi}{\sqrt\gamma}
 \left(1-\frac{\sqrt{(1+\sqrt\gamma s)(1+\sqrt\gamma t)}}
 {1+\sqrt\gamma x}\right).
\end{align*}
Substituting this expression for $\ell$ gives the density formula \eqref{eq:lcden}.

Let us evaluate \eqref{eq:defigg}. We follow the convention that whenever a coefficient is zero, then we regard this term zero without any further computations.  Using the density \eqref{eq:lcden}, we have
\begin{align*}
 &\int V_\gamma^L(x)
 \mu_{\zeta_1,\zeta_2}^L(\dd x)\\
= & \int_s^t
 \frac{\sqrt\gamma x-(1-\gamma)\log(1+\sqrt\gamma x)}
 {4\pi\gamma^2 \sqrt{(t-x)(x-s)}}
 \bigg(
 1+\gamma+\sqrt\gamma(m-x)
 -\frac{(1-\gamma)
 \sqrt{(1+\sqrt\gamma s)(1+\sqrt\gamma t)}}
 {1+\sqrt\gamma x}
 \bigg)\dd x\\
 =&
 \frac1{4\gamma^2}\bigg(
 2(1+\sqrt\gamma m)-2\gamma-\frac{\gamma h^2}{2}
 -2(1-\gamma)
 \sqrt{(1+\sqrt\gamma s)(1+\sqrt\gamma t)}\\
 &
 -2(1-\gamma^2)
 \log\frac{\sqrt{1+\sqrt\gamma s}
 +\sqrt{1+\sqrt\gamma t}}2
 +2(1-\gamma)^2
 \log\frac{2\sqrt{(1+\sqrt\gamma s)(1+\sqrt\gamma t)}}
 {\sqrt{1+\sqrt\gamma s}+\sqrt{1+\sqrt\gamma t}}
 \bigg),
\end{align*}
where the last equality follows from integration by parts and the identity from Lemma \ref{le:intlog}, 
$$
 \frac1\pi\int_s^t
 \frac{\log(1+\sqrt\gamma x)}
 {(1+\sqrt\gamma x)\sqrt{(t-x)(x-s)}}\dd x
 =
 \frac2{\sqrt{(1+\sqrt\gamma s)(1+\sqrt\gamma t)}}
 \log\frac{2\sqrt{(1+\sqrt\gamma s)(1+\sqrt\gamma t)}}
 {\sqrt{1+\sqrt\gamma s}+\sqrt{1+\sqrt\gamma t}}.
$$
Evaluating \eqref{eq:fzeta} at $x=t$ and using \eqref{eq:lcden}, we obtain
\begin{align*}
F_{\zeta_1,\zeta_2}
&=V_\gamma^L(t)
-\int_s^t\log(t-x)\mu_{\zeta_1,\zeta_2}^L(\dd x)\\
&=\frac{m}{2\sqrt\gamma}-\log h+\frac{\log2}{\gamma}
-\frac{1-\gamma}{\gamma}
\log\bigl(\sqrt{1+\sqrt\gamma s}+\sqrt{1+\sqrt\gamma t}\bigr),
\end{align*}
where we used Lemma \ref{le:intlog} in the second equality and $t=m+h$ in the last.
Plugging in the constant $C_{V_\gamma^L}$ from \eqref{eq:lagfun} and the above expressions for
$\int V_\gamma^L\dd\mu_{\zeta_1,\zeta_2}^L$ and
$F_{\zeta_1,\zeta_2}$ into \eqref{eq:defigg} yields 
$ I_{V_\gamma^L}(\mu_{\zeta_1,\zeta_2}^L)
 =\Phi_\gamma^L(s,t). $
With $\Phi_\gamma^L$ given by \eqref{eq:clphi}, this proves \eqref{eq:lcrate} for $0<\gamma\leq1$, with the value at $\gamma=1$ understood by continuity.

It remains to consider $\gamma=0$, where $V_0^L(x)=x^2/4$.  For
$\mu\in\mathcal M_1([\zeta_1,\zeta_2])$, define
$\nu(A)=\mu(\sqrt2 A)$ for every Borel set $A\subset\mathbb R$.  Then
$\nu\in\mathcal M_1([\zeta_1/\sqrt2,\zeta_2/\sqrt2])$ and
$$
 \int V_0^L\dd\mu=\int V^H\dd\nu,
 \qquad
 \Sigma(\nu)=\Sigma(\mu)-\frac12\log2.
$$
Since $C_{V_0^L}=3/8$ and
$C_{V^H}=3/8+(1/4)\log2$, we find
\begin{align*}
 I_{V_0^L}(\mu)
 &=\int V^H\dd\nu-\frac12
   \left(\Sigma(\nu)+\frac12\log2\right)-\frac38
 =I_{V^H}(\nu).
\end{align*}
Therefore, \eqref{eq:lcon} is equivalent  to \eqref{eq:hcon} with constraints $\zeta_1/\sqrt2$ and $\zeta_2/\sqrt2$.  By uniqueness,
$$
 I_{V_0^L}(\mu_{\zeta_1,\zeta_2}^L)
 =I_{V^H}(\mu_{\zeta_1/\sqrt2,\zeta_2/\sqrt2}^H),
 \qquad
 \rho_{0;\zeta_1,\zeta_2}^L(x)=\frac1{\sqrt2}
 \rho^H_{\zeta_1/\sqrt2,\zeta_2/\sqrt2}
 \left(\frac{x}{\sqrt2}\right).
$$
Since the corresponding support of the Hermite ensemble is $[s/\sqrt2,t/\sqrt2]$, we have
$a_0(z)=\sqrt2a(z/\sqrt2)$ and
$b_0(z)=\sqrt2b(z/\sqrt2)$.  These identities prove the $\gamma=0$ part of \eqref{eq:lcsupp},
\eqref{eq:lcden}, and
\eqref{eq:lcrate}.  Finally, substituting  $m/\sqrt2$ and $h/\sqrt2$ into the  rate function of the Hermite ensemble, \eqref{eq:lcgam0} follows from direct simplification.
\end{proof}

\subsection{The Jacobi ensemble} 
In this section, we justify Examples \ref{ex:JacempLDP} and \ref{ex:jac-max}, and prove Theorem \ref{thm:jaclow}.

\begin{proof}[Proof of Example \ref{ex:JacempLDP}]
To prove (i), it suffices to verify the assumptions of Theorem \ref{ldpln}. The constant $C_{V_{\gamma,\sigma}^J}$ is given in \cite{Ma23J}.
Take $D=\mathbb R$ and set
\begin{equation}\label{e:JacDn}
D_{n,J}=\varphi_n([0,1])
=\left[-\sqrt{\frac{p_1}{n}},\frac{p_2}{\sqrt{np_1}}\right].
\end{equation} 
For $y$ in the interior of $D_{n,J}$, define
\begin{align}
     W_{n, J}(y)
 &=-\frac{p_1-n+1}{2n}
   \log\left(1+\sqrt{\frac n{p_1}}y\right)
   -\frac{p_2-n+1}{2n}
   \log\left(1-\frac{\sqrt{np_1}}{p_2}y\right), \notag \\
 \omega_{n, J}(y)
 &=\frac1{\left(1+\sqrt{n/p_1} y\right)
           \left(1-\sqrt{np_1} y/p_2\right)}. \label{e:jacwgt}
\end{align}
Let $W_{n,J}=+\infty$ and $\omega_{n,J}=1$ elsewhere on $\mathbb R$. Let $(\lambda_1,\ldots,\lambda_n)$ have the $\beta$-Jacobi density \eqref{pdf} with potential \eqref{def:vnj}. Then  $(\varphi_n(\lambda_1),\ldots,\varphi_n(\lambda_n))$ has the weighted density \eqref{eq:wpdf} with $W_n=W_{n,J}$ and $\omega_n=\omega_{n,J}$.

Assumption \ref{A0} is immediate.
The local uniform convergence in \ref{A1}(ii) follows directly from \eqref{e:jacwgt} when $\gamma\sigma>0$ and from a Taylor expansion when $\gamma\sigma=0$.
Here the set $\mathcal O$ in Assumption \ref{ass:ldp} is the interior of $D_{\gamma,\sigma}=[-\gamma^{-1/2},(\sigma\sqrt\gamma)^{-1}]$.
Thus $D\setminus\mathcal O$ is empty when $\gamma=0$. For $\gamma>0$, it remains to prove
\begin{equation}\label{e:Jaclim}
\liminf_{n\to\infty}W_{n,J}(y_n)\geq V_{\gamma,\sigma}^J(y)
\qquad\text{whenever }y_n\to y\in D\setminus\mathcal O.
\end{equation}
If $y\notin D_{\gamma,\sigma}$, then $W_{n,J}(y_n)=+\infty$ for all sufficiently large $n$. It therefore suffices to consider $y_n$ in the interior of $D_{n,J}$ approaching a finite endpoint of $D_{\gamma,\sigma}$.
For $0<\gamma<1$ and $y_n\to-\gamma^{-1/2}$, $W_{n,J}(y_n)\to+\infty=V_{\gamma,\sigma}^J(-\gamma^{-1/2})$.
The same conclusion holds at $(\sigma\sqrt\gamma)^{-1}$ when $0<\gamma\sigma<1$.
For $\gamma=1$ and $y_n\to-1$, the first logarithmic term in
\eqref{e:jacwgt} is nonnegative for all sufficiently large $n$.
Hence
$$
\liminf_{n\to\infty}W_{n, J}(y_n)
\geq\lim_{n\to\infty}\left[
-\frac{p_2-n+1}{2n}
\log\left(1-\frac{\sqrt{np_1}}{p_2}y_n\right)\right]
=V_{1,\sigma}^J(-1),
$$
where the limit for $\sigma=0$ follows from
$\log(1-u)=-u+o(u)$.
For $\gamma\sigma=1$ and $y_n\to\sqrt\gamma$, the second logarithmic term is eventually nonnegative, so
$$
\liminf_{n\to\infty}W_{n, J}(y_n)
\geq\lim_{n\to\infty}\left[
-\frac{p_1-n+1}{2n}
\log\left(1+\sqrt{\frac n{p_1}}y_n\right)\right]
=V_{\gamma,1/\gamma}^J(\sqrt\gamma).
$$
Thus \eqref{e:Jaclim} holds in all cases, proving \ref{A1}(i).

Assumption \ref{A2} is trivial for $V_{\gamma,\sigma}^J$. If $\gamma\sigma>0$, then $W_{n,J}=+\infty$ outside $\cup_n D_{n,J}$. It suffices to consider $\gamma\sigma=0$.
 If $\gamma=0$, applying \eqref{e:logpos} and \eqref{e:logneg} to the two logarithms in \eqref{e:jacwgt}, we deduce
\begin{align}
 W_{n, J}(y)
 \geq\frac{n-1}{2n}\sqrt{\frac n{p_1}}
 \left(1-\frac{p_1}{p_2}\right)y
 +\frac{(p_1-n+1)y^2}{4p_1(1+\sqrt{n/p_1}|y|)}
 \geq\frac{c y^2}{1+|y|}-C \label{e:Jacgr0}
\end{align}
for all sufficiently large $n$ and every $y$ with $W_{n,J}(y)<\infty$.
If $\sigma=0<\gamma$, the left endpoints of $D_{n, J}$ remain bounded, whereas the right endpoints tend to $+\infty$.
 By \eqref{e:logneg}, for $y\geq0$ and all sufficiently large $n$,
$$ W_{n, J}(y)
 \geq-\frac{p_1-n+1}{2n}
 \log\left(1+\sqrt{\frac n{p_1}}y\right)
 +\frac{p_2-n+1}{2n}\frac{\sqrt{np_1}}{p_2}y
 \geq-C_1\log(1+y)+c_1y.$$
These bounds and the boundedness of each $D_{n,J}$ verify \ref{A2}.

By \eqref{e:jacwgt}, $\omega_{n,J}$ is positive, finite, and measurable, with $\log\omega_{n,J}=O(1)=o(\beta n)$ locally uniformly on $\mathcal O$. Thus \ref{A3}(i) holds. For \ref{A3}(ii), take $\theta_n=(\beta n)^{-1/2}$. Equations \eqref{eq:def-frac} and \eqref{e:jacwgt} give
\begin{align*}
 \mathcal R_{n, J}(\theta_n)
 ={}&\int_{D_{n, J}^n}
 \prod_{1\leq i<j\leq n}|y_i-y_j|^{\sqrt{\beta n}/n}
 \prod_{i=1}^n
 \left(1+\sqrt{\frac n{p_1}}y_i\right)^{r_{1,n}}
 \left(1-\frac{\sqrt{np_1}}{p_2}y_i\right)^{r_{2,n}}
  \dd\boldsymbol y\\
 ={}&C_{\mathcal R_{n, J}}\int_{(0,1)^n}
 \prod_{1\leq i<j\leq n}|x_i-x_j|^{\sqrt{\beta n}/n}
 \prod_{i=1}^n x_i^{r_{1,n}}(1-x_i)^{r_{2,n}}
  \dd\boldsymbol x\\
 ={}&C_{\mathcal R_{n, J}}\prod_{j=1}^n
 \frac{\Gamma\bigl(1+j\sqrt{\beta n}/(2n)\bigr)
       \Gamma\bigl(\frac{\sqrt{\beta n}(p_1-n+j)}{2n}\bigr)
       \Gamma\bigl(\frac{\sqrt{\beta n}(p_2-n+j)}{2n}\bigr)}
      {\Gamma\bigl(1+\sqrt{\beta n}/(2n)\bigr)
       \Gamma\bigl(\frac{\sqrt{\beta n}(p_1+p_2-n+j)}{2n}\bigr)},
\end{align*}
where $r_{i,n}=\frac{\sqrt{\beta n}(p_i-n+1)}{2n}-1, i=1,2,$ and
$$
C_{\mathcal R_{n, J}}=\left(\frac{p_1+p_2}{p_1}\right)^{
 \frac{\sqrt{\beta n}}2(p_1-n+1)}
 \left(\frac{p_1+p_2}{p_2}\right)^{
 \frac{\sqrt{\beta n}}2(p_2-n+1)}
 \left(\frac{p_2\sqrt{p_1/n}}{p_1+p_2}\right)^n
 \left(\frac{p_1+p_2}{\sqrt{np_1}}\right)^{
 \sqrt{\beta n}(n-1)/2}.
$$
Here the second equality follows from the change of variables $y_i=(p_1+p_2)x_i/\sqrt{np_1}-\sqrt{p_1/n}$, and the last from the Selberg integral formula; see, e.g., \cite{AGZ10}*{Eq.~(2.5.9)}. 
For all sufficiently large $n$, 
\begin{align*}
C_{LS}:={}&C\sum_{j=1}^n\Bigg[
 3+\left(-\log\frac{\sqrt{\beta n}(p_1-n+j)}{2n}\right)_+
 +\left(-\log\frac{\sqrt{\beta n}(p_2-n+j)}{2n}\right)_+\\
 &\qquad
 +\left(-\log\frac{\sqrt{\beta n}(p_1+p_2-n+j)}{2n}\right)_+
 \Bigg]\\
 \leq{}&
 C\sum_{j=1}^n\left(1+\log\frac{2n}{j}\right)
 ={}Cn+C\bigl(n\log(2n)-\log(n!)\bigr)
 \leq Cn.
\end{align*}
Using $-(1-u)\log(1-u)\leq u$ and $p_1,p_2\geq n$, we have
\begin{align*}
 J_{1,j}:={}&\frac{p_1-n+j}{n}\log\frac{p_1-n+j}{p_1}
 +\frac{p_2-n+j}{n}\log\frac{p_2-n+j}{p_2}\\
 &-\frac{p_1+p_2-n+j}{n}
  \log\frac{p_1+p_2-n+j}{p_1+p_2}
 \leq1-\frac jn,\\
 J_{2,j}:={}&-\log\frac{p_1-n+j}{p_1}
 -\log\frac{p_2-n+j}{p_2}
 +\log\frac{p_1+p_2-n+j}{p_1+p_2}
 \leq2\log\frac nj.
\end{align*}
Applying \eqref{bd:stir2} and substituting the expression for $C_{\mathcal R_{n,J}}$ now gives
\begin{align*}
 &\log\left[
 C_{\mathcal R_{n, J}}\prod_{j=1}^n
 \frac{
  \Gamma\bigl(\frac{\sqrt{\beta n}(p_1-n+j)}{2n}\bigr)
  \Gamma\bigl(\frac{\sqrt{\beta n}(p_2-n+j)}{2n}\bigr)}
 {\Gamma\bigl(\frac{\sqrt{\beta n}(p_1+p_2-n+j)}{2n}\bigr)}
 \right]\\
\leq {}& \log C_{\mathcal R_{n, J}}
 +\sum_{j=1}^n\Bigg[
 \left(\frac{\sqrt{\beta n}(p_1-n+j)}{2n}-\frac12\right)
 \log\frac{\sqrt{\beta n}(p_1-n+j)}{2n}
 -\frac{\sqrt{\beta n}(p_1-n+j)}{2n}\\
 &
 +\left(\frac{\sqrt{\beta n}(p_2-n+j)}{2n}-\frac12\right)
 \log\frac{\sqrt{\beta n}(p_2-n+j)}{2n}
 -\frac{\sqrt{\beta n}(p_2-n+j)}{2n}\\
 &
 -\left(\frac{\sqrt{\beta n}(p_1+p_2-n+j)}{2n}-\frac12\right)
 \log\frac{\sqrt{\beta n}(p_1+p_2-n+j)}{2n}\\
 &
 +\frac{\sqrt{\beta n}(p_1+p_2-n+j)}{2n}
 \Bigg]+C_{LS}\\
\leq{} & \frac{\sqrt{\beta n}(n-1)}4
 -\left(\frac{\sqrt{\beta n}(n-1)}4+\frac n2\right)
  \log\frac{\sqrt{\beta n}}2\\
 &
 +\left(\frac{\sqrt{\beta n}(n-1)}4+\frac n2\right)
   \log\frac{p_2}{p_1+p_2}
 +\frac{\sqrt{\beta n}}2\sum_{j=1}^nJ_{1,j}
 +\frac12\sum_{j=1}^nJ_{2,j}+Cn\\
\leq{}  & Cn(1+\sqrt{\beta n}).
\end{align*}
Together with \eqref{bd:gamj}, these estimates give, uniformly in $p_1,p_2\geq n$,
$$
 \log\mathcal R_{n, J}(\theta_n)
 \leq Cn(1+\sqrt{\beta n})\log(2+\sqrt{\beta n})
 =o(\beta n^2),
$$
as $\beta n\to\infty$. This proves \ref{A3}.

Every $\mu\in\mathcal M_1(\mathbb R)$ with
$I_{V_{\gamma,\sigma}^J}(\mu)<\infty$ satisfies
$\mu(D_{\gamma,\sigma})=1$.
Lemma \ref{le:symreg} with $D=\mathbb R$ and $E=D_{\gamma,\sigma}$
therefore verifies Assumption \ref{ass:d}.
 Theorem \ref{ldpln} now proves (i).

To prove (ii), take $D=[0,1]$ and
$$
 W_n(x)=-\frac{p_1-n+1}{2n}\log x
        -\frac{p_2-n+1}{2n}\log(1-x),
 \qquad \omega_n(x)=\frac1{x(1-x)},\qquad 0<x<1,
$$
with $W_n=+\infty$ and $\omega_n=1$ at the endpoints.
The change of variables $y=\varphi_n(x)$ relates this decomposition to \eqref{e:jacwgt}. Since $\gamma\sigma>0$, the preceding estimates also verify Assumption \ref{ass:ldp} for the unrescaled model on $[0,1]$. 
Lemma \ref{le:symreg} verifies Assumption \ref{ass:d} for every $\mu\in\mathcal M_1([0,1])$ with $I_{\widehat V_{\gamma,\sigma}^J}(\mu)<\infty$. 
Theorem \ref{ldpln} therefore proves (ii). Substituting $y=\varphi(x)$ in \eqref{defiv} gives \eqref{eq:joremp}.
\end{proof}

By \cite{Ma23J}, the unique unconstrained minimizer $\widetilde\nu_{\gamma,\sigma}$ of $I_{V_{\gamma,\sigma}^J}$ over $\mathcal M_1(D_{\gamma,\sigma})$ has density 
\begin{align}\label{eq:tildeh}
 \tilde{h}_{\gamma,\sigma}(x)=
\begin{cases}
\frac{\sigma \sqrt{\gamma}}{1+\sigma}
 \rho_{\gamma,\sigma}^J\left(\frac{\sigma}{1+\sigma}(\sqrt{\gamma} x+1)\right),
& 0<\sigma\gamma\le 1, \\
\sqrt{\gamma} \rho_\gamma^L(1+\sqrt{\gamma}x),
& \sigma=0, 0<\gamma\le 1, \\
\frac{1+\sigma}{2\pi}\sqrt{\frac{4}{1+\sigma}-x^2},
& \gamma=0, \sigma\ge 0,
\end{cases}
\end{align}
where $\rho_\gamma^L$ and $\rho_{\gamma,\sigma}^J$ are given in
\eqref{eq:hgam} and \eqref{eq:hgamsig}, respectively.  Its support is $[\widetilde u_-,\widetilde u_+]$, where
$
 \widetilde u_\pm
 =\frac{\sqrt\gamma(1-\sigma)
 \pm2\sqrt{1+\sigma-\sigma\gamma}}{1+\sigma}.
$
Let us now prove the LDP of the largest eigenvalue.

\begin{proof}[Proof of Example \ref{ex:jac-max}]
To prove (i), use \eqref{e:jacwgt} on $D=\mathbb R$, with $\mathcal O$ the interior of $D_{\gamma,\sigma}$.
The proof of Example \ref{ex:JacempLDP} verifies
Assumption \ref{ass:ldp} and Assumption \ref{ass:d} for every $\mu\in\mathcal M_1(\mathbb R)$ with $I_{V_{\gamma,\sigma}^J}(\mu)<\infty$.
By Proposition \ref{pr:maxldp}, it remains to verify Assumption \ref{ass:maximum}.
Condition \ref{C0} holds by assumption.

Suppose first that $0<\gamma\sigma<1$.
The density in \eqref{eq:tildeh} gives $b_{V_{\gamma,\sigma}^J}=\max\supp(\widetilde\nu_{\gamma,\sigma})
=\widetilde u_+$.
By \eqref{eq:mxrate},
\begin{equation}\label{e:Jacpot}
J_{V_{\gamma,\sigma}^J}(y)
=V_{\gamma,\sigma}^J(y)
-\int_{\widetilde u_-}^{\widetilde u_+}
\log|y-t| \tilde h_{\gamma,\sigma}(t)\dd t
-F_{V_{\gamma,\sigma}^J},
\qquad \widetilde u_+\leq y<\frac1{\sigma\sqrt\gamma}.
\end{equation}
Equations \eqref{eq:arcavg} and \eqref{eq:jpot}, together with
Lemma \ref{le:intlog}, give
\begin{align*}
F_{V_{\gamma,\sigma}^J}
&=\frac1\pi\int_{\widetilde u_-}^{\widetilde u_+}
\frac{V_{\gamma,\sigma}^J(t)}
{\sqrt{(\widetilde u_+-t)(t-\widetilde u_-)}}\dd t
-\log\frac{\widetilde u_+-\widetilde u_-}{4}
=\frac{1+\sigma}{2\sigma\gamma}\log(1+\sigma)
-\frac{1+\sigma-\sigma\gamma}{2\sigma\gamma}
\log(1+\sigma-\sigma\gamma).
\end{align*}
Differentiating \eqref{e:Jacpot} and using \eqref{eq:tildeh} and Lemma \ref{le:Ia}, we obtain
\begin{align*}
J'_{V_{\gamma,\sigma}^J}(y)
&=-\frac{1-\gamma}{2\sqrt\gamma(1+\sqrt\gamma y)}
+\frac{1-\sigma\gamma}{2\sqrt\gamma(1-\sigma\sqrt\gamma y)}
-\int_{\widetilde u_-}^{\widetilde u_+}
\frac{\tilde h_{\gamma,\sigma}(t)}{y-t}\dd t\\
&=\frac{1+\sigma}{2}
\frac{\sqrt{(y-\widetilde u_-)(y-\widetilde u_+)}}
{(1+\sqrt\gamma y)(1-\sigma\sqrt\gamma y)}>0,
\qquad \widetilde u_+<y<\frac1{\sigma\sqrt\gamma}.
\end{align*}
Since $J_{V_{\gamma,\sigma}^J}(\widetilde u_+)=0$, integration gives
\begin{equation}\label{e:Jacrate}
J_{V_{\gamma,\sigma}^J}(y)
=\frac{1+\sigma}{2}\int_{\widetilde u_+}^y
\frac{\sqrt{(t-\widetilde u_-)(t-\widetilde u_+)}}
{(1+\sqrt\gamma t)(1-\sigma\sqrt\gamma t)}\dd t,
\qquad \widetilde u_+\leq y<\frac1{\sigma\sqrt\gamma}.
\end{equation}
Together with $b_{V_{\gamma,\sigma}^J}\in\mathcal O$ and the boundary limits in \eqref{eq:jpot}, this verifies \ref{C1}.

To verify \ref{C2}, we check \eqref{eq:mxend} at the
two endpoints of $D_{\gamma,\sigma}$.
The endpoint convergence in \eqref{e:JacDn} ensures that the
integrand vanishes on every compact subset of
$\mathbb R\setminus D_{\gamma,\sigma}$ for all sufficiently large $n$.
For fixed $0<\delta<\min\{\gamma^{-1/2},(\sigma\sqrt\gamma)^{-1}\}$,
\eqref{e:jacwgt} gives
\begin{align}\label{e:Jacsbd}
&\limsup_{n\to\infty}\frac1{\beta n}
\log\int_{-\gamma^{-1/2}-\delta}^{-\gamma^{-1/2}+\delta}
e^{-\beta nW_{n,J}(y)}\omega_{n,J}(y)\dd y
\leq \limsup_{n\to\infty}\frac1{\beta n}
\log\int_{-\infty}^{-\gamma^{-1/2}+\delta}
e^{-\beta nW_{n,J}(y)}\omega_{n,J}(y)\dd y\notag\\
\leq{} & \limsup_{n\to\infty}\frac1{\beta n}\log\left[
\left(1+\frac{p_1}{p_2}\right)^{\beta(p_2-n+1)/2-1}
\int_{-\sqrt{p_1/n}}^{-\gamma^{-1/2}+\delta}
\left(1+\sqrt{\frac n{p_1}}y\right)^{\beta(p_1-n+1)/2-1}\dd y
\right]\notag\\
={} & \limsup_{n\to\infty}\frac1{\beta n}\log\left[
\frac{2\sqrt{p_1/n}}{\beta(p_1-n+1)}
\left(1+\frac{p_1}{p_2}\right)^{\beta(p_2-n+1)/2-1}
\left(1-\sqrt{\frac n{\gamma p_1}}+\sqrt{\frac n{p_1}}\,\delta\right)^{\beta(p_1-n+1)/2}
\right]\notag\\
\leq{} &  \frac{1-\gamma}{2\gamma}\log(\sqrt\gamma\,\delta)
+\frac{1-\gamma\sigma}{2\gamma\sigma}\log(1+\sigma),\notag\\
&\limsup_{n\to\infty}\frac1{\beta n}
\log\int_{(\sigma\sqrt\gamma)^{-1}-\delta}^{(\sigma\sqrt\gamma)^{-1}+\delta}
e^{-\beta nW_{n,J}(y)}\omega_{n,J}(y)\dd y \leq \limsup_{n\to\infty}\frac1{\beta n}
\log\int_{(\sigma\sqrt\gamma)^{-1}-\delta}^\infty
e^{-\beta nW_{n,J}(y)}\omega_{n,J}(y)\dd y\notag\\
\leq{} & \limsup_{n\to\infty}\frac1{\beta n}\log\left[
\left(1+\frac{p_2}{p_1}\right)^{(\beta(p_1-n+1)/2-1)_+}
\int_{(\sigma\sqrt\gamma)^{-1}-\delta}^{p_2/\sqrt{np_1}}
\left(1-\frac{\sqrt{np_1}}{p_2}y\right)^{\beta(p_2-n+1)/2-1}\dd y
\right]\notag\\
={} &\limsup_{n\to\infty}\frac1{\beta n}\log\left[
\frac{2p_2 \left(1+{p_2}/{p_1}\right)^{(\beta(p_1-n+1)/2-1)_+}}{\beta(p_2-n+1)\sqrt{np_1}}
\left(1-\frac{\sqrt{np_1}}{\sigma\sqrt\gamma p_2}
+\frac{\sqrt{np_1}}{p_2}\delta\right)^{\beta(p_2-n+1)/2}
\right]\notag\\
\leq{} & \frac{1-\gamma}{2\gamma}\log(1+\sigma^{-1})
+\frac{1-\gamma\sigma}{2\gamma\sigma}\log(\sigma\sqrt\gamma\,\delta),
\end{align}
where $a_+=\max\{a,0\}$. In the endpoint estimates we used
$\beta(p_2-n+1)/2-1\geq0$ for all sufficiently large $n$;
in the final inequalities we also used $p_i-n+1\geq1$ and
$\log^+(1/\beta)=o(\beta n)$ from \ref{C0}.
As $\delta\downarrow0$, the two bounds tend to
$-V_{\gamma,\sigma}^J(-\gamma^{-1/2})$ and
$-V_{\gamma,\sigma}^J((\sigma\sqrt\gamma)^{-1})$.
Thus \eqref{eq:mxend} holds at both endpoints, and
Remark \ref{re:mxendpt} gives \ref{C2}. 
Condition \ref{C3} follows from the uniform boundedness of the intervals $D_{n,J}$ in \eqref{e:JacDn}.

Next suppose that $\sigma=0<\gamma$.
Then $V_{\gamma,0}^J=V_\gamma^L$.
Their equilibrium measures and constants therefore coincide, so \eqref{eq:mxrate} gives $J_{V_{\gamma,0}^J}=J_{V_\gamma^L}.$
This proves the second identity in \eqref{eq:jcmax}, while Example \ref{ex:lag-max} verifies \ref{C1}.

For \ref{C2}, we prove \eqref{eq:mxend} at $-\gamma^{-1/2}$.
The finite $n$ left endpoint estimate in \eqref{e:Jacsbd} remains valid when $\sigma=0$, because $\beta(p_2-n+1)/2-1\geq0$ for all sufficiently large $n$.
Thus, for fixed $0<\delta<\gamma^{-1/2}$, the same estimate gives
\begin{align*}
&\limsup_{n\to\infty}\frac1{\beta n}
\log\int_{-\gamma^{-1/2}-\delta}^{-\gamma^{-1/2}+\delta}
e^{-\beta nW_{n,J}(y)}\omega_{n,J}(y)\dd y
\leq\frac{1-\gamma}{2\gamma}\log(\sqrt\gamma\,\delta)
+\frac1{2\gamma}.
\end{align*}
As in the Laguerre case, letting $\delta\downarrow0$ verifies \eqref{eq:mxend}. Remark \ref{re:mxendpt} therefore gives \ref{C2}.

The left endpoints $-\sqrt{p_1/n}$ are bounded below, so \ref{C3} reduces to a right tail estimate.
For all large $n$ and $0\leq y<p_2/\sqrt{np_1}$,
\eqref{e:jacwgt} and \eqref{e:logneg} give
\begin{align*}
W_{n,J}(y)-\frac1{\beta n}\log\omega_{n,J}(y)
\geq -\frac{p_1-n+1}{2n}
\log\left(1+\sqrt{\frac n{p_1}}y\right)+\left(\frac{p_2-n+1}{2n}-\frac1{\beta n}\right)
\frac{\sqrt{np_1}}{p_2}y 
\geq cy-C.
\end{align*}
Since $W_{n,J}=+\infty$ for $y\geq p_2/\sqrt{np_1}$, Remark \ref{rem:mxpow} gives \ref{C3}.

If $\gamma=0$, then $\mathcal O=D=\mathbb R$ and $V_{0,\sigma}^J(y)=V^H(\sqrt{(1+\sigma)/2} y)$.
The densities in \eqref{eq:tildeh} and \eqref{eq:semic} give $b_{V_{0,\sigma}^J}=2/\sqrt{1+\sigma}$ and
$$
\int_{\mathbb R}\log|y-t| \tilde h_{0,\sigma}(t)\dd t
=\int_{\mathbb R}
\log\left|\sqrt{\frac{1+\sigma}{2}} y-u\right|\rho^H(u)\dd u
-\frac12\log\frac{1+\sigma}{2}.
$$
Evaluating \eqref{eq:vlogmu} at $y=b_{V_{0,\sigma}^J}$ therefore gives $F_{V_{0,\sigma}^J}=F_{V^H}+\frac12\log((1+\sigma)/2)$.
Hence \eqref{eq:mxrate} yields, for $y\geq b_{V_{0,\sigma}^J}$,
\begin{align*}
J_{V_{0,\sigma}^J}(y)
&=V_{0,\sigma}^J(y)
-\int_{\mathbb R}\log|y-t| \tilde h_{0,\sigma}(t)\dd t
-F_{V_{0,\sigma}^J}\\
&=V^H\left(\sqrt{\frac{1+\sigma}{2}} y\right)
-\int_{\mathbb R}
\log\left|\sqrt{\frac{1+\sigma}{2}} y-u\right|\rho^H(u)\dd u
-F_{V^H}
=J_{V^H}\left(\sqrt{\frac{1+\sigma}{2}} y\right).
\end{align*}
For $y<b_{V_{0,\sigma}^J}$, $J_{V_{0,\sigma}^J}(y)
=J_{V^H}(\sqrt{(1+\sigma)/2} y)=+\infty$.
This proves the last identity in \eqref{eq:jcmax}.
By \eqref{eq:Hermax}, $J_{V_{0,\sigma}^J}$ is zero at $b_{V_{0,\sigma}^J}$ and positive above it.
Together with the continuity of $V_{0,\sigma}^J$, this verifies \ref{C1}.

Condition \ref{C2} is automatic since $D\setminus\mathcal O=\varnothing$.
For \ref{C3}, the coefficients $\frac{p_i-n+1}{2n}-\frac1{\beta n}$, $i=1,2$, are positive for all sufficiently large $n$.
Applying \eqref{e:logpos} and \eqref{e:logneg} to the two logarithmic terms in \eqref{e:jacwgt} gives
\begin{align*}
 W_{n,J}(y)-\frac1{\beta n}\log\omega_{n,J}(y) 
\geq{} & \frac{(n-1+2/\beta)(p_2-p_1)y}{2p_2\sqrt{np_1}}
+\frac{(p_1-n+1-2/\beta)y^2}
{4p_1\bigl(1+\sqrt{n/p_1}|y|\bigr)}\\
\geq{} & c\frac{y^2}{1+|y|}-C
\geq\frac c2|y|-C,\qquad |y|\geq1,
\end{align*}
for all sufficiently large $n$ and $y$ in the interior of $D_{n,J}$.
Here the coefficients of the linear and quadratic terms tend to $0$ and $1/4$, respectively, as in \eqref{e:Jacgr0}.
The bound also holds outside this interval, where $W_{n,J}=+\infty$.
Remark \ref{rem:mxpow} therefore gives \ref{C3}. 
Combining the preceding verifications with Proposition \ref{pr:maxldp} proves (i) for all $\gamma\in[0,1]$, $\gamma \sigma \in [0,1).$

To prove (ii), take $D=[0,1]$ and use
$$
W_n(x)=-\frac{p_1-n+1}{2n}\log x
       -\frac{p_2-n+1}{2n}\log(1-x),\qquad
\omega_n(x)=\frac1{x(1-x)},\qquad 0<x<1,
$$
with $W_n=+\infty$ and $\omega_n=1$ at $0$ and $1$.
The limiting potential is $\widehat V_{\gamma,\sigma}^J$ and $\mathcal O=(0,1)$.
As in (i) for $\gamma\sigma>0$, the assumptions of Proposition \ref{pr:maxldp} hold. Thus $\lambda_{\max}$ satisfies an LDP with speed $\bt n$ and good rate function $J_{\widehat V_{\gamma,\sigma}^J}$ of the form \eqref{eq:mxrate}. Using the same calculation for \eqref{e:Jacrate}, we arrive at
$$
J_{\widehat V_{\gamma,\sigma}^J}(x)
=\frac{1+\sigma}{2\sigma\gamma}\int_{u_+}^x
\frac{\sqrt{(t-u_-)(t-u_+)}}{t(1-t)}\dd t,
\qquad u_+\leq x<1,
$$
and $J_{\widehat V_{\gamma,\sigma}^J}(x)=+\infty$ otherwise. We have proved (ii).
For $\widetilde u_+\leq y<1/(\sigma\sqrt\gamma)$, plugging $t=\sigma(1+\sqrt\gamma v)/(1+\sigma)$ in \eqref{eq:jomax} and using \eqref{e:Jacrate}, we find
\begin{align*}
J_{\widehat V_{\gamma,\sigma}^J}
\left(\frac{\sigma(1+\sqrt\gamma y)}{1+\sigma}\right)
=\frac{1+\sigma}{2}\int_{\widetilde u_+}^y
\frac{\sqrt{(v-\widetilde u_-)(v-\widetilde u_+)}}
{(1+\sqrt\gamma v)(1-\sigma\sqrt\gamma v)}\dd v
=J_{V_{\gamma,\sigma}^J}(y).
\end{align*}
Both sides are $+\infty$ outside this range. The first identity in \eqref{eq:jcmax} follows.

For (iii), from the density \eqref{def:vnj} we see that the change of variables $\lambda_i\mapsto1-\lambda_i$ transforms the $\beta$-Jacobi ensemble with parameters $(n,p_1,p_2)$ into the same ensemble with parameters $(n,p_2,p_1)$.
The limiting parameters satisfy
$$
\frac n{p_2}\to\gamma\sigma\in(0,1],
\qquad
\frac{p_2}{p_1}\to\sigma^{-1},
\qquad
(\gamma\sigma)\sigma^{-1}=\gamma\in(0,1).
$$
By (ii), $\max_i(1-\lambda_i)=1-\lambda_{\min}$ satisfy an LDP  with speed $\beta n$ and good rate function $J_{\widehat V_{\gamma\sigma,\sigma^{-1}}^J}$.
Applying the contraction principle to $x\mapsto1-x$ yields the rate function $x\mapsto J_{\widehat V_{\gamma\sigma,\sigma^{-1}}^J}(1-x)$ for $\lambda_{\min}$ and (iii) follows.
\end{proof}

Next, we define the endpoint functions used in the following. For $\gamma\sigma>0$, let
\begin{align}\label{eq:jacH}
 \mathcal H_z(w)
 &=
 \frac{1-\gamma}{2\gamma}
 \left(
 \sqrt{\frac{1+\sqrt\gamma z}{1+\sqrt\gamma w}}-1
 \right)
 +\frac{1-\sigma\gamma}{2\sigma\gamma}
 \left(
 \sqrt{\frac{1-\sigma\sqrt\gamma z}
 {1-\sigma\sqrt\gamma w}}-1
 \right).
\end{align}
When $z<\widetilde u_+$ and $0<\gamma<1$, define
$a_{\gamma,\sigma}(z)$ as the unique
$a\in(-\gamma^{-1/2},z)$ satisfying $\mathcal H_z(a)=1$; If
$\gamma=1$, let $a_{1,\sigma}(z)=-1$.  Similarly, for
$z>\widetilde u_-$ and $0<\sigma\gamma<1$, define
$b_{\gamma,\sigma}(z)$ as the unique
$b\in(z,(\sigma\sqrt\gamma)^{-1})$ satisfying
$\mathcal H_z(b)=1$; when $\sigma\gamma=1$, set
$b_{\gamma,\sigma}(z)=(\sigma\sqrt\gamma)^{-1}$.
The boundary functions $a_{0,\sigma}$, $b_{0,\sigma}$,
$a_{\gamma,0}$, and $b_{\gamma,0}$ are defined by continuous extension
as $\gamma\downarrow0$ or $\sigma\downarrow0$, respectively.  Specifically, for $\gamma=0$ and $\sigma\geq0$, writing
$c_\sigma=\sqrt{(1+\sigma)/2}$, we find
\begin{align}\label{eq:jacher}
 a_{0,\sigma}(z)
 =c_\sigma^{-1}a(c_\sigma z)
 =\frac{z-2\sqrt{z^2+12/(1+\sigma)}}3, \quad
 b_{0,\sigma}(z)
 =c_\sigma^{-1}b(c_\sigma z)
 =\frac{z+2\sqrt{z^2+12/(1+\sigma)}}3,
\end{align}
where $a$ and $b$ are the Hermite endpoint functions in Example \ref{ex:her2};  for $\sigma=0$ and $0<\gamma\leq1$, we have
\begin{equation}\label{eq:jaclag}
 a_{\gamma,0}(z)=a_\gamma(z),\qquad
 b_{\gamma,0}(z)=b_\gamma(z),
\end{equation}
where $a_\gamma$ and $b_\gamma$ are defined in Example 
\ref{ex:Lag2}. Indeed, for fixed $z,w$ in the interior of the limiting domain, direct expansion gives locally uniformly in $w$,
\begin{align*}
 \mathcal H_z(w)
 & \to
 \frac{1+\sigma}{16}(w-z)(z+3w),
 &&\gamma\downarrow0,\quad \sigma>0,\\
 \mathcal H_z(w)
 & \to
 \frac{1-\gamma}{2\gamma}
 \left(
 \sqrt{\frac{1+\sqrt\gamma z}{1+\sqrt\gamma w}}-1
 \right)
 +\frac{w-z}{4\sqrt\gamma},
 &&\sigma\downarrow0,\quad 0<\gamma\leq1.
\end{align*}
Setting the first limit equal to one gives
$
 w=\frac{z\pm2\sqrt{z^2+12/(1+\sigma)}}3,
$
namely $c_\sigma^{-1}a(c_\sigma z)$ and
$c_\sigma^{-1}b(c_\sigma z)$, while the second limiting function with value at one is precisely the Laguerre endpoint equation.  Uniqueness of the solutions
on the two sides of $z$ proves the stated identities.

\begin{proposition}\label{p:jacmin}
Under the assumptions of Theorem \ref{thm:jaclow}, let
 $\mu_{\zeta_1,\zeta_2}^J
 =\argmin_{\mu\in\mathcal M_1([\zeta_1,\zeta_2])}
 I_{V_{\gamma,\sigma}^J}(\mu).$
Then the following statements hold.
\begin{enumerate}
\item[(i)] The support of $\mu_{\zeta_1,\zeta_2}^J$ is an interval
$[s,t]$, given by
\begin{equation}\label{eq:jacsup}
 [s,t]=
 \begin{cases}
 [\widetilde u_-,\widetilde u_+],
 &\zeta_1\leq\widetilde u_-,\quad
  \zeta_2\geq\widetilde u_+,\\
 [\zeta_1,b_{\gamma,\sigma}(\zeta_1)],
 &\zeta_1>\widetilde u_-,\quad
  \zeta_2\geq b_{\gamma,\sigma}(\zeta_1),\\
 [a_{\gamma,\sigma}(\zeta_2),\zeta_2],
 &\zeta_2<\widetilde u_+,\quad
  \zeta_1\leq a_{\gamma,\sigma}(\zeta_2),\\
 [\zeta_1,\zeta_2],&\text{otherwise}.
 \end{cases}
\end{equation}

\item[(ii)] If $\gamma\sigma=0$, then, for $x\in(s,t)$,
$\mu_{\zeta_1,\zeta_2}^J$ has density
$$
 \rho_{\gamma,\sigma;\zeta_1,\zeta_2}^J(x)=
 \begin{cases}
  
 \sqrt{\frac{1+\sigma}{2}} 
 \rho^H_{\sqrt{(1+\sigma)/2} \zeta_1, 
          \sqrt{(1+\sigma)/2} \zeta_2}
 \left(\sqrt{\frac{1+\sigma}{2}} x\right),
 &\gamma=0,\quad\sigma\geq0,\\
 \rho_{\gamma;\zeta_1,\zeta_2}^L(x),
 &\sigma=0,\quad0<\gamma\leq1.
 \end{cases}
$$
Here $\rho^H_{\zeta_1,\zeta_2}$ and
$\rho_{\gamma;\zeta_1,\zeta_2}^L$ are the densities
defined in Examples \ref{ex:her2} and
\ref{ex:Lag2}, respectively.

\item[(iii)] If $\gamma\sigma>0$, then
$\mu_{\zeta_1,\zeta_2}^J$ has density
\begin{align}
     \rho_{\gamma,\sigma;\zeta_1,\zeta_2}^J(x)
 ={}&\frac{1}{\pi\sqrt{(t-x)(x-s)}}
 \Bigg[1+\frac{1-\gamma}{2\gamma}
 \left(1-\frac{\sqrt{(1+\sqrt\gamma s)(1+\sqrt\gamma t)}}
 {1+\sqrt\gamma x}\right) \notag
 \\
 &+\frac{1-\sigma\gamma}{2\sigma\gamma}
 \left(1-\frac{\sqrt{(1-\sigma\sqrt\gamma s)
 (1-\sigma\sqrt\gamma t)}}
 {1-\sigma\sqrt\gamma x}\right)\Bigg],
 \qquad x\in(s,t). \label{eq:jacden}
\end{align}
\end{enumerate}
If $\zeta_1\leq\widetilde u_-$ and
$\zeta_2\geq\widetilde u_+$, then the density in part~(ii) when
$\gamma\sigma=0$, or in part~(iii) when $\gamma\sigma>0$,
coincides with the unconstrained density $\tilde h_{\gamma,\sigma}$
in \eqref{eq:tildeh}.
\end{proposition}

\begin{proof}
If $\zeta_1\leq\widetilde u_-$ and
$\zeta_2\geq\widetilde u_+$, then
$
 \widetilde\nu_{\gamma,\sigma}\in
 \mathcal M_1([\zeta_1,\zeta_2]).
$
Since $\widetilde\nu_{\gamma,\sigma}$ is the unconstrained
equilibrium measure, it is also the unique minimizer of the constrained
problem.  Hence $[s,t]=[\widetilde u_-,\widetilde u_+]$, which proves
the first case in part~(i).  Assume that
$\zeta_1>\widetilde u_-$ or $\zeta_2<\widetilde u_+$ below.

We first consider $\gamma=0$ and
$\sigma\geq0$.  Let $c_\sigma=\sqrt{(1+\sigma)/2}$.  Then
$V_{0,\sigma}^J(x)=V^H(c_\sigma x)$.  For any
$\mu\in\mathcal M_1([\zeta_1,\zeta_2])$, define
$
 \nu(A)=\mu(c_\sigma^{-1}A)
$
for every Borel set $A\subset\mathbb R$.  Then
$\nu\in\mathcal M_1([c_\sigma\zeta_1,c_\sigma\zeta_2])$ and
\begin{align}\label{eq:v0sjh}
 \int V_{0,\sigma}^J\dd\mu=\int V^H\dd\nu,
 \qquad
 \Sigma(\nu)=\Sigma(\mu)+\log c_\sigma.	
\end{align}
Thus the constrained problem \eqref{eq:con-J} for
$V_{0,\sigma}^J$ on $[\zeta_1,\zeta_2]$ is equivalent to the Hermite
constrained problem \eqref{eq:hcon} on
$[c_\sigma\zeta_1,c_\sigma\zeta_2]$.  Let $[s_H,t_H]$ be the support
of $\mu^H_{c_\sigma\zeta_1,c_\sigma\zeta_2}$.  By uniqueness and
Example \ref{ex:her2},
$$
 \rho_{0,\sigma;\zeta_1,\zeta_2}^J(x)
 =c_\sigma
 \rho^H_{c_\sigma\zeta_1,c_\sigma\zeta_2}(c_\sigma x),
 \qquad [s,t]=[c_\sigma^{-1}s_H,c_\sigma^{-1}t_H].
$$
We have
$ c_\sigma^{-1}[-\sqrt2,\sqrt2]
 = \left[-2/\sqrt{1+\sigma}, 2/\sqrt{1+\sigma}\right]
 = [\widetilde u_-,\widetilde u_+],$
and \eqref{eq:jacher} gives
\begin{align*}
 c_\sigma^{-1}[c_\sigma\zeta_1,b(c_\sigma\zeta_1)]
 =[\zeta_1,b_{0,\sigma}(\zeta_1)], \quad 
 c_\sigma^{-1}[a(c_\sigma\zeta_2),c_\sigma\zeta_2]
=[a_{0,\sigma}(\zeta_2),\zeta_2], \quad
 c_\sigma^{-1}[c_\sigma\zeta_1,c_\sigma\zeta_2]
 =[\zeta_1,\zeta_2].
\end{align*}
The four Hermite support alternatives are mapped exactly to those in
\eqref{eq:jacsup}.  This proves \eqref{eq:jacsup}
when $\gamma=0$.

Next assume that $\sigma=0$ and $0<\gamma\leq1$.  Since
$V_{\gamma,0}^J=V_\gamma^L$ on
$D_{\gamma,0}=D_\gamma$, the Jacobi constrained problem
\eqref{eq:con-J} is identical to the  Laguerre constrained problem \eqref{eq:lcon} with the same walls.  
Therefore, $\mu_{\zeta_1,\zeta_2}^J=\mu_{\zeta_1,\zeta_2}^L,$ and $\rho_{\gamma,0;\zeta_1,\zeta_2}^J =\rho_{\gamma;\zeta_1,\zeta_2}^L.$
The identities in \eqref{eq:jaclag} prove
\eqref{eq:jacsup} when $\sigma=0$.

We now assume $\gamma\sigma>0$.  By Lemma \ref{remr},
$\mu_{\zeta_1,\zeta_2}^J$ is supported on a single interval $[s,t]$.
Put $m=(s+t)/2$, $h=(t-s)/2$, and
\begin{equation}\label{eq:jacdelta}
 \Delta_+=\sqrt{(1+\sqrt\gamma s)(1+\sqrt\gamma t)},
 \qquad
 \Delta_-=\sqrt{(1-\sigma\sqrt\gamma s)
                       (1-\sigma\sqrt\gamma t)}.
\end{equation}
The density formula in Lemma \ref{remr} gives
\begin{equation}\label{eq:jactri}
 \rho_{\gamma,\sigma;\zeta_1,\zeta_2}^J(x)=\frac{\pi-\ell(x)}
 {\pi^2\sqrt{(t-x)(x-s)}},
 \qquad
 \ell(x)=\pv\int_s^t\sqrt{(t-y)(y-s)}
 \frac{(V_{\gamma,\sigma}^J)'(y)}{x-y}\dd y.
\end{equation}
Since
$$
 (V_{\gamma,\sigma}^J)'(y)
 =-\frac{1-\gamma}{2\sqrt\gamma(1+\sqrt\gamma y)}
 +\frac{1-\sigma\gamma}
 {2\sqrt\gamma(1-\sigma\sqrt\gamma y)},
$$
applying Lemma \ref{le:Ia} with the shifts $-x$,
$\gamma^{-1/2}$, and $-(\sigma\sqrt\gamma)^{-1}$ yields
\begin{align*}
 J_0(x)
 &:=\pv\int_s^t\frac{\sqrt{(t-y)(y-s)}}{x-y}\dd y
 =\pi(x-m),\\
 J_+(x)
 &:=\pv\int_s^t\frac{\sqrt{(t-y)(y-s)}}
 {(1+\sqrt\gamma y)(x-y)}\dd y
 =\frac{\sqrt\gamma}{1+\sqrt\gamma x}
   \frac\pi\gamma(1+\sqrt\gamma m-\Delta_+)
   +\frac{J_0(x)}{1+\sqrt\gamma x} 
 =\frac\pi{\sqrt\gamma} -
   \frac{\pi \Delta_+}{\sqrt\gamma(1+\sqrt\gamma x)},\\
 J_-(x)
 &:=\pv\int_s^t\frac{\sqrt{(t-y)(y-s)}}
 {(1-\sigma\sqrt\gamma y)(x-y)}\dd y
 =-\frac{\sigma\sqrt\gamma}{1-\sigma\sqrt\gamma x}
   \frac\pi{\sigma^2\gamma}
   (1-\sigma\sqrt\gamma m-\Delta_-)
   +\frac{J_0(x)}{1-\sigma\sqrt\gamma x}\\
 &=-\frac\pi{\sigma\sqrt\gamma}
   \left(1-\frac{\Delta_-}{1-\sigma\sqrt\gamma x}\right).
\end{align*}
Consequently,
\begin{align*}
 \ell(x)
 =\frac{(\gamma-1)J_+(x)+(1-\sigma\gamma)J_-(x)}{2\sqrt\gamma}
 =\frac{(\gamma-1)\pi}{2\gamma}
   \bigg[1-\frac{\Delta_+}{1+\sqrt\gamma x}\bigg]
   -\frac{\pi(1-\sigma\gamma)}{2\sigma\gamma}
   \bigg[1-\frac{\Delta_-}{1-\sigma\sqrt\gamma x}\bigg].
\end{align*}
Substituting this into \eqref{eq:jactri} proves \eqref{eq:jacden}. 

It remains to determine the support endpoints $s$ and $t$ when
$\gamma\sigma>0$.  We first derive the equations for the free
endpoints.  Let $0<\gamma<1$, $0<\sigma\gamma\leq1$, and suppose that
$t=z$ while the left endpoint $a$ is free.  The endpoint condition
\eqref{eqa} gives
$$
 \frac1\pi\int_a^z
 \left[
 -\frac{1-\gamma}
 {2\sqrt\gamma(1+\sqrt\gamma x)}
 +\frac{1-\sigma\gamma}
 {2\sqrt\gamma(1-\sigma\sqrt\gamma x)}
 \right]
 \sqrt{\frac{z-x}{x-a}} \dd x=-1.
$$
By Lemma \ref{le:Ia},
\begin{align*}
 \int_a^z\frac{\sqrt{(z-x)/(x-a)}}{1-\sigma\sqrt\gamma x}
 \dd x
 &=
\frac1{1-\sigma\sqrt\gamma a}\int_a^z \sqrt{(z-x)(x-a)}\left(
 \frac1{x-a}
 +\frac{\sigma\sqrt\gamma}{1-\sigma\sqrt\gamma x}
 \right)\dd x \\
& =
 \frac\pi{\sigma\sqrt\gamma}
 \left(
 1-\sqrt{\frac{1-\sigma\sqrt\gamma z}
 {1-\sigma\sqrt\gamma a}}
 \right),
\end{align*}
and the same identity holds with $\si$ replaced by $-1$.
It follows that $\mathcal H_z(a)=1$, where $\mathcal H_z$ is
defined in \eqref{eq:jacH}. Similarly, let $0<\gamma\leq1$, $0<\sigma\gamma<1$, and suppose that
$s=z$ while the right endpoint $b$ is free.  The endpoint condition
\eqref{eqb} becomes
$$
 \frac1\pi\int_z^b
 \left[
 -\frac{1-\gamma}
 {2\sqrt\gamma(1+\sqrt\gamma x)}
 +\frac{1-\sigma\gamma}
 {2\sqrt\gamma(1-\sigma\sqrt\gamma x)}
 \right]
 \sqrt{\frac{x-z}{b-x}} \dd x=1.
$$
The same calculation using Lemma \ref{le:Ia} gives
$\mathcal H_z(b)=1$. The solutions are unique. Indeed, we have $\mathcal H_z(z)=0$,
$$
 \lim_{w\downarrow-\gamma^{-1/2}}\mathcal H_z(w)=+\infty
 \quad\text{if }\gamma<1,
 \qquad
 \lim_{w\uparrow(\sigma\sqrt\gamma)^{-1}}\mathcal H_z(w)=+\infty
 \quad\text{if }\sigma\gamma<1,
$$
and
$$
 \mathcal H_z''(w)
 =
 \frac{3(1-\gamma)\sqrt{1+\sqrt\gamma z}}
 {8(1+\sqrt\gamma w)^{5/2}}
 +
 \frac{3\sigma(1-\sigma\gamma)
 \sqrt{1-\sigma\sqrt\gamma z}}
 {8(1-\sigma\sqrt\gamma w)^{5/2}}>0.
$$
Thus $\mathcal H_z$ is strictly convex in the presence of a free endpoint.  Since $\mathcal H_z(z)=0<1$, the equation
$\mathcal H_z(w)=1$ has a unique solution in
$(-\gamma^{-1/2},z)$ when $\gamma<1$, and a unique solution in
$(z,(\sigma\sqrt\gamma)^{-1})$ when $\sigma\gamma<1$.

Recall that
$
 D_{\gamma,\sigma}
 =\bigl[-\gamma^{-1/2},(\sigma\sqrt\gamma)^{-1}\bigr],
 \widetilde u_\pm
 =\frac{\sqrt\gamma(1-\sigma)
 \pm2\sqrt{1+\sigma-\sigma\gamma}}{1+\sigma}.
$
We now determine which constraints are in effect.

\noindent\textit{Case
$\zeta_2<\widetilde u_+$:}
The upper wall is active, so $t=\zeta_2$.  If $\gamma=1$, then
$
 \widetilde u_-=-1=-\gamma^{-1/2},
$
and the lower natural edge coincides with the left hard edge of $D_{\gamma,\sigma}$.  Since $\zeta_1\geq-1$, it follows that
$[s,t]=[\zeta_1,\zeta_2]$.

If $0<\gamma<1$ and the lower wall is inactive, then $s$ is free and
$$
 \mathcal H_{\zeta_2}(s)=1,
 \qquad
 s=a_{\gamma,\sigma}(\zeta_2).
$$
The left endpoint $s=a_{\gamma,\sigma}(\zeta_2)$ satisfies the constraint $s\geq\zeta_1$ precisely when $\zeta_1\leq a_{\gamma,\sigma}(\zeta_2)$, in which case
$
 [s,t]=[a_{\gamma,\sigma}(\zeta_2),\zeta_2].
$
If $\zeta_1>a_{\gamma,\sigma}(\zeta_2)$, both walls are active and
$[s,t]=[\zeta_1,\zeta_2]$.

\noindent\textit{Case
$\zeta_1>\widetilde u_-$:}  The lower wall is active and $s=\zeta_1$.  If
$\sigma\gamma=1$, then
$
 \widetilde u_+=\sqrt\gamma
 =(\sigma\sqrt\gamma)^{-1},
$
so the unconstrained upper edge coincides with the right endpoint of $D_{\gamma,\sigma}$.  Since
$\zeta_2\leq(\sigma\sqrt\gamma)^{-1}$, we have
$[s,t]=[\zeta_1,\zeta_2]$.

If $0<\sigma\gamma<1$ and the upper wall is inactive, then $t$ is free
and
$$
 \mathcal H_{\zeta_1}(t)=1,
 \qquad
 t=b_{\gamma,\sigma}(\zeta_1).
$$
This endpoint is admissible if and only if
$\zeta_2\geq b_{\gamma,\sigma}(\zeta_1)$, and we have
$ [s,t]=[\zeta_1,b_{\gamma,\sigma}(\zeta_1)]$ in this case. If $\zeta_2<b_{\gamma,\sigma}(\zeta_1)$, both walls are active and
$[s,t]=[\zeta_1,\zeta_2]$.

We have shown \eqref{eq:jacsup} for $\gamma\sigma>0$ and the proof is complete.
\end{proof}
\begin{proof}[Proof of Theorem \ref{thm:jaclow}]
The proof of Example \ref{ex:JacempLDP} verifies the assumptions of Proposition \ref{pr:lowgap} for the rescaled density and both minimizers. Consequently,
\begin{align}\label{eq:jacen}
 I_{V_{\gamma,\sigma}^J}(\mu_{\zeta_1,\zeta_2}^J)
 &=\frac12\int V_{\gamma,\sigma}^J(x)
       \mu_{\zeta_1,\zeta_2}^J(\dd x)
   +\frac12F_{s,t}-C_{V_{\gamma,\sigma}^J},
\end{align}
where $F_{s,t}$ is determined by \eqref{eq:fzeta} with
$V=V_{\gamma,\sigma}^J$ and
$\mu_{\zeta_1,\zeta_2}=\mu_{\zeta_1,\zeta_2}^J$.

Using the support and density of
$\mu_{\zeta_1,\zeta_2}^J$ obtained in Proposition
\ref{p:jacmin}, we now compute its constrained rate
$I_{V_{\gamma,\sigma}^J}(\mu_{\zeta_1,\zeta_2}^J)$. We first consider $\gamma=0$ and $\sigma\geq0$.   By \eqref{eq:jcrdef},
\eqref{eq:v0sjh}, and \eqref{eq:herfun},
\begin{align*}
 I_{V_{0,\sigma}^J}(\mu_{\zeta_1,\zeta_2}^J)
 ={}&\int V_{0,\sigma}^J\dd\mu_{\zeta_1,\zeta_2}^J
 -\frac12\Sigma(\mu_{\zeta_1,\zeta_2}^J)-C_{V_{0,\sigma}^J}\\
 ={}&\int V^H\dd\mu_{c_\sigma\zeta_1,c_\sigma\zeta_2}^H
 -\frac12\Sigma(\mu_{c_\sigma\zeta_1,c_\sigma\zeta_2}^H)
 +\frac12\log c_\sigma-C_{V_{0,\sigma}^J}
 =I_{V^H}(\mu_{c_\sigma\zeta_1,c_\sigma\zeta_2}^H).
\end{align*}
Here $\mu^H_{\cdot}, V^H$ and $C_{V^H}$ are given in Example
\ref{ex:her2}; we have also   used the identity
$$
 C_{V_{0,\sigma}^J}
 =\frac38+\frac14\log(1+\sigma)
 =C_{V^H}+\frac12\log c_\sigma.
$$

Next assume that $\sigma=0$ and $\gamma\in(0,1]$.  By \eqref{eq:jcrdef} and
\eqref{eq:lagfun},
\begin{align*}
 I_{V_{\gamma,0}^J}(\mu_{\zeta_1,\zeta_2}^J)
 ={}&\int V_{\gamma,0}^J\dd\mu_{\zeta_1,\zeta_2}^J
 -\frac12\Sigma(\mu_{\zeta_1,\zeta_2}^J)-C_{V_{\gamma,0}^J}=\int V_\gamma^L\dd\mu_{\zeta_1,\zeta_2}^L
 -\frac12\Sigma(\mu_{\zeta_1,\zeta_2}^L)-C_{V_\gamma^L}
 =I_{V_\gamma^L}(\mu_{\zeta_1,\zeta_2}^L).
\end{align*}
Here $\mu_{\zeta_1,\zeta_2}^L, V_\gamma^L$ and $C_{V_\gamma^L}$ are given in Example
\ref{ex:Lag2}; we also used
$V_{\gamma,0}^J=V_\gamma^L$ on $D_\gamma$ and
$C_{V_{\gamma,0}^J}=C_{V_\gamma^L}$.  This proves the second case of
\eqref{eq:Jac_t}.

It remains to compute the rate for $\sigma\gamma>0$.  
Recall $\mathcal J_{s,t}$ from Lemma \ref{le:intlog} and $\Delta_\pm$ from \eqref{eq:jacdelta}.
Applying the lemma to the density \eqref{eq:jacden}, we obtain
\begin{align*}
\int_s^t\log(a+bx)\,\mu_{\zeta_1,\zeta_2}^J(\dd x)
&=\frac{1+\sigma}{2\sigma\gamma}
\mathcal J_{s,t}(a,b;1,0)
-\frac{(1-\gamma)\Delta_+}{2\gamma}
\mathcal J_{s,t}(a,b;1,\sqrt\gamma)\\
&\quad-\frac{(1-\sigma\gamma)\Delta_-}{2\sigma\gamma}
\mathcal J_{s,t}(a,b;1,-\sigma\sqrt\gamma).
\end{align*}
When $\gamma=1$ or $\sigma\gamma=1$, we omit the terms multiplied by $1-\gamma$ or $1-\sigma\gamma$, respectively, before evaluating the integrals and endpoint values below. 
For $(a,b)=(1,\sqrt\gamma)$, $(1,-\sigma\sqrt\gamma)$, and $(t,-1)$, we denote the corresponding integrals on the left-hand side by $M_+$, $M_-$, and $H_t$, respectively.
As a result, we find
\begin{align*}
 \int V_{\gamma,\sigma}^J\dd
 \mu_{\zeta_1,\zeta_2}^J
 &=-\frac{1-\gamma}{2\gamma}
   \int\log(1+\sqrt\gamma x)
       \mu_{\zeta_1,\zeta_2}^J(\dd x)-\frac{1-\sigma\gamma}{2\sigma\gamma}
   \int\log(1-\sigma\sqrt\gamma x)
       \mu_{\zeta_1,\zeta_2}^J(\dd x)\\
 &=-\frac{1-\gamma}{2\gamma}M_+
   -\frac{1-\sigma\gamma}{2\sigma\gamma}M_-.
\end{align*}
Taking $x=t$ in \eqref{eq:fzeta} gives
\begin{align*}
 F_{s,t}
 =V_{\gamma,\sigma}^J(t)-H_t
 =-\frac{1-\gamma}{2\gamma}\log(1+\sqrt\gamma t)
  -\frac{1-\sigma\gamma}{2\sigma\gamma}
   \log(1-\sigma\sqrt\gamma t)-H_t.
\end{align*}
Substituting these into \eqref{eq:jacen}, we obtain
\begin{align}\label{eq:jacgen}
 &I_{V_{\gamma,\sigma}^J}(\mu_{\zeta_1,\zeta_2}^J)
 =\Phi_{\gamma,\sigma}^J(s,t)=\frac12\int V_{\gamma,\sigma}^J
       \dd\mu_{\zeta_1,\zeta_2}^J
 +\frac12F_{s,t}-C_{V_{\gamma,\sigma}^J}\notag\\
 ={}&-\frac{1-\gamma}{4\gamma}
 \left[M_++\log(1+\sqrt\gamma t)\right]
 -\frac{1-\sigma\gamma}{4\sigma\gamma}
 \left[M_-+\log(1-\sigma\sqrt\gamma t)\right]
 -\frac12H_t-C_{V_{\gamma,\sigma}^J}.
\end{align}
This completes the proof.
\end{proof}

\appendix

\section{Proof of the largest eigenvalue LDP}\label{app:maxldp}
We prove Proposition \ref{pr:maxldp} using the following two lemmas.
\begin{lemma}
\label{le:maxub}
Under the assumptions of Proposition \ref{pr:maxldp},
\begin{equation}\label{eq:mxtigh}
 \lim_{R\to\infty}\limsup_{n\to\infty}\frac1{\beta n}
 \log\mathbb P(|\lambda_{\max}|>R)=-\infty.
\end{equation}
Moreover, for every $x\in\mathbb R$,
\begin{equation}\label{eq:mxlu}
 \lim_{\delta\downarrow0}\limsup_{n\to\infty}\frac1{\beta n}
 \log\mathbb P(|\lambda_{\max}-x|<\delta)\leq-J_V(x).
\end{equation}
\end{lemma}

\begin{lemma}
\label{le:maxlb}
Under the assumptions of Proposition \ref{pr:maxldp},
for every $x\in \R$,
\begin{equation}\label{eq:mxll}
 \lim_{\delta\downarrow0}\liminf_{n\to\infty}\frac1{\beta n}
 \log\mathbb P(|\lambda_{\max}-x|<\delta)\geq-J_V(x).
\end{equation}
\end{lemma}

\begin{proof}[Proof of Proposition \ref{pr:maxldp}]
By \ref{C1}, every $\mu\in\mathcal M_1(D)$ with $I_V(\mu)<\infty$ is supported on $\overline{\mathcal O}$.
Applying Lemma \ref{le:symreg} on $\overline{\mathcal O}$ verifies Assumption \ref{ass:d} for every $\mu\in\mathcal M_1(D)$ with $I_V(\mu)<\infty$.
Thus the empirical measure LDP in Theorem \ref{ldpln} applies\footnote{In fact, we only need LDP upper bound so that Assumption \ref{ass:d} for the unconstrained minimizer suffices.}. 
Assumptions \ref{A1}(ii), \ref{A2} and \ref{C1} ensure that $e^{-V}$ is a continuous admissible weight on $\overline{\mathcal O}$. 
By Remark \ref{rem:qepw}, \eqref{eq:vlogmu} holds pointwise on $\overline{\mathcal O}$ and hence on $D$, since $V=+\infty$ on $D\setminus\overline{\mathcal O}$.
The logarithmic potential is finite and continuous on $D$.
Consequently, $J_V\geq0$ and $J_V(b_V)=0$, while \ref{C1} gives $J_V(x)>0$ for $x>b_V$, $x\in D$. Moreover, $J_V$ is lower semicontinuous on $\mathbb R$. If $D$ is unbounded to the right, \ref{A2} gives $J_V(x)\to+\infty$ as $x\to+\infty$.
Therefore, $J_V$ is a good rate function.

Lemmas \ref{le:maxub} and \ref{le:maxlb} establish exponential tightness and the local upper and lower bounds for every $x\in \R$. The LDP therefore follows from \cite{AGZ10}*{Corollary D.6 and Theorem D.4(a)}.

Finally, for $x>b_V$ in the interior of $D \cap \{V<+\infty\} $, continuity of $J_V$ at $x$ gives
$$
\inf_{y\in D, y>x}J_V(y)=\inf_{y\in D, y\geq x}J_V(y).
$$
The LDP lower bound on $(x,\infty)$ and upper bound on $[x,\infty)$ therefore establish \eqref{eq:mxrt}.
\end{proof}

\begin{proof}[Proof of Lemma \ref{le:maxub}]
Using the empirical measure LDP in Theorem \ref{ldpln}, we follow the standard argument in \cite{AGZ10} to show that
\begin{equation}\label{eq:mxleft}
 \limsup_{n\to\infty}\frac1{\beta n}
 \log\mathbb P(\lambda_{\max}\leq x)=-\infty,\qquad x<b_V.
\end{equation}
Since
$ \mathbb P(|\lambda_{\max}|>R)  =\mathbb P(\lambda_{\max}>R)  +\mathbb P(\lambda_{\max}< -R),
$
Assumption \ref{C3} and equation \eqref{eq:mxleft} imply \eqref{eq:mxtigh}.

For $x<b_V$ and $0<\delta<b_V-x$,
$ \mathbb P(|\lambda_{\max}-x|<\delta)
 \leq\mathbb P(\lambda_{\max}\leq x+\delta),$
hence \eqref{eq:mxleft} gives \eqref{eq:mxlu}.
At $x=b_V$, \eqref{eq:mxlu} follows from $\mathbb P(|\lambda_{\max}-b_V|<\delta)\leq1$ and $J_V(b_V)=0$.
For $x\notin D$, the probability in \eqref{eq:mxlu} vanishes
for all sufficiently small $\delta>0$, since $D$ is closed.
It remains to prove \eqref{eq:mxlu} for $x>b_V$ in $D$. For every Borel set $A\subset D$,
\begin{equation}\label{eq:mxun}
\mathbb P(\lambda_{\max}\in A)
\leq\sum_{i=1}^n\mathbb P(\lambda_i\in A)
=n \mathbb P(\lambda_n\in A)
=n\int_{D^{n-1}}
\mathbb P(\lambda_n\in A\mid\by) \widetilde{\mathbb P}(\dd\by),
\end{equation}
where $\widetilde{\mathbb P}$ is marginal law of $(\lambda_1,\ldots,\lambda_{n-1})$. Fix $\varepsilon>0$.
By \ref{C1}, \eqref{eq:vlogmu}, and continuity of
$u\mapsto V(u)-\int\log|u-y| \mu_V(\dd y)$ at $b_V$,
choose $h>0$ such that $[b_V-h,b_V+h]$ lies in the interior of $D\cap \{V < +\infty\}$ and
\begin{equation}\label{eq:mxedge}
 V(u)-\int\log|u-y| \mu_V(\dd y)
 \leq F_V+\varepsilon,\qquad u\in[b_V-h,b_V+h].
\end{equation}
Let $\eta$ be the uniform probability measure on $[b_V-h,b_V+h]$. 
For $z,y\in D$, let
\begin{equation}\label{eq:mxgdef}
 g_z(y)=\log|z-y|-\int\log|u-y| \eta(\dd u).
\end{equation}
Conditioning on $(\lambda_1,\ldots,\lambda_{n-1})=\by$ in \eqref{eq:wpdf} gives
\begin{align}\label{eq:max-cond}
 \mathbb P(\lambda_n\in A\mid\by)
 &=\frac1{B_n(\by)}\int_A
 e^{-\beta nW_n(z)+\beta\sum_{i<n}g_z(y_i)}\omega_n(z)\dd z,
\end{align}
where
\begin{align*}
 B_n(\by)&=\int_D
 e^{-\beta nW_n(z)+\beta\sum_{i<n}g_z(y_i)}\omega_n(z)\dd z.
\end{align*}
The idea in this proof is to pin the logarithmic interaction at the upper edge of the equilibrium measure so that the regularized interaction is bounded above and no truncation is needed for $y_i$.
By restricting the integral
to $[b_V-h,b_V+h]$ and using $\dd x=2h \eta(\dd x)$, we obtain uniformly for $\by\in D^{n-1}$,
\begin{align}\label{eq:mxdl}
 B_n(\by) 
&\geq2h\int_{b_V-h}^{b_V+h}
e^{-\beta nW_n(x)+\log\omega_n(x)+\beta\sum_{i<n}g_x(y_i)}
\eta(\dd x)
\geq2h e^{-\beta n\int W_n\dd\eta
+\int\log\omega_n\dd\eta}\notag\\
&=e^{-\beta n\int V\dd\eta+o(\beta n)},
\end{align}
where in the second inequality we used Jensen's inequality and $\int g_x(y)\eta(\dd x)=0$, and in the last equality we used \ref{A0}, \ref{A1}(ii), and \ref{A3}(i).
Using \eqref{e:pairbd} and the continuity of $y\mapsto\int\log|u-y|\eta(\dd u)$, we obtain
\begin{equation}\label{eq:mxlog}
 g_x(y)=\log|x-y|-\int\log|u-y| \eta(\dd u)
 \leq\log(1+|x|)+C,\qquad x,y\in D,
\end{equation}
for some $C>0$ depending only on $\eta$. Moreover, for fixed $x$, $\lim_{y\to\8} g_x(y)=0$.

Now fix $x>b_V$ in $D$ and $0<\delta<x-b_V$. Let $K_{x,\delta}=[x-\delta,x+\delta]\cap D$,  $L_n=n^{-1}\sum_{i=1}^n\delta_{\lambda_i}$ and
$\widehat L_{n-1}=(n-1)^{-1}\sum_{i<n}\delta_{\lambda_i}$.
For $r>0$, let
$$A_{n-1}^r=\left\{\by=(y_1,\ldots,y_{n-1})\in D^{n-1}:
\frac1{n-1}\sum_{i<n}\delta_{y_i}\in B(\mu_V,r)\right\}.$$
Applying \eqref{eq:mxun} with $A=K_{x,\delta}$
and splitting the conditional integral yield
\begin{align}
 \mathbb P(|\lambda_{\max}-x|<\delta)
 &\leq n\mathbb P(\widehat L_{n-1}\notin B(\mu_V,r))+n\sup_{\by\in A_{n-1}^r}
 \mathbb P(\lambda_n\in K_{x,\delta}\mid\by). \label{eq:lamkxd}
\end{align}
Since $L_n=\frac{n-1}{n}\widehat L_{n-1} +\frac1n\delta_{\lambda_n}$, we have 
$d(L_n,\widehat L_{n-1}) \leq \frac2{n-1}$. The empirical measure LDP in Theorem \ref{ldpln}  gives
\begin{align}\label{eq:mxrest}
 \limsup_{n\to\infty}\frac1{\beta n}
 \log\mathbb P(\widehat L_{n-1}\notin B(\mu_V,r))
 \leq\limsup_{n\to\infty}\frac1{\beta n}
 \log\mathbb P\bigl(d(L_n,\mu_V)\geq r/2\bigr) =-\infty.
\end{align}

We claim that for sufficiently small $r>0$ and all sufficiently large $n$,
\begin{equation}\label{eq:mxkern}
\frac1n\sum_{i<n}g_z(y_i)
\leq\int g_z \dd\mu_V+\varepsilon,
\qquad \by\in A_{n-1}^r,\quad z\in K_{x,\delta}.
\end{equation}
Indeed, for fixed $M>0$, the functions
$
 g_{M,z}(y)=\max\{\log|z-y|,-M\}
 -\int\log|u-y| \eta(\dd u),  z\in K_{x,\delta},
$
are continuous and uniformly bounded over $z\in K_{x,\delta}$ and $y\in D$ by an argument similar to \eqref{eq:mxlog}. Moreover, $\lim_{|y|\to\infty}\sup_{z\in K_{x,\delta}}g_{M,z}(y)=0$,  $g_z\leq g_{M,z}$ and for every $z\in K_{x,\delta}$,
$$
 \int g_{M,z}\dd\mu_V\downarrow\int g_z\dd\mu_V
 \qquad\text{as }M\uparrow\infty.
$$
Both integrals are continuous functions of $z$. 
For the latter integral, this follows from the continuity of the logarithmic potential as stated in Remark \ref{rem:qepw}. Dini's theorem then allows us to fix $M$ so large that
\begin{equation}\label{eq:mxtrun}
 0\leq\int g_{M,z}\dd\mu_V-\int g_z\dd\mu_V
 \leq\frac{\varepsilon}{2},\qquad z\in K_{x,\delta}.
\end{equation}
Since $\|g_{M,z}-g_{M,z'}\|_\infty\leq e^M|z-z'|$,
compactness of $K_{x,\delta}$ gives finitely many points
$z_1,\ldots,z_m\in K_{x,\delta}$ such that every
$z\in K_{x,\delta}$ satisfies
$$
\|g_{M,z}-g_{M,z_j}\|_\infty<\frac{\varepsilon}{12}
$$
for some $j$. From the construction, each $g_{M,z_j}$ is bounded and continuous. Using the definition of $A_{n-1}^r$ with weak continuity, we may find $r>0$ such that
$$
\frac1{n-1}\sum_{i=1}^{n-1} g_{M,z_j}(y_i) 
\leq\int g_{M,z_j}\dd\mu_V+\frac{\varepsilon}{12},
\qquad \by\in A_{n-1}^r,\quad 1\leq j\leq m.
$$
By choosing the closest $z_j$ for each $z\in K_{x,\delta}$, we have uniformly for $\by\in A_{n-1}^r$ and $z\in K_{x,\delta}$,
\begin{align*}
\frac1{n-1}\sum_{i=1}^{n-1} g_{M,z}(y_i) 
\leq \frac1{n-1}\sum_{i=1}^{n-1} g_{M,z_j}(y_i) 
+\frac{\varepsilon}{12}
\leq\int g_{M,z_j}\dd\mu_V+\frac{\varepsilon}{6}
\leq\int g_{M,z}\dd\mu_V+\frac{\varepsilon}{4}.
\end{align*}
Since $g_z\leq g_{M,z}$ and
$\sup_{z\in K_{x,\delta}}\|g_{M,z}\|_\infty<\infty$,
this estimate and \eqref{eq:mxtrun} yield for all
sufficiently large $n$,
\begin{align*}
\frac1n\sum_{i<n}g_z(y_i)
&\leq \frac1{n}\sum_{i<n} g_{M,z}(y_i)
\leq\int g_{M,z}\dd\mu_V+\frac{\varepsilon}{2}
\leq\int g_z\dd\mu_V+\varepsilon, \quad \by\in A_{n-1}^r, \ \ z\in K_{x,\delta}.
\end{align*}

We claim that for every compact $K\subset D$ and continuous $\psi:K\to\mathbb R$,
\begin{equation}\label{eq:mxct}
 \limsup_{n\to\infty}\frac1{\beta n}
 \log\int_K e^{-\beta n[W_n(z)-\psi(z)]}\omega_n(z)\dd z
 \leq-\inf_{z\in K}[V(z)-\psi(z)].
\end{equation}
Indeed, for $\alpha>0$, cover $K$ by finitely many compact sets $K_1,\ldots,K_m\subset K$ with
$\sup_{K_j}\psi-\inf_{K_j}\psi\leq\alpha$.
Applying \eqref{eq:mxcomp} on these sets gives
\begin{align*}
 \limsup_{n\to\infty}\frac1{\beta n}
 \log\int_K e^{-\beta n[W_n(z)-\psi(z)]}\omega_n(z)\dd z
 \leq\max_{1\leq j\leq m}
 \left\{\sup_{K_j}\psi-\inf_{K_j}V\right\}
 \leq-\inf_K(V-\psi)+\alpha.
\end{align*}
Letting $\alpha\downarrow0$ proves \eqref{eq:mxct}.
Taking $K=K_{x,\delta}$ and $\psi(z)=\int g_z\dd\mu_V$, and using \eqref{eq:max-cond}, \eqref{eq:mxkern}, and \eqref{eq:mxdl}, we obtain
\begin{align}\label{eq:mxcond}
& \limsup_{n\to\infty}\frac1{\beta n}\log
 \sup_{\by\in A_{n-1}^r}
 \mathbb P(\lambda_n\in K_{x,\delta}\mid\by)
\leq \int V\dd\eta+\varepsilon
 -\inf_{z\in K_{x,\delta}}\left[V(z)-\int g_z\dd\mu_V\right]\notag\\
={}&-\inf_{K_{x,\delta}}J_V-F_V
 +\int\left[V(u)-\int\log|u-y| \mu_V(\dd y)\right]\eta(\dd u)
 +\varepsilon
\leq-\inf_{K_{x,\delta}}J_V+2\varepsilon.
\end{align}
Here the equality uses \eqref{eq:mxgdef} and \eqref{eq:mxrate}, and the last inequality follows from integrating \eqref{eq:mxedge} against $\eta$.
Combining \eqref{eq:mxcond}, \eqref{eq:lamkxd}, \eqref{eq:mxrest}, and \ref{C0}, we obtain
$$
 \limsup_{n\to\infty}\frac1{\beta n}
 \log\mathbb P(|\lambda_{\max}-x|<\delta)
 \leq-\inf_{K_{x,\delta}}J_V+2\varepsilon.
$$
Letting $\varepsilon\downarrow0$ and then $\delta\downarrow0$
proves \eqref{eq:mxlu}, since $J_V$ is lower semicontinuous.
\end{proof}

\begin{remark}
Under Assumption \ref{ass:ldp} with $D\subset\mathbb R$, for every $R>0$ and all sufficiently large $n$, we have the following estimate that will be used several times:
\begin{align}\label{eq:mxwtail}
 \mathbb P(\lambda_{\max}>R)
 \leq n\mathbb P(\lambda_n>R)
 \leq n e^{K\beta n+o(\beta n)}
 \int_{D\cap(R,\infty)}e^{-\beta nW_n(x)}
 \omega_n(x)(1+x)^{\beta(n-1)}\dd x,
\end{align}
where $K$ is independent of $n, R$. Indeed, choose a nondegenerate compact interval $I\subset\mathcal O$ and let $\eta$ be the uniform probability measure on $I$. Using Assumptions \ref{A0}, \ref{A1}(ii), and \ref{A3}(i), we can verify the estimates \eqref{eq:mxdl} and \eqref{eq:mxlog}.
Plugging into \eqref{eq:max-cond}, we find a conditional bound uniform in $\by$. Then the assertion follows from integrating over $\by$, the union bound and exchangeability.
\end{remark}

\begin{proof}[Proof of Lemma \ref{le:maxlb}]
The lower bound is immediate when $J_V(x)=+\infty$.
Using the local upper bound in Lemma \ref{le:maxub} and exponential tightness in \eqref{eq:mxtigh} together with a finite covering argument, we have
\begin{equation}\label{eq:mxtailub}
 \limsup_{n\to\infty}\frac1{\beta n}
 \log\mathbb P(\lambda_{\max}>T)
 \leq-\inf_{u\geq T}J_V(u)<0,\qquad T>b_V,
\end{equation}
where in the last inequality we used that $J_V$ is a good rate
function with unique zero $b_V$.
Combining this with \eqref{eq:mxleft}, we derive
$\lambda_{\max}\to b_V$ in probability, proving
\eqref{eq:mxll} at $b_V$.

Now fix $x>b_V$ in $\mathcal O$ and $\varepsilon>0$.
Choose $0<\delta<x-b_V$ with
$[x-\delta,x+\delta]\subset\mathcal O$, and let $\nu$ be the
uniform probability measure on $(x-\delta,x+\delta)$.
Recall $\eta$, $g_z$, $B_n$, and $A_{n-1}^r$ from the proof of
Lemma \ref{le:maxub}, with $r>0$ to be chosen sufficiently small.
Choose $R$ such that $\supp\mu_V\cup\supp\eta\subset(-R,R)$
and define
$$
 B_{n,R}(\by)=\int_{D\cap[-R,R]}
 e^{-\beta nW_n(u)+\beta\sum_{i<n}g_u(y_i)}\omega_n(u)\dd u,
 \qquad
 G_n=\{\by:B_{n,R}(\by)\geq B_n(\by)/2\}.
$$
For a continuous function $0\leq f\leq1$ that vanishes on
$\supp\mu_V$ and equals one outside $[-R,R]$, the empirical
measure LDP and exchangeability give
$$
 \mathbb P(|\lambda_n|>R)
 =\mathbb E L_n(D\setminus[-R,R])
 \leq\mathbb E\int f\dd L_n\to0.
$$
Since $B_{n,R}(\by)/B_n(\by)=\mathbb P(|\lambda_n|\leq R\mid\by)$ and $B_n(\by)\leq2B_{n,R}(\by)$ for $\by\in G_n$, Markov's inequality yields
\begin{equation}\label{eq:mxdt}
 \widetilde{\mathbb P}(G_n^c)=\widetilde{\mathbb P} \big(\mathbb P(|\lambda_n|> R\mid\by)>1/2\big)
 \leq2\mathbb P(|\lambda_n|>R)\to0,
 .
\end{equation}
For $E_n=A_{n-1}^r\cap G_n\cap\{\max_{i<n}y_i<x+\delta\}$,
\eqref{eq:mxrest}, \eqref{eq:mxdt}, and \eqref{eq:mxtailub} give
\begin{align*}
 \widetilde{\mathbb P}(E_n^c)
 &\leq\mathbb P(\widehat L_{n-1}\notin B(\mu_V,r))
 +\widetilde{\mathbb P}(G_n^c)
 +\mathbb P(\lambda_{\max}>x+\delta/2) \to 0.
\end{align*}

For sufficiently small $r$ and large $n$, uniformly on
$A_{n-1}^r\cap G_n$,
\begin{align}\label{eq:mxdu}
 \log B_n(\by)
 &\leq\log2+\log B_{n,R}(\by)
 \leq\log2+\beta n\varepsilon+
 \log\int_{D\cap[-R,R]}
 e^{-\beta n[W_n(u)-\int g_u\dd\mu_V]}\omega_n(u)\dd u\notag\\
 &\leq-\beta n\left[F_V+
 \iint\log|v-y|\eta(\dd v)\mu_V(\dd y)-\varepsilon\right]
 +o(\beta n),
\end{align}
where in the first two steps we used \eqref{eq:mxdt} and \eqref{eq:mxkern}, respectively and in the last step, we applied \eqref{eq:mxct} together with \eqref{eq:mxgdef} and
\eqref{eq:vlogmu}, since
$u\mapsto\int g_u\dd\mu_V$ is continuous.

By \eqref{eq:max-cond} and Jensen's inequality,
\begin{align}
 \log\mathbb P(x-\delta<\lambda_n<x+\delta\mid\by)
 \geq{}&\log(2\delta)-\log B_n(\by)
 -\beta n\int W_n\dd\nu+\int\log\omega_n\dd\nu
 +\beta\sum_{i<n}\int g_u(y_i)\nu(\dd u)\notag\\
 \geq{}&-\log B_n(\by)-\beta n\left[
 \int V\dd\nu-\iint g_u(y)\nu(\dd u)\mu_V(\dd y)
 +\varepsilon\right]+o(\beta n).
 \label{eq:mxjens}
\end{align}
Here the last inequality follows from \ref{A0}, \ref{A1}(ii), \ref{A3}(i), and weak continuity, since the function
$$
y\mapsto\int g_u(y)\nu(\dd u)
=\int\log|u-y| \nu(\dd u)
-\int\log|v-y| \eta(\dd v)$$
is bounded and continuous.
Combining \eqref{eq:mxdu} and \eqref{eq:mxjens}, and using
\eqref{eq:mxgdef} and \eqref{eq:mxrate}, we find
$$
 \frac1{\beta n}\log
 \mathbb P(x-\delta<\lambda_n<x+\delta\mid\by)
 \geq-\int J_V\dd\nu-2\varepsilon+o(1),
 \qquad \by\in E_n.
$$
On $E_n$, the event $x-\delta<\lambda_n<x+\delta$ implies
$x-\delta<\lambda_{\max}<x+\delta$. Hence
$$
 \mathbb P(|\lambda_{\max}-x|<\delta)
 \geq\int_{E_n}\mathbb P(|\lambda_n-x|<\delta\mid\by)
 \widetilde{\mathbb P}(\dd\by).
$$
Since $\widetilde{\mathbb P}(E_n)\to1$, we conclude that
$$
 \liminf_{n\to\infty}\frac1{\beta n}
 \log\mathbb P(|\lambda_{\max}-x|<\delta)
 \geq-\frac1{2\delta}\int_{x-\delta}^{x+\delta}J_V(u)\dd u
 -2\varepsilon.
$$
Letting $\delta\downarrow0$ and then $\varepsilon\downarrow0$, using the continuity of $J_V$, we have verified \eqref{eq:mxll} for $x>b_V$ in $\mathcal O$.

It remains to consider $e=\sup\mathcal O<+\infty$ when
$J_V(e)<+\infty$.
For $0<\delta<e-b_V$, let $y_\delta=e-\delta/2$.
Observing that
$(y_\delta-a,y_\delta+a)\subset(e-\delta,e+\delta)$ for
$0<a<\delta/2$, we derive 
$$
 \liminf_{n\to\infty}\frac1{\beta n}
 \log\mathbb P(|\lambda_{\max}-e|<\delta)
 \geq\lim_{a\downarrow0}\liminf_{n\to\infty}\frac1{\beta n}
 \log\mathbb P(|\lambda_{\max}-y_\delta|<a)
 \geq-J_V(y_\delta).
$$
By \ref{C1} and continuity of the logarithmic potential,
$J_V(y_\delta)\to J_V(e)$ as $\delta\downarrow0$.
This proves \eqref{eq:mxll} at $e$ and completes the proof.
\end{proof}

\section{Evaluation of some integrals}

\begin{lemma}\label{le:intloglem}
Let $A_2,B_2\in\mathbb R$ satisfy $A_2>|B_2|$, and set $\Delta_2=\sqrt{A_2^2-B_2^2}.$ Then
\begin{equation}\label{eq:coslog}
\int_{0}^{2\pi}\cos t\log(A_2+B_2\cos t)\dd t
=
2\pi\frac{B_2}{A_2+\Delta_2}.
\end{equation}
More generally, let $A_1,B_1\in\mathbb R$ satisfy $
A_1\ge |B_1|, 
A_1>0,
$
and set
$
\Delta_1=\sqrt{A_1^2-B_1^2}.
$
Then, with the integrals understood as improper integrals when
$A_1=|B_1|$, one has
\begin{equation}\label{eq:logaibi}
  \int_0^\pi
\frac{\log(A_1+B_1\cos\theta)}
{A_2+B_2\cos\theta} \dd\theta =
\frac{\pi}{\Delta_2}
\log\left[
\frac{A_1+\Delta_1}{2}
\left(
1-\frac{B_1B_2}
{(A_1+\Delta_1)(A_2+\Delta_2)}
\right)^2
\right].
\end{equation}
\end{lemma}

\begin{proof}
For $i=1,2$, whenever
$A_i>|B_i|$, set
$r_i=-\frac{B_i}{A_i+\Delta_i}.$
Then $|r_i|<1$ and
\begin{equation}\label{eq:daibi}
A_i+B_i\cos\theta
=
\frac{A_i+\Delta_i}{2}
\left(1-2r_i\cos\theta+r_i^2\right).
\end{equation} Since
$1-2r_i\cos\theta+r_i^2
=
(1-r_i e^{\ii\theta})(1-r_i e^{-\ii\theta}),$
the Taylor expansion of $\log(1-z)$ gives
\begin{equation}\label{eq:logfrier}
\log(1-2r_i\cos\theta+r_i^2)
=
-2\sum_{n=1}^{\infty}\frac{r_i^n}{n}\cos(n\theta),
\qquad |r_i|<1.
\end{equation}
Taking $i=2$, we get
$$
\int_0^{2\pi}\cos t\log(A_2+B_2\cos t)\dd t
=
-2r_2\int_0^{2\pi}\cos^2 t\dd t
=
-2\pi r_2= 2 \pi 
\frac{B_2}{A_2+\Delta_2}.
$$
Thus,  we obtain \eqref{eq:coslog}.

We next prove \eqref{eq:logaibi}. First assume that $A_1>|B_1|$. From \eqref{eq:daibi},
$$
\int_0^\pi
\frac{\log(A_1+B_1\cos\theta)}
{A_2+B_2\cos\theta}\dd\theta =
\frac{2}{A_2+\Delta_2}
\int_0^\pi
\frac{
\log\frac{A_1+\Delta_1}{2}
+
\log(1-2r_1\cos\theta+r_1^2)
}
{1-2r_2\cos\theta+r_2^2}\dd\theta.
$$
For $|r_2|<1$, the Taylor expansion gives
\begin{equation}\label{eq:pfourier}
\frac{1-r_2^2}{1-2r_2\cos\theta+r_2^2}
=
1+\frac{r_2e^{\ii \theta}}{1-r_2e^{\ii \theta}}
+\frac{r_2e^{-\ii \theta}}{1-r_2e^{-\ii \theta}}
=
1+2\sum_{n=1}^{\infty}r_2^n\cos(n\theta).
\end{equation}
Therefore,
\begin{equation}\label{eq:fcos}
\int_0^\pi
\frac{\dd\theta}
{1-2r_2\cos\theta+r_2^2}
=
\frac{\pi}{1-r_2^2}.
\end{equation}
Combining \eqref{eq:logfrier} and \eqref{eq:pfourier}, and
using orthogonality on $[0,\pi]$, we get
\begin{align*}
\int_0^\pi
\frac{\log(1-2r_1\cos\theta+r_1^2)}
{1-2r_2\cos\theta+r_2^2}\dd\theta 
=
-\frac{4}{1-r_2^2}
\sum_{n=1}^{\infty}
\frac{(r_1r_2)^n}{n}
\int_0^\pi \cos^2(n\theta)\dd\theta 
=
\frac{2\pi}{1-r_2^2}\log(1-r_1r_2).
\end{align*}
Thus
\begin{align*}
\int_0^\pi
\frac{\log(A_1+B_1\cos\theta)}
{A_2+B_2\cos\theta}\dd\theta  
=
\frac{\pi}{\Delta_2}
\log\left[
\frac{A_1+\Delta_1}{2}
\left(
1-\frac{B_1B_2}
{(A_1+\Delta_1)(A_2+\Delta_2)}
\right)^2
\right].
\end{align*}
This proves \eqref{eq:logaibi} when $A_1>|B_1|$. It remains to consider the case $A_1=|B_1|$ with $A_1>0$. For
$\varepsilon>0$, replace $A_1$ by $A_1+\varepsilon$ in the formula just
proved. Let
$\Delta_{1,\varepsilon}
=
\sqrt{(A_1+\varepsilon)^2-B_1^2}.$  Since $A_1+\varepsilon>|B_1|$, we find
$$
\int_0^\pi
\frac{\log(A_1+\varepsilon+B_1\cos\theta)}
{A_2+B_2\cos\theta}\dd\theta  =
\frac{\pi}{\Delta_2}
\log\left[
\frac{A_1+\varepsilon+\Delta_{1,\varepsilon}}{2}
\left(
1-\frac{B_1B_2}
{(A_1+\varepsilon+\Delta_{1,\varepsilon})(A_2+\Delta_2)}
\right)^2
\right].
$$
Sending $\varepsilon\downarrow0$, the right hand side converges to  \eqref{eq:logaibi}. For the left hand side, since $|A_2+B_2\cos\theta|>A_2-|B_2|$ and   $\log(A_1+B_1\cos\theta)$ is integrable on $(0,\pi)$, \eqref{eq:logaibi} follows from the dominated convergence theorem.
\end{proof}

\begin{lemma}\label{le:intlog}
Let $s<t$, $m=(s+t)/2$, and $h=(t-s)/2$.  Suppose that
$\alpha+\beta x\geq0$ on $[s,t]$ and is positive on $(s,t)$, and that
$\eta+\tau x>0$ on $[s,t]$.  Then
\begin{align}\label{eq:intlog}
 \mathcal J_{s,t}(\alpha,\beta;\eta,\tau)
 &:=
 \frac1\pi\int_s^t
 \frac{\log(\alpha+\beta x)}
 {(\eta+\tau x)\sqrt{(t-x)(x-s)}}\dd x \notag\\
 &=
 \frac1{\Delta_{\eta,\tau}}
 \log\left[
 \frac{\alpha+\beta m+\Delta_{\alpha,\beta}}2
 \left(
 1-\frac{\beta\tau h^2}
 {(\alpha+\beta m+\Delta_{\alpha,\beta})
  (\eta+\tau m+\Delta_{\eta,\tau})}
 \right)^2
 \right],
\end{align}
where
$\Delta_{\alpha,\beta}=\sqrt{(\alpha+\beta s)(\alpha+\beta t)}$
and
$\Delta_{\eta,\tau}=\sqrt{(\eta+\tau s)(\eta+\tau t)}$.
The integral is understood as an improper integral when
$\alpha+\beta s=0$ or $\alpha+\beta t=0$.
\end{lemma}

\begin{proof}
Under $x=m+h\cos\theta$, one has
$$
 \frac{\dd x}{\sqrt{(t-x)(x-s)}}=-\dd\theta,
$$
while
$ \alpha+\beta x=\alpha+\beta m+\beta h\cos\theta,
 \eta+\tau x=\eta+\tau m+\tau h\cos\theta.$
Thus \eqref{eq:intlog} follows directly from
Lemma \ref{le:intloglem}.  
\end{proof}

Recall that $\pv$ denotes the Cauchy principal value whenever the integrand has a simple pole inside the integration interval.
\begin{lemma}\label{le:Ia}
Let $\alpha<\beta$ and $x\in\mathbb R$. 
Then
$$
\pv\int_\alpha^\beta
\frac{\sqrt{(\beta-y)(y-\alpha)}}{y+x} \dd y
=
\begin{cases}
\pi\left(
x+\frac{\alpha+\beta}{2}
-
\sqrt{(x+\alpha)(x+\beta)}
\right),
& x>-\alpha,\\
\pi\left(
x+\frac{\alpha+\beta}{2}
\right),
& -\beta\le x\le -\alpha,\\
\pi\left(
x+\frac{\alpha+\beta}{2}
+
\sqrt{(x+\alpha)(x+\beta)}
\right),
& x<-\beta.
\end{cases}
$$
Here the principal value is needed  only when $-x\in(\alpha,\beta)$. 
\end{lemma}

\begin{proof}
Set $c=\frac{\alpha+\beta}{2}$ and $d=\frac{\beta-\alpha}{2}.$
By the change of variables $y=c+dt,$ we obtain
$$
J(x)
=\pv\int_\alpha^\beta
\frac{\sqrt{(\beta-y)(y-\alpha)}}{y+x} \dd y=
d \pv\int_{-1}^1
\frac{\sqrt{1-t^2}}{t-\eta} \dd t,
\qquad
\eta=-\frac{x+c}{d}.
$$
Using the Stieltjes transform of the semicircle law \cite{BAG97}, we have
\begin{equation}\label{eq:Gform}
G(z)=\int_{-1}^1 \frac{\sqrt{1-t^2}}{z-t}\dd t=\pi\Bigl(z-\sqrt{z^2-1}\Bigr),\quad z\in\C\setminus[-1,1],
\end{equation}
where the branch is chosen so that $\sqrt{z^2-1}\sim z$ as $z\to\infty$. By the Sokhotski--Plemelj formula, for $\alpha\in(-1,1)$, 
$$
\pv\int_{-1}^1\frac{\sqrt{1-t^2}}{t-\alpha}\dd t
=-\Re G(\alpha+i0)=-\pi\alpha.
$$
If $|\alpha|>1$, then \eqref{eq:Gform} gives
$$
\int_{-1}^1\frac{\sqrt{1-t^2}}{t-\alpha}\dd t
=-G(\alpha)
=-\pi\left(\alpha-\sgn(\alpha)\sqrt{\alpha^2-1}\right).
$$
Since
$$
d\eta=-(x+c),
\qquad
d\sqrt{\eta^2-1}
=
\sqrt{(x+c)^2-d^2}
=
\sqrt{(x+\alpha)(x+\beta)},
$$
and the three regimes $\eta<-1$, $-1\le\eta\le1$, and $\eta>1$ correspond respectively to $x>-\alpha,$ $-\beta\le x\le -\alpha,$ $x<-\beta,$ the stated formula follows.
\end{proof}

 We state the following elementary inequalities without proof.

\begin{lemma}The following logarithmic inequalities hold:
\begin{align}
 t-\log(1+t)&\geq\frac{t^2}{2(1+t)}, &t\geq0,\label{e:logpos}\\
 -t+\frac{t^2}{2}\leq-\log(1+t)&\leq-\frac{t}{1+t},
 &-1<t\leq0.\label{e:logneg}
\end{align}
For $z,w\in\mathbb C$,
\begin{equation}\label{e:pairbd}
 |z-w|\leq(1+|z|)(1+|w|), \qquad
 \log_+|z-w|\leq\log(1+|z|)+\log(1+|w|).
\end{equation}
Here $\log_+ x =\log x \mathbf1_{\{x\ge1\} }$. For every $\alpha,\varepsilon>0$, there is $C>0$ such that
\begin{align}
 \log(1+t)&\leq\varepsilon t^\alpha+C,
 &t\geq0,\label{eq:Frelb}\\
 \log|z-w|&\leq\varepsilon(|z|^\alpha+|w|^\alpha)+C,
 &z,w\in\mathbb C.\label{e:logpair}
\end{align}
\end{lemma}

\begin{lemma}
 There exists $C>0$ such that
\begin{align}
 \left|\log\Gamma(x)-\left(x-\frac12\right)\log x+x\right|
 &\leq C\bigl(1+(-\log x)_+\bigr), &x>0,\label{bd:stir2}\\
 -C\leq\log\Gamma(1+x)
 & \leq C(1+x)\log(2+x), &x\geq0.\label{e:gambd}
\end{align}
Here $(-\log x)_+ = -\log x \mathbf1_{\{x\le 1\} }$.
Moreover, for every integer $m\geq1$,
\begin{equation}\label{e:factbd}
 0\leq m\log m-\log(m!)\leq m-1.
\end{equation}
\end{lemma}

\bibliography{references}

\end{document}